%% file: 2005-39.tex
\documentclass{gtart_h}  

\input gtoutput

\lognumber{515}
\received{8 November 2004}
\volumenumber{9}\papernumber{39}\volumeyear{2005}
\pagenumbers{1689}{1774}   
\revised{1 September 2005}
\published{14 September 2005}
\accepted{18 August 2005}
\proposed{Benson Farb}
\seconded{Leonid Polterovich, David Gabai}

\usepackage{amsmath,amssymb,graphicx, color}

\def\fref#1{\hyperlink{#1anchor}{\ref*{#1}}}
\def\figref#1{\hyperlink{#1anchor}{Figure~\ref*{#1}}}
\def\anchor#1{\noindent\hypertarget{#1anchor}{\smash{$\phantom{99}$}}\newline}

\def\ie{ie }
\def\cf{cf }
\def\cad{c'est-\`a-dire}
\def\proofname{Preuve}

\def\adhe{\operatorname{Adh}}
\def\inte{\operatorname{Int}}
\newcommand{\bbR}{{\mathbb{R}}}
\newcommand{\bbZ}{{\mathbb{Z}}}
\def\sing{\mathrm{Sing}}
\def\wns{W_{N \rightarrow S}}
\def\wsn{W_{S \rightarrow N}}
\def\wss{W_{S \rightarrow S}}
\def\wnn{W_{N \rightarrow N}}
\def\wnf{W_{N \rightarrow }}
\def\wfs{W_{\rightarrow S}}
\def\wfn{W_{\rightarrow N}}
\def\wsf{W_{S \rightarrow}}

\def\pali{\operatorname{Palin}}

\def\ua{\uparrow}
\def\da{\downarrow}
\def\ra{\rightarrow}
\def\la{\leftarrow}

\def\t{\widetilde}
\def\enfs{E_{N \rightarrow S}}
\def\esfn{E_{S \rightarrow N}}

\newtheorem*{ques*}{Question}
\newtheorem*{prop*}{Proposition}
\newtheorem*{conj*}{Conjecture}
\newtheorem*{theo*}{Th\'eor\`eme}
\newtheorem{coro}{Corollaire}[section]
\newtheorem{theo2}[coro]{Th\'eor\`eme}
\newtheorem{affi}[coro]{Affirmation}
\newtheorem{prop}[coro]{Proposition}
\newtheorem{lemm}[coro]{Lemme}
\newtheorem{theoprin}{Th\'eor\`eme}
\newtheorem{propprin}[theoprin]{Proposition}
\newtheorem{propbis}{Proposition}
\newtheorem{theobis}{Th\'eor\`eme}
\newtheorem{theoter}{Th\'eor\`eme}
\theoremstyle{remark}
\newtheorem{defi}[coro]{D\'efinition}
\newtheorem{rema}[coro]{Remarque}

\begin{document}
\title{Structure des hom\'eomorphismes de Brouwer}
\covertitle{Structure des hom\noexpand\'eomorphismes de Brouwer}
\asciititle{Structure des homeomorphismes de Brouwer}
\author{Fr\'ed\'eric Le Roux}
\coverauthors{Fr\noexpand\'ed\noexpand\'eric Le Roux}
\asciiauthors{Frederic Le Roux}
\address{Universit\'e Paris Sud, Bat.\ 425,
91405 Orsay Cedex, FRANCE}
\asciiaddress{Universite Paris Sud, Bat. 425\\
91405 Orsay Cedex, FRANCE}
\email{frederic.le-roux@math.u-psud.fr.}


\begin{abstract}
For every Brouwer (ie planar, fixed point free, orientation preserving)
homeomorphism $h$ there exists a covering of the plane by
\emph{translation domains}, invariant simply-connected open subsets on
which $h$ is conjugate to an affine translation.  We introduce a
distance $d_h$ on the plane that counts the minimal number of
translation domains connecting a pair of points. This allows us to
describe a combinatorial conjugacy invariant, and to show the
existence of a finite family of generalised Reeb components separating
any two points $x,y$ such that $d_h(x,y) >1$.

\textbf{R\'esum\'e}

Tout hom\'eomorphisme de Brouwer s'obtient en recollant des \emph{domaines
de translation} (ouverts simplement connexes, invariants, en
restriction auxquels la dynamique est conjugu\'ee \`a une translation).
On introduit  une distance $d_h$ sur le plan qui  compte le nombre
minimal de domaines de translation dont la r\'eunion connecte deux points.
Ceci nous permet de d\'ecrire un invariant combinatoire de conjugaison,
qui d\'ecrit tr\`es grossi\`erement la mani\`ere dont les
domaines de translation se recollent. On montre \'egalement
  l'existence de structures dynamiques qui g\'en\'eralisent la
  pr\'esence de
 composantes de Reeb dans les feuilletages non triviaux du plan.
\end{abstract}

\asciiabstract{%
For every Brouwer (ie planar, fixed point free, orientation preserving)
homeomorphism h there exists a covering of the plane by translation
domains, invariant simply-connected open subsets on which h is
conjugate to an affine translation.  We introduce a distance d_h on
the plane that counts the minimal number of translation domains
connecting a pair of points.  This allows us to describe a
combinatorial conjugacy invariant, and to show the existence of a
finite family of generalised Reeb components separating any two points
x,y such that d_h(x,y)>1.

Resume

Tout homeomorphisme de Brouwer s'obtient en recollant des domaines de
translation (ouverts simplement connexes, invariants, en restriction
auxquels la dynamique est conjuguee a une translation).  On introduit
une distance d_h sur le plan qui compte le nombre minimal de domaines
de translation dont la reunion connecte deux points.  Ceci nous permet
de decrire un invariant combinatoire de conjugaison, qui decrit tres
grossierement la maniere dont les domaines de translation se
recollent. On montre egalement l'existence de structures dynamiques
qui generalisent la presence de composantes de Reeb dans les
feuilletages non triviaux du plan.}

\maketitlepage

\section{Introduction}
\subsection{Les feuilletages r\'eguliers du plan}\label{ss.feuilletages}
Consid\'erons, sur le plan,  un champ de vecteurs qui ne s'annule
pas. En l'int\'egrant, on obtient un feuilletage orient\'e $\cal F$
sans singularit\'e.
Ces objets sont parmi les premiers sur lesquels les dynamiciens se
soient pench\'es, et on peut dire qu'ils sont tr\`es bien
compris. Rappelons rapidement quelques r\'esultats. La th\'eorie de
Poincar\'e--Bendixson dit qu'il n'y a pas de r\'ecurrence: toute
orbite du champ de vecteurs part \`a l'infini. D'un point de vue un
peu plus
global, on montre que si l'on prend une carte du feuilletage $\cal F$,
l'union des orbites qui passe par cette carte est encore une carte de
$\cal F$. En d'autres termes, on peut recouvrir le plan par des ouverts
connexes, simplement connexes, invariants par la dynamique, sur
lesquels la restriction
 du feuilletage est triviale, c'est-\`a-dire
hom\'eomorphe au feuilletage horizontal du plan.

\begin{figure}[htp]\anchor{fig.Reeb}
\centerline{\hbox{\input{fig-Reeb2.pstex_t}}}
\caption{\label{fig.Reeb}Le feuilletage de Reeb}
\end{figure}
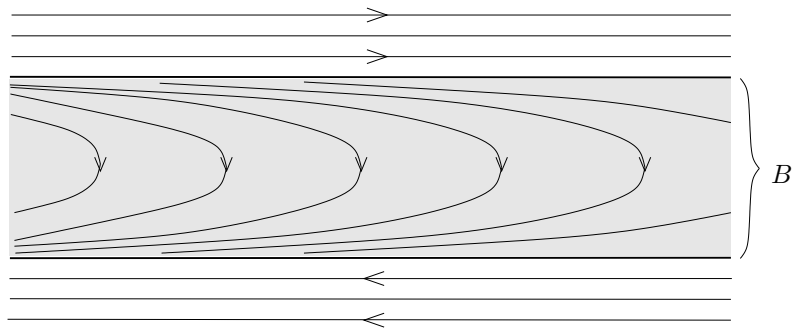

\begin{figure}[htp]\anchor{fig.feuilletage}
\centerline{\hbox{\input{fig-feuilletage4.pstex_t}}}
\caption{\label{fig.feuilletage} Tout couple de points est s\'epar\'e
par un nombre fini de composantes de Reeb}
\end{figure}
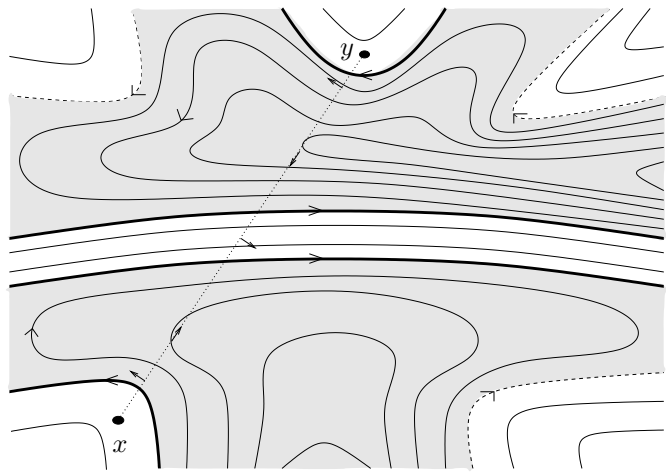

Ces r\'esultats de ``trivialit\'e semi-globale'' ont encourag\'e la
recherche d'une
classification compl\`ete de ces objets. W. Kaplan a ainsi d\'ecrit
les feuilletages topologiques r\'eguliers du plan au moyen d'objets un peu
exotiques appel\'es ``syst\`emes cordaux''\cite{kapl1,kapl2}. Plus
tard, G. Reeb, A. Haefliger et C. Godbillon ont \'etudi\'e
l'espace des feuilles, qui est une vari\'et\'e de dimension 1
\textit{non s\'epar\'ee}, le feuilletage s'obtenant comme fibr\'e
en droites au-dessus de cette vari\'et\'e \cite{haef1,gore,godb1}. On
peut extraire de ces
th\'eories les r\'esultats suivants. Le feuilletage non trivial le
plus simple est le \emph{feuilletage de Reeb}, repr\'esent\'e sur la
\figref{fig.Reeb}; il est unique \`a hom\'eomorphisme pr\`es.
Soit maintenant $\cal F$ un feuilletage r\'egulier quelconque.
Disons qu'un couple de feuilles $(\Delta_0,\Delta_1)$ de $\cal F$ est une
\emph{composante de Reeb}  s'il existe
un hom\'eomorphisme $\phi$ de la bande $B=\bbR \times [0,1] $ vers une
partie $\phi(B)$ du plan, qui  envoie $\{0\} \times \bbR$ sur $\Delta_0$
et
$\{1\} \times \bbR$ sur $\Delta_1$, et qui envoie le feuilletage
de Reeb de la bande $B$, repr\'esent\'e \figref{fig.Reeb}, sur la
restriction de $\cal F$ \`a $\phi(B)$
(remarquons que l'ensemble $\phi(B)$ n'est pas n\'ecessairement ferm\'e,
voir les exemples de la  \figref{fig.feuilletage}). Les feuilles
$\Delta_0$ et
$\Delta_1$ sont appel\'ees \emph{bords} de la composante de Reeb.
La notion de composante de Reeb est fondamentale pour comprendre les
feuilletages du
plan: en effet, tout feuilletage non trivial admet des composantes de
Reeb. On peut m\^eme pr\'eciser ce  r\'esultat de la mani\`ere suivante.
\begin{theo*}
Soit $\cal F$ un feuilletage du plan sans singularit\'e.
On se donne deux points $x$ et $y$, et on note $d$ le nombre minimal de
cartes du feuilletage dont la r\'eunion est connexe et contient $x$
et $y$. Alors il existe exactement $d-1$ composantes de Reeb
dont chacun des deux bords s\'epare (faiblement)\footnote{Un ensemble
 \emph{s\'epare (faiblement)} deux points si ceux-ci ne sont
pas dans la m\^eme composante connexe de son compl\'ementaire. Ici,
les points $x$ ou $y$ peuvent \'eventuellement appartenir au bord d'une
composante de Reeb; cependant, cette situation est non-g\'en\'erique,
au sens o\`u la r\'eunion des bords de toutes les composantes de Reeb
du feuilletage est un ensemble \emph{maigre}  (r\'eunion d\'enombrable
de ferm\'es d'int\'erieurs vides).
Par ailleurs, notons que cet \'enonc\'e n'appara\^it pas explicitement
dans les articles cit\'es ci-dessus.}
$x$ et $y$.
\end{theo*}

Le th\'eor\`eme est illustr\'e par la \figref{fig.feuilletage},
avec $d=3$.
En particulier, ce th\'eor\`eme permet de d\'efinir un
invariant combinatoire qui d\'ecrit totalement le
feuilletage ``entre $x$ et $y$''. Par exemple, la partie du
feuilletage entre $x$ et $y$ sur la \figref{fig.feuilletage}
pourrait \^etre repr\'esent\'ee par le mot $(\uparrow, \rightarrow,
\downarrow,
\leftarrow, \uparrow)$, obtenu de la fa\c con suivante (voir la figure):
on se d\'eplace sur le segment joignant $x$ \`a $y$, et on y \'ecrit
les fl\`eches au fur et \`a mesure;
les fl\`eches verticales
codent la dynamique sur les bords des composantes de Reeb, les
fl\`eches horizontales codent la dynamique \`a l'int\'erieur d'une
composante
(ce proc\'ed\'e est pr\'ecis\'e plus bas).

\subsection{Les hom\'eomorphismes de Brouwer}
Un \emph{hom\'eomorphisme de Brouwer} est un hom\'eomorphisme du
plan, pr\'eservant l'orientation, et sans point fixe. Ces objets ont
 \'et\'e introduits par Brouwer, puis \'etudi\'es notamment par
 Ker\'ekj\'art\'o,
Homma--Terasaka, Andrea \cite{brou1,kere2,homm1,andr1,andr3}.
Dans une certaine mesure, on peut les voir  comme une
g\'en\'eralisation des champs de vecteurs qui ne s'annulent pas: en
effet, d'apr\`es le th\'eor\`eme de Poincar\'e--Bendixson, tout
temps $\Phi_t$ du flot obtenu en int\'egrant un tel champ de vecteurs est
sans point fixe, c'est donc un hom\'eomorphisme de Brouwer. En fait,
le th\'eor\`eme de Poincar\'e--Bendixson se g\'en\'eralise aux
hom\'eomorphismes de Brouwer: l\`a aussi, il y a absence de
r\'ecurrence, et toute orbite part \`a l'infini. L'analogie culmine
avec le c\'el\`ebre th\'eor\`eme des translations planes de
Brouwer: on peut recouvrir le plan par des \emph{domaines de
translation}, qui sont des ouverts, connexes, simplement connexes,
invariants par $h$, sur lesquels la restriction de $h$ est triviale,
c'est-\`a-dire
conjugu\'ee \`a une translation affine du plan. En raccourci, tout
hom\'eomorphisme de Brouwer peut se construire en recollant des
translations.

Peut-on pousser plus loin cette analogie avec les
feuilletages? La litt\'erature regorge d'exemples
d'hom\'eomorphismes de Brouwer exotiques, qui semblent tr\`es
diff\'erents de ceux provenant d'un champ de vecteurs
(voir~\cite{kere2,homm1,brow2,daw1,lero8})
M\^eme si on se limite aux dynamiques obtenues en recollant deux
translations, il faut abandonner l'id\'ee d'obtenir une
classification compl\`ete (voir~\cite{naka1}, \cite{lero8}, et
\cite{begu1} o\`u l'on
tentait de d\'ecrire le plus finement possible ce recollement).
Le point de vue de cet article est en quelque sorte compl\'ementaire
de celui de \cite{begu1}. Il s'agit ici de d\'ecrire de mani\`ere
vague le recollement d'un grand nombre de translations, en oubliant
toutes les diff\'erences qu'on analysait l\`a en d\'etail. Avec ce
point de vue, nous allons g\'en\'eraliser le th\'eor\`eme d'existence de
composantes de Reeb \'enonc\'e ci-dessus pour les feuilletages, et
en d\'eduire la construction d'un invariant
de conjugaison combinatoire qui \'etend le codage esquiss\'e plus haut.

\subsection{Description des r\'esultats}
Nous d\'ecrivons maintenant les r\'esultats de l'article  de
mani\`ere plus pr\'ecise.
\paragraph{Distance de translation}
Soit $h$ un hom\'eomorphisme de Brouwer. \'Etant donn\'es deux
points distincts $x$ et $y$, on note $d_h(x,y)$ le nombre minimal
 de domaines de translation de $h$ dont la r\'eunion est connexe et
contient $x$ et $y$. Ceci d\'efinit une distance sur le plan.
On d\'efinit \'egalement la $h$--longueur d'un arc, c'est le nombre
minimal de morceaux dans un d\'ecoupage en sous-arcs \emph{libres}
(un ensemble $E$ est \emph{libre} si $h(E) \cap E = \emptyset$).
On montre alors que la distance $d_h(x,y)$ co\"\i ncide avec la plus
courte longueur des arcs joignant $x$ \`a $y$. Les arcs r\'ealisant
le minimum sont nomm\'es arcs \emph{g\'eod\'esiques}.
Dans le cas o\`u l'hom\'eomorphisme $h$ est obtenu en int\'egrant un champ
de vecteurs, le nombre $d_{h}(x,y)$ co\"incide avec le nombre $d$ qui
appara\^it dans le th\'eor\`eme de la Section~\ref{ss.feuilletages}
(nombre minimum de cartes du feuilletage dont la r\'eunion est connexe et
contient $x$ et $y$). Sur l'exemple de la \figref{fig.feuilletage}, le
segment $[xy]$ est un arc g\'eod\'esique de $h$--longueur \'egale \`a $3$.

\paragraph{Composantes de Reeb}
D'autre part, nous g\'en\'eralisons la notion de composante de Reeb
de la mani\`ere suivante. Essentiellement, on dit qu'un couple
$(F,G)$ de parties disjointes, ferm\'ees et
connexes du plan est une \emph{composante de Reeb} pour un couple de
points $(x,y)$ si
\begin{itemize}
\item $F$ s\'epare $x$ et $G \cup \{y\}$ (ou bien $F$ contient $x$),
\item $G$ s\'epare $y$ et $F \cup \{x\}$ (ou bien $G$ contient $y$),
\item le produit $F \times G$ est \emph{positivement} ou bien
\emph{n\'egativement singulier}, \textit{i.e.}  pour tout couple $(V,W)$
d'ouverts
rencontrant respectivement $F$ et $G$, il existe des temps $n$
\emph{positifs}
arbitrairement grands tels que $h^n(V)$ rencontre $W$, ou bien
pour tout couple $(V,W)$ d'ouverts
rencontrant respectivement $F$ et $G$, il existe des temps $n$
\emph{n\'egatifs}
arbitrairement grands tels que $h^n(V)$ rencontre $W$.
\end{itemize}
 On peut v\'erifier que cette
d\'efinition est compatible avec celle donn\'ee dans le cadre des
feuilletages. On prouve alors le th\'eor\`eme suivant, qui g\'en\'eralise
l'\'enonc\'e de la Section~\ref{ss.feuilletages}.
\begin{theoprin}
\label{theo.crgm}
Soient $x$ et $y$ deux points du plan, et $d=d_h(x,y)$.
Alors il existe exactement $d-1$  composantes de Reeb
pour $(x,y)$ qui sont minimales pour l'inclusion.
\end{theoprin}
De plus, ces composantes sont donn\'ees par une formule utilisant
la distance $d_h$ et la topologie du plan: en particulier, pour chaque
 composante $(F,G)$, le bord $G$ est situ\'e sur la
fronti\`ere d'une boule pour la distance $d_h$, centr\'ee en $x$.

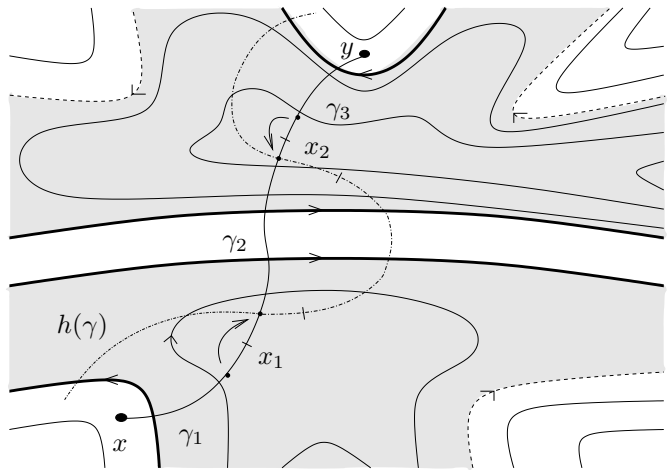
\begin{figure}[htp]\anchor{fig.feuilletage-geo}
\centerline{\hbox{\input{fig-feuilletage-geo2.pstex_t}}}
\caption{\label{fig.feuilletage-geo}
Le  mot $(\ra \la)$ obtenu ne d\'epend que de $(x,y)$ }
\end{figure}

\paragraph{Fl\`eches horizontales}
Comme pour les feuilletages, ceci va  permettre de d\'efinir un invariant
de conjugaison combinatoire, d\'ependant de $x$ et de $y$,
et qui d\'ecrit (partiellement cette fois-ci) la dynamique
de $h$ ``entre $x$ et $y$''. Cet invariant sera, l\`a encore, un mot sur
l'alphabet $\{\leftarrow, \rightarrow,\uparrow,\downarrow\}$, dont la
longueur est li\'ee \`a la distance
$d_h(x,y)$, et dans lequel les fl\`eches horizontales et verticales
alternent.
 La d\'efinition se fait en deux temps.
Voici la d\'efinition des fl\`eches horizontales.
Notons $d$ la distance $d_h(x,y)$.
On consid\`ere une \emph{d\'ecomposition minimale} d'un arc
g\'eod\'esique $\gamma$ joignant $x$ \`a $y$, c'est-\`a-dire un
d\'ecoupage en $d$ sous-arcs libres (voir la
\figref{fig.feuilletage-geo}, o\`u l'arc $\gamma$ est d\'ecoup\'e
par les points $x_{1}$ et $x_{2}$). Consid\'erons deux sous-arcs
adjacents $\gamma_i$ et $\gamma_{i+1}$ dans la d\'ecomposition.
Par minimalit\'e, leur r\'eunion n'est pas libre; on en d\'eduit
que $h(\gamma_i)$ rencontre $\gamma_{i+1}$ ou bien que
 $h(\gamma_{i+1})$ rencontre $\gamma_{i}$. Le point cl\'e est que ces
deux possibilit\'es ne peuvent pas avoir lieu simultan\'ement: en
effet, ceci est une cons\'equence de l'absence de
``quasi-orbite'' p\'eriodique pour les hom\'eomorphismes de Brouwer
(l'\'enonc\'e pr\'ecis est un lemme de J. Franks, rappel\'e \`a la
Section~\ref{ssec.rappels}). On peut donc
associer \`a cette d\'ecomposition de l'arc g\'eod\'esique
$\gamma$ un mot $(m_1 \dots m_{d-1})$ sur l'alphabet $\{\rightarrow,
\leftarrow\}$, avec  $m_i=\rightarrow$ si
$h(\gamma_i)$ rencontre $\gamma_{i+1}$ et $m_i=\leftarrow$ si c'est le
contraire. Il reste \`a prouver l'\'enonc\'e suivant.
\begin{theoprin}
\label{theo.fleches-hori}
Le mot ainsi d\'efini ne d\'epend que de $x$ et de $y$, et pas du
choix de l'arc g\'eod\'esique $\gamma$ joignant $x$ \`a $y$ (ni de
sa d\'ecomposition).
\end{theoprin}

\paragraph{Fl\`eches verticales}
La d\'efinition des fl\`eches verticales est plus compliqu\'ee,
 nous nous restreignons au cas des \emph{suites g\'eod\'esiques
 infinies}\footnote{Dans le cas des suites g\'eod\'esiques finies,
 il y a des
 complications suppl\'ementaires dues \`a des probl\`emes de bord,
 et il n'est pas clair que les fl\`eches initiales et finales soient
 bien d\'efinies.},
c'est-\`a-dire des suites de points $(x_k)_{k \in \bbZ}$ pour
lesquelles $d_h(x_k,x_l)=|l-k|$. Le Th\'eor\`eme~\ref{theo.crgm} permet
 d'associer \`a la suite g\'eod\'esique $(x_k)$ une suite de composantes
 de Reeb minimales $(F_k,G_k)_{k \in \bbZ}$ (voir la partie sup\'erieure
 de la \figref{fig.releve}). Et le
 Th\'eor\`eme~\ref{theo.fleches-hori} associe \`a $(x_k)$
 une suite infinie de fl\`eches horizontales, li\'ees \`a la
 dynamique \`a l'int\'erieur de chaque composante. Nous allons
 intercaler une fl\`eche verticale entre chaque couple de fl\`eches
 horizontales successives, et cette fl\`eche sera li\'ee cette
 fois-ci \`a la dynamique sur les bords des composantes; ceci
 donnera un mot infini $M$ sur l'alphabet de quatre fl\`eches. Nous
 proc\'edons de la mani\`ere suivante.
On peut compactifier chaque bord $F_k$ (ou $G_k$)
en lui ajoutant deux points  \`a l'infini, les points $N$ (Nord)  et
 $S$ (Sud), de mani\`ere \`a ce que l'orbite de tout point de $F_k$
 tende vers $N$
 ou vers $S$ (cf \figref{fig.releve}). Dans ce cadre,
 nous montrons deux
 propri\'et\'es: d'une part,  la dynamique
dans  $F_k$ (ou dans $G_k$) est ``\`a sens unique''; d'autre part,
la dynamique
 ``circule dans le m\^eme sens'' dans les bords adjacents $G_k$ et
 $F_{k+1}$:
\begin{propprin}
\label{propprin.continus-minimaux}
De deux choses l'une: ou bien il existe des orbites de
$F_k$ allant de $N$ vers $S$, ou bien il existe des orbites de $F_k$
allant de $S$ vers $N$.
\end{propprin}
\begin{theoprin}
\label{theo.fleches-verticales}
Le bord $G_k$ contient des orbites allant de $N$ vers $S$ si et
seulement si le bord $F_{k+1}$ en contient \'egalement.
\end{theoprin}
 Ces deux propri\'et\'es permettent de
 d\'efinir les fl\`eches verticales: la fl\`eche correspondant \`a $G_{k}$
 et \`a $F_{k+1}$ vaut $\uparrow$ si la dynamique dans ces deux bords se
 fait du Sud vers le Nord, et $\downarrow$ si c'est le contraire.

La partie sup\'erieure de la \figref{fig.releve} illustre ces
notions. Le point $N$ est ``\`a l'infini vers le haut'', le point $S$
``\`a l'infini vers le bas''.
Comme l'exemple dessin\'e vient d'un champ de vecteurs, les bords des
composantes de Reeb sont des droites, et les fl\`eches verticales sont
alternativement $\ua$ et $\da$; ces propri\'et\'es ne seront pas toujours
vrai pour un hom\'eomorphisme qui n'est pas issu d'un champ de vecteurs
(des exemples sont donn\'es dans l'Appendice~\ref{sec.exemples}).

\subsection{Application: un indice pour les points fixes}
\begin{figure}[htp]\anchor{fig.releve}
\centerline{\hbox{\input{fig-releve6.pstex_t}}}
\caption{\label{fig.releve}
 Un hom\'eomorphisme du plan $h$ ayant un unique point fixe, et son
 relev\'e canonique $\t h$; le mot infini p\'eriodique $M$ est \'ecrit
 au-dessus de la droite g\'eod\'esique  $\Gamma$;  le mot cyclique $M(h)$
 vaut $(\da \la \ua \la \da \la \ua \ra)$}
\end{figure}
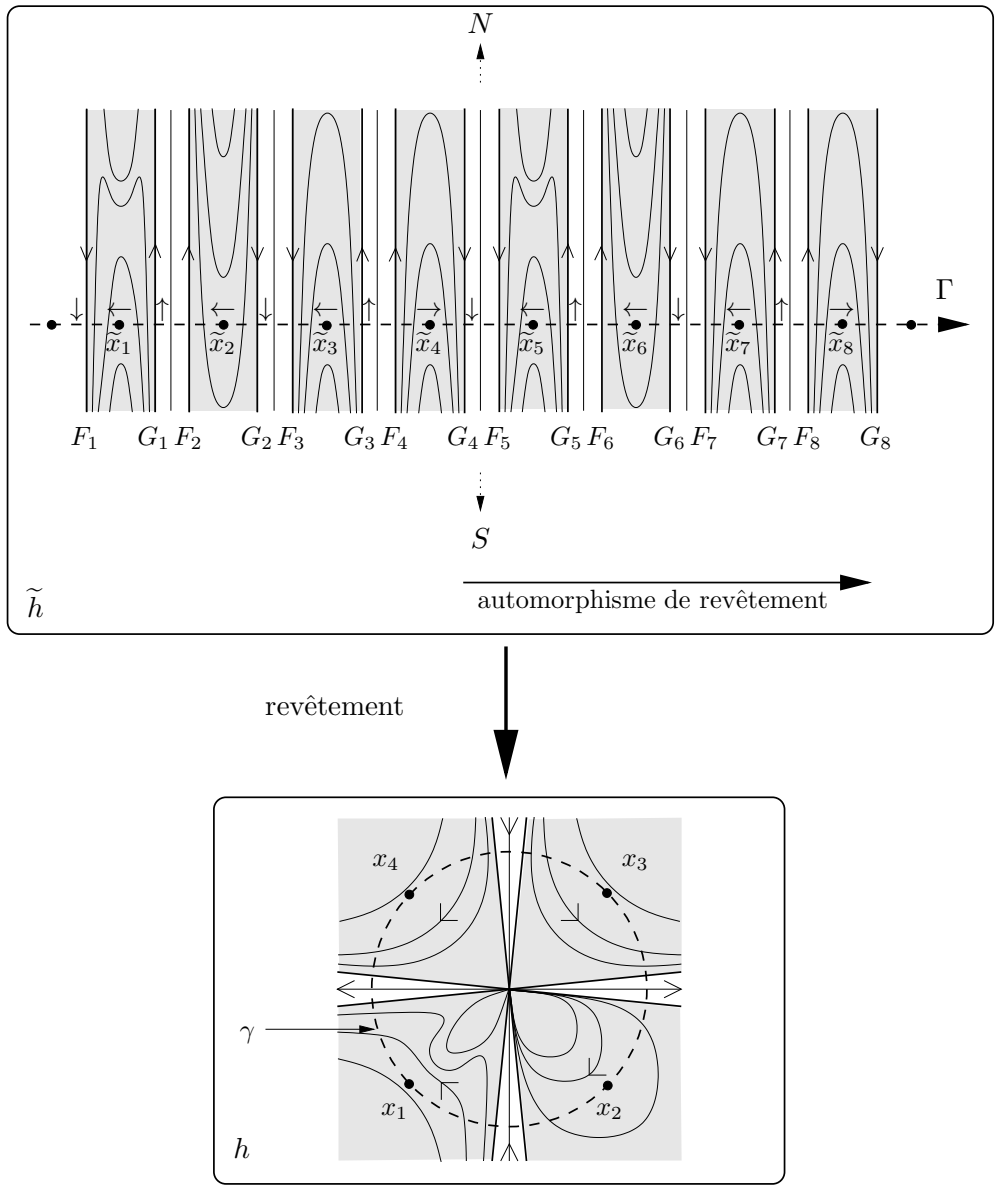

Toutes ces constructions sont motiv\'ees par l'\'etude de la
dynamique des hom\'eo\-mor\-phismes de surfaces autour d'un point
fixe isol\'e, et sont appliqu\'ees \`a cette \'etude dans
l'article~\cite{leroinf}. \'Evoquons rapidement les r\'esultats
correspondants.
 On se donne un hom\'eomorphisme $h$ (pr\'eservant l'orientation)
d'une surface (orientable)
$S$ et un point fixe isol\'e $x_0$. Nous associons \`a ces
donn\'ees un indice qui affine l'indice de Poincar\'e--Lefschetz, et ce
nouvel indice est un mot cyclique sur l'alphabet
$\{\leftarrow, \rightarrow,\uparrow,\downarrow\}$. R\'esumons la
construction dans le cas particulier o\`u $h$ est un hom\'eomorphisme du
plan dont $x_{0}$ est l'unique point fixe (voir l'exemple repr\'esent\'e
sur la partie inf\'erieure de la \figref{fig.releve}).
Le rev\^etement universel de $\bbR^2 \setminus \{0\}$
est un plan, et les relev\'es de $h$ \`a ce plan sont clairement des
hom\'eomorphismes sans
point fixe: ce sont des  hom\'eomorphismes de Brouwer. Un point essentiel
est que si l'indice de $x_0$ est diff\'erent
de $1$, parmi tous les relev\'es de $h$, il en existe un unique qui
ne soit pas conjugu\'e \`a une translation.
On se place sous cette hypoth\`ese d'indice, et on consid\`ere ce
relev\'e canonique $\widetilde h$ (repr\'esent\'e sur la partie
sup\'erieure de la \figref{fig.releve}).
On peut  alors trouver  une  courbe de Jordan $\gamma$, qui entoure
 $x_0$, et dont l'un des relev\'es $\Gamma$ est  une \emph{droite
 g\'eod\'esique}
pour $\tilde h$ (\ie\ $\Gamma$ contient une suite
 g\'eod\'esique infinie $(\t x_{k})$ comme d\'efinie ci-dessus).

Soit maintenant $M$ le mot infini sur l'alphabet de quatre fl\`eches
qui est
associ\'e \`a $\Gamma$. Alors $M$ ne d\'epend pas du choix de la
courbe $\gamma$, et pr\'esente une p\'eriodicit\'e de longueur
paire $2d$. On obtient ainsi un mot cyclique $M(h)$, de longueur $2d$,
qui ne
d\'epend que des donn\'ees $S$, $x_0$ et $h$.
Nous expliquons aussi dans ~\cite{leroinf} comment retrouver l'indice de
Poincar\'e--Lefschetz \`a partir de ce mot: ``on compte le
nombre de tours faits par la fl\`eche lorsqu'on lit le mot $M(h)$''.
Et on en d\'eduira quelques cons\'equences dynamiques. Par exemples, on
pr\'ecisera, dans ce cadre, la version topologique du  th\'eor\`eme de la
fleur de Leau--Fatou donn\'ee dans \cite{lero2001}.

\subsection{Plan du texte}
La Section~\ref{sec.preliminaires} propose des d\'efinitions et des
rappels.
La Section~\ref{sec.distance} raconte la d\'efinition de la distance de
translation $d_{h}$, des fl\`eches horizontales, et donne l'\'enonc\'e
correspondant au  Th\'eor\`eme~\ref{theo.fleches-hori} d'invariance.
La Section~\ref{sec.crg} expose la d\'efinition des composantes de Reeb,
pr\'ecise l'\'enonc\'e du  Th\'eor\`eme~\ref{theo.crgm} d'ex\-ist\-ence,
et en d\'eduit le Th\'eor\`eme~\ref{theo.fleches-hori}.
Les Sections~\ref{sec.topologie-disque} et~\ref{sec.preuve-crg}
fournissent la preuve du Th\'eor\`eme~\ref{theo.crgm}.
La Section~\ref{sec.geodesiques-bout} introduit les suites g\'eod\'esiques
infinies et leurs bouts, et leur associe une topologie sur le plan
augment\'e des points $N$ et $S$.
La Section~\ref{sec.fleches-verticales} prouve le
Th\'eor\`eme~\ref{theo.fleches-verticales} qui conduit \`a la d\'efinition
des fl\`eches verticales.
La Section~\ref{sec.continus-minimaux} contient la preuve de la
Proposition~\ref{propprin.continus-minimaux} d\'ecrivant la dynamique
dans les bords des composantes de Reeb. Le premier appendice rappelle
quelques \'enonc\'es de topologie, le second illustre le texte par trois
exemples, qui font voir en quoi la vie est plus compliqu\'ee pour les
hom\'eomorphismes de Brouwer qui ne sont pas issus d'un champ de vecteurs.

\subsection{Remerciements}
Ce texte a b\'en\'efici\'e de discussions avec François B\'eguin,
Christian Bonatti, Sylvain Crovisier, John Franks,
Lucien Guillou, Vincent Guirardel et Patrice Le Calvez.

\setcounter{theo}{0}
\section{Pr\'eliminaires}
\label{sec.preliminaires}
\subsection{D\'efinitions, notations}
La droite r\'eelle est not\'ee $\bbR$, et $\bbZ$ d\'esigne l'anneau
des entiers relatifs.
\subsubsection*{Topologie \ldots{}}
On notera respectivement $\inte(E)$, $\adhe(E)$, $\partial E$
l'int\'erieur,
l'adh\'erence, la fronti\`ere
d'une partie $E$ dans un espace topologique. On utilisera souvent
l'inclusion g\'en\'erale $\partial \adhe(E) \subset \partial E$.

La notion de s\'eparation joue un r\^ole essentiel dans ce texte:
\begin{defi}
Soient $X$ un espace topologique (qui sera le plus souvent le plan
$\bbR^2$) et  $F$,  $A$ et $B$ trois parties
de $X$. On dira que \emph{$F$ s\'epare $A$ et $B$} si
\begin{itemize}
\item $A$ est incluse dans une composante connexe $O_A(F)$ du
compl\'ementaire
de $F$ dans $X$;
\item $B$ est incluse dans une composante connexe $O_B(F)$ du
compl\'ementaire
de $F$ dans $X$;
\item $O_A(F) \neq O_B(F)$.
\end{itemize}
\end{defi}

Un \emph{arc} dans un espace topologique $X$ est une application
continue injective de $[0,1]$ dans $X$. Puisque le param\'etrage de
l'arc importe peu en g\'en\'eral, on identifiera $\gamma$ \`a
$\gamma \circ \phi$ o\`u $\phi$ est un hom\'eomorphisme croissant de
l'intervalle $[0,1]$; on
confondra m\^eme souvent un arc $\gamma$
et son image $\gamma([0,1])$. On dira que $\gamma$ \emph{joint}
$\gamma(0)$ \`a $\gamma(1)$, ou que ces deux points sont les
\emph{extr\'emit\'es} de $\gamma$.
Un \emph{sous-arc} de $\gamma$ est un arc inclus dans $\gamma$.
Si $x$ et $y$ sont deux points sur un arc $\gamma$, on notera
$[xy]_\gamma$ le sous-arc de $\gamma$ joignant $x$ \`a $y$; on notera
\'egalement $]xy]_\gamma$ ce sous-arc priv\'e du point $x$,
\textit{etc.}. L'\emph{int\'erieur} de l'arc $\gamma$ est
$]\gamma(0), \gamma(1)[_\gamma=\gamma(]0,1[)=\gamma
\setminus\{\gamma(0),\gamma(1)\}$.

Soient $\gamma_1$ et $\gamma_2$ deux arcs tels que
$\gamma_1(1)=\gamma_2(0)$. On appellera \emph{concat\'enation de
$\gamma_1$ et $\gamma_2$}, et on notera $\gamma_1 * \gamma_2$, la
courbe d\'efinie par $\gamma_1 * \gamma_2 (t) = \gamma_1(2t)$ pour $t
\in [0,1/2]$ et  $\gamma_1 * \gamma_2 (t) = \gamma_2(2t-1)$ pour $t
\in [1/2,1]$. Cette loi induit une loi sur les classes d'\'equivalence
modulo reparam\'etrages croissants, loi qui est associative.

Une \emph{droite topologique} (du plan) est une application $\Gamma$
continue, injective,
propre\footnote{La propret\'e revient \`a
dire que $\Gamma$ se prolonge contin\^ument aux compactifi\'es
d'Alexandroff en
envoyant l'infini sur l'infini, \ie en un plongement du cercle $\bbR\cup
\{\infty\}$ dans la sph\`ere $\bbR^2 \cup \{\infty\}$.},
de $\bbR$ dans $\bbR^2$.
Une \emph{droite topologique orient\'ee} est une classe
d'\'equivalence de  droites topologiques modulo reparam\'etrage croissant.

Un \emph{disque topologique ferm\'e} est une partie du plan
hom\'eomorphe au disque unit\'e ferm\'e du plan.

Nous utiliserons souvent diff\'erentes variantes du th\'eor\`eme de
Schoenflies (\cf\ \cite{cair51}); celles-ci affirment en particulier
qu'un arc, une droite
topologique ou un disque topologique ferm\'e sont, \`a
hom\'eomorphisme du plan pr\`es, respectivement un segment euclidien, une
droite euclidienne et un disque euclidien. Ainsi, une droite topologique
orient\'ee $\Gamma$ s\'epare le plan en deux ouverts, l'un situ\'e \`a
gauche de $\Gamma$ et l'autre \`a sa droite.

\subsubsection*{\ldots et dynamique}
Soit $h$ un hom\'eomorphisme d'un espace topologique $X$. Un
\emph{point fixe} de $h$ est un point $x$ tel que $h(x)=x$.
L'\emph{orbite} de $x$ est la suite $(h^n(x))_{n \in \bbZ}$, son
\emph{orbite positive} est $(h^n(x))_{n \geq 0}$.
Une partie $E$ de $X$ est \emph{invariante} si $h(E)=E$.

Un hom\'eomorphisme du plan \emph{pr\'eserve l'orientation} s'il
envoie l'ouvert \`a gauche d'une droite topologique orient\'ee
$\Gamma$ sur l'ouvert \`a gauche de la droite top\-o\-log\-ique orient\'ee
$h \circ \Gamma$. Dans ce cas, on peut montrer que cette
propri\'et\'e est vraie pour toute droite topologique orient\'ee.

\subsection{Rappels sur les hom\'eomorphismes de
Brouwer}\label{ssec.rappels}
Dans tout le texte, sauf mention contraire, $h$ d\'esigne
un hom\'eomorphisme de Brouwer, \ie un hom\'eomorphisme du plan
$\bbR^2$, sans point fixe, pr\'eservant l'orientation.
On rappelle ici quelques r\'esultats sur le comportement dynamique d'un
hom\'eomorphisme de Brouwer. Dans la section suivante, nous traduirons
ces r\'esultats en associant \`a $h$ une distance sur le plan.
Les r\'ef\'erences pour cette partie (et notamment pour les
d\'emonstrations) sont \cite{brou1,guil1,fran88,fran92,leca3,lero4}.

\subsubsection*{Domaines de translation}
\begin{defi}
Un \emph{domaine de translation} est une partie $O$ du plan, ouverte,
connexe et simplement connexe, invariante par $h$, sur laquelle la
restriction de $h$ est
conjugu\'ee \`a la translation du plan $\tau\co (x,y) \mapsto (x+1,y)$
(\cad\
qu'il existe un hom\'eomorphisme $\Phi$ de $O$ sur $\bbR^2$ v\'erifiant
$\Phi h \Phi^{-1}=\tau$).
\end{defi}

Le r\'esultat suivant est un morceau du c\'el\`ebre
``th\'eor\`eme des translations pl\-anes'' de Brouwer; c'est
un r\'esultat fondamental concernant la structure des
hom\-\'e\-o\-morph\-ismes de Brouwer.
\begin{theo*}[Brouwer]
Tout point du plan est dans un domaine de translation.
\end{theo*}

\begin{defi}
Un sous-ensemble $E$ du plan est \emph{libre}  si $h(E) \cap
E = \emptyset$. Il est \emph{essentiellement libre} si il existe un
entier $n$ tel que  $h^n(E) \cap E = \emptyset$.
\end{defi}
Remarquons que tout point du plan appartient \`a  un disque
 libre (il suffit de prendre un disque  assez petit). Les
 preuves du th\'eor\`eme des translations planes montrent en
 r\'ealit\'e que \emph{tout compact connexe libre est inclus dans un
 domaine de
 translation}. On peut encore g\'en\'eraliser ce r\'esultat:
\begin{lemm}[\cite{lero4} compacts connexes libres]
\label{lemm.ETDS}
 On consid\`ere  une par\-tie $K$ \!com\-pacte et connexe du plan. Alors
 $K$ est incluse dans un domaine de
 translation si et seulement si elle est essentiellement libre.
\end{lemm}

\subsubsection*{Le lemme de Franks}

L'absence de r\'ecurrence est r\'esum\'ee par le ``lemme de Franks'',
qui joue un r\^ole central dans ce texte. La plupart du temps, nous nous
contenterons de la version \`a deux disques:

\emph{\'Etant donn\'es $D_1$ et $D_2$ deux disques topologiques
ferm\'es, d'int\'erieurs disjoints,
d'int\'erieurs libres, on ne peut pas
avoir simultan\'ement, pour deux entiers \emph{positifs} $n_1$ et $n_2$,
$$
h^{n_1}(\inte(D_1)) \cap \inte(D_2) \neq \emptyset \mbox{ et }
\inte(D_1) \cap h^{n_2}(\inte(D_2)) \neq \emptyset.
$$}
Ce premier \'enonc\'e admet une reformulation imm\'ediate.
\begin{lemm}{\rm\cite{lero4}}\qua
\label{lemm.intervalle}
Soit $D$ et $V$ deux disques topologiques ferm\'es d'int\'e\-rieurs
libres. Alors l'ensemble
des entiers $n$ tels que $h^n(\inte(D))$ rencontre $\inte(V)$ est un
intervalle de $\bbZ$.
\end{lemm}

Nous appliquerons une seule fois\footnote{\cf preuve du
Lemme~\ref{lemm.limitesup}.}
 une version plus g\'en\'erale, qui
utilise la notion suivante.
\begin{defi}
 Une \emph{cha\^\i ne de disques p\'eriodique}
 est une suite
  $(D_i)_{i=1..k}$, ($k \geq 2$),  de disques topologiques ferm\'es
 tels que:
\begin{enumerate}
  \item les disques ouverts $\inte(D_i)$ sont disjoints deux \`a deux;
  \item chaque disque ouvert $\inte(D_i)$ est libre;
  \item il existe des entiers strictement positifs $n_1, \dots, n_k$ tels
  que
\begin{eqnarray*}
h^{n_1}(\inte(D_1)) & \cap & \inte(D_2) \neq \emptyset \\
& \vdots & \\
h^{n_{k-1}}(\inte(D_{k-1})) & \cap & \inte(D_k) \neq \emptyset \\
h^{n_k}(\inte(D_k)) & \cap & \inte(D_1) \neq \emptyset.
\end{eqnarray*}
\end{enumerate}
\end{defi}
\begin{lemm}[Franks \cite{fran88}]\label{lemm.franks}
Il n'existe pas de cha\^\i ne de disques p\'eriodique.
\end{lemm}
Ce lemme a d\'ej\`a de nombreuses applications (voir
\cite{fran88,fran89,leca3,lero2001}).
Pr\'ecisons que l'\'enonc\'e que nous utilisons ici est contenu
dans l'\'enonc\'e d'origine de \cite{fran88}; en particulier, nous n'avons
pas besoin des sophistications du lemme 3.18 de \cite{lero2001}.

En appliquant le lemme pr\'ec\'edent \`a la suite $(\inte(D),\inte(h(D)))$
o\`u
$D$ est un disque topologique ferm\'e libre, on obtient le corollaire
suivant.
\begin{coro}\label{coro.liberte}
Un disque topologique ferm\'e libre $D$  est disjoint de l'int\'erieur
de tous ses it\'er\'es $h^n(D)$ pour $n \neq 0$.
\end{coro}
\`A l'aide du th\'eor\`eme de Schoenflies, on peut voir que tout arc
libre est inclus dans l'int\'erieur d'un disque libre: par cons\'equent,
le Corollaire~\ref{coro.liberte} entra\^\i ne que \emph{tout arc libre
est disjoint de tous ses it\'er\'es}.
Une autre cons\'equence est:
\begin{coro}\label{coro.erre}
Tout  point $x$ du plan erre sous l'action de $h$, et l'orbite
positive de $x$  tend vers l'infini; autrement dit
\begin{itemize}
\item il existe un voisinage $V$ de $x$ tel que pour tout $n \neq 0$,
$h^n(V) \cap V = \emptyset$;
\item pour toute partie  compacte $K$ du plan, il existe un entier $n_0$
tel que pour tout $n \geq n_0$, $h^n(x) \not \in K$.
\end{itemize}
\end{coro}

Rappelons qu'en g\'en\'eral, contrairement \`a l'intuition, les
it\'er\'es du voisinage $V$, eux, ne sortent pas du compact $K$. Autrement
dit: la suite des it\'er\'es de $h$ converge simplement vers
l'infini mais non uniform\'ement. Cette remarque conduit \`a la
notion classique de \emph{point singulier}, que nous rappellerons
\`a la Section~\ref{ssec.singulier}. Le comportement de la suite des
itérés d'un disque libre sera étudié à la
Section~\ref{sec.topologie-disque} (et est
représenté sur la \figref{fig.iteres-disque}).

\section{Distance de translation, fl\`eches horizontales}
\label{sec.distance}
Les propri\'et\'es des hom\'eomorphismes de Brouwer, rappel\'ees
dans la partie pr\'ec\'edente, vont nous permettre d'associer \`a
l'hom\'eomorphisme de Brouwer $h$ une distance $d_h$ sur $\bbR^2$; cette
distance peut se comprendre en termes de domaines de translation ou
bien en termes d'arcs libres (Section~\ref{ssec.distance}). \`A l'aide
de cette distance, nous construisons ensuite un invariant de
conjugaison associ\'e \`a la donn\'ee d'un couple de points
(Section~\ref{ssec.fleches}, Th\'eor\`eme~\ref{theo.fleches-hori}).

\subsection{Distance de translation de $h$}
\label{ssec.distance}

\subsubsection*{Domaines de translation et distance}
\begin{defi}
On appelle \emph{$h$--distance} ou \emph{distance de translation de
$h$}, et on note $d_h$, la distance d\'efinie
sur le plan de la mani\`ere suivante.
Si $x$ et $y$ sont deux points distincts du plan, $d_h(x,y)$ est le
plus petit
entier $k$ tel qu'il existe des domaines de translation $O_1, \dots
O_k$ dont l'union est connexe et contient $x$ et $y$.
\end{defi}
Remarquons que l'existence de l'entier $k$ est imm\'ediate: d'apr\`es le
th\'eor\`eme des translations planes de Brouwer,
on peut  recouvrir le
segment euclidien $[xy]$ par un nombre fini de domaines de
translation. Le fait que l'application $d_h$ soit une distance d\'ecoule
facilement de la d\'efinition.

\begin{defi}\label{defi.suite-geo}
Une \emph{suite g\'eod\'esique} (joignant $x_0$ \`a $x_d$) est une suite
de points
$(x_0,\dots, x_d)$
du plan telle que $d_h(x_0,x_d)=d$ et $d_h(x_i,x_{i+1})=1$ pour tout
$i=0, \dots, d-1$.
\end{defi}
Il est clair que pour tous points $x,y$, il
existe une suite g\'eod\'esique joignant $x$ \`a $y$. D'autre part,
l'in\'egalit\'e triangulaire entra\^ine que si $(x_0,\dots, x_d)$ est une
suite g\'eod\'esique, on a $d_h(x_i,x_j)=j-i$ pour tous $0 \leq i \leq j
\leq d$.

\subsubsection*{Distance et arcs libres}

\begin{defi}[\rm(\figref{fig.decomposition})]
\label{defi.decomposition}
Soit $\gamma$ un arc du plan. Une \emph{d\'ecomposition} de $\gamma$ est
une \'ecriture
$$
\gamma=\gamma_1 * \dots * \gamma_k
$$ o\`u les sous-arcs $\gamma_i$ sont libres. Les \emph{sommets} de la
d\'ecomposition sont les extr\'emit\'es  $x_0, \dots x_k$ des arcs
$\gamma_i$, comme sur la \figref{fig.decomposition}.
L'entier $k$ est appel\'e
\emph{longueur} de la d\'ecomposition. La \emph{$h$--longueur} de
$\gamma$ est le plus petit entier $k$ tel qu'il existe une
d\'ecomposition de $\gamma$ de longueur $k$. Une d\'ecomposition
r\'ealisant
ce minimum est appel\'ee \emph{d\'ecomposition minimale} de $\gamma$.\footnote{On obtient une d\'efinition alternative   en rempla\c{c}ant
``libre'' par ``essentiellement libre'' dans cette d\'efinition. Ceci
conduirait \`a la m\^eme th\'eorie.}
\end{defi}
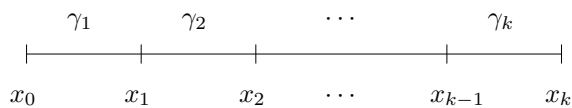
\begin{figure}[htp]\anchor{fig.decomposition}
\centerline{\hbox{\input{fig-decomposition.pstex_t}}}
\caption{\label{fig.decomposition}D\'ecomposition d'un arc $\gamma$}
\end{figure}

L'affirmation suivante, dont la preuve est donn\'ee plus bas,  fait le
lien entre arcs libres et domaines de
translation. En corollaire, elle permet de donner une d\'efinition de
la distance $d_h$ en termes d'arcs.
\begin{affi}[Domaines de translation et arcs
libres]\label{affi.arcs-et-domaines}
Tout arc libre est inclus dans un domaine de
translation. R\'eciproquement, si
 deux points $x$ et $y$ appartiennent \`a un m\^eme domaine de
translation $O$,
alors ou bien $x$ et $y$ sont dans la m\^eme orbite de $h$, ou bien il
existe un arc libre joignant $x$ \`a $y$ et inclus dans $O$.
\end{affi}
\begin{coro}\label{coro.definition-arcs}
Soient $x$ et $y$ deux points du plan qui ne sont pas dans la m\^eme
orbite. Alors le nombre
$d_h(x,y)$ est le plus petit entier $k$ tel qu'il existe un arc
d'extr\'emit\'es $x$ et $y$ et de $h$--longueur $k$.
\end{coro}
Remarquons que deux points distincts dans la m\^eme orbite sont
clairement \`a
$h$--distance $1$.
\begin{defi}
Un arc r\'ealisant le minimum du Corollaire~\ref{coro.definition-arcs}
est appel\'e arc
\emph{g\'eod\'esique}.
\end{defi}
En particulier, la suite des sommets d'une d\'ecomposition minimale
d'un arc g\'eod\'esique est une suite g\'eod\'esique (au sens de la
D\'efinition~\ref{defi.suite-geo}). Dans la suite du texte, toutes les
d\'ecompositions consid\'er\'ees seront des d\'ecompositions minimales
d'arcs
g\'eod\'esiques. Un exemple d'arc g\'eod\'esique de longueur $3$ est
visible dans l'introduction (\figref{fig.feuilletage-geo}).

\begin{proof}[\proofname\ de l'Affirmation~\ref{affi.arcs-et-domaines}]
La premi\`ere partie est une cons\'equence imm\-\'e\-di\-ate du
Lemme~\ref{lemm.ETDS} sur les ensembles compacts connexes libres
 (et m\^eme de sa version simple, \cf les
commentaires le pr\'ec\'edant).

Soient $x$ et $y$ deux points quelconques du plan. Si $x$ et $y$ ne
sont pas dans la m\^eme orbite de la translation horizontale $\tau\co
(x,y) \mapsto (x+1,y)$,
alors il existe un arc joignant $x$ \`a $y$ et libre pour $\tau$:
en effet, tout segment euclidien non
horizontal est libre pour $\tau$, et on peut traiter le cas o\`u $x$
et $y$ sont sur la m\^eme horizontale en raisonnant dans l'anneau
quotient $\bbR^2/\tau$ (ou bien en conjuguant pour \'eviter cette
situation).
Par conjugaison, on en d\'eduit la deuxi\`eme partie de l'affirmation.
\end{proof}

\begin{proof}[\proofname\ du Corollaire~\ref{coro.definition-arcs}]
Soient $x$ et $y$ deux points du plan qui ne sont pas dans la m\^eme
orbite, et $p=d_h(x,y)$. Nous allons construire un arc reliant $x$ \`a
$y$ et de $h$--longueur inf\'erieure ou \'egale \`a $p$.
Soient $O_1, \dots, O_p$ des domaines de translation dont l'union est
connexe et contient $x$ et $y$. Quitte \`a r\'eindexer les ouverts $O_i$,
on peut supposer que $O_1$ contient $x$, que $O_p$ contient $y$, et
que $O_i$ rencontre $O_{i+1}$ ($i=1, \dots, p-1$). Posons $x_0=x$,
$x_p=y$, et choisissons pour chaque $i$ un point $x_i$ dans $O_i \cap
O_{i+1}$. Quitte \`a bouger un peu les points  $x_1, \dots, x_{p-1}$
on peut supposer de
plus que les orbites des points $x_0, \dots, x_p$ sont deux \`a deux
distinctes (on utilise ici que $x$ et $y$ ne sont pas dans la m\^eme
orbite). On applique alors la deuxi\`eme partie de
l'Affirmation~\ref{affi.arcs-et-domaines}:
pour chaque $i=1 \dots p$, il existe un arc libre $\gamma_{i}$ reliant
$x_{i-1}$ \`a $x_{i}$ et  inclus dans $O_i$.
Remarquons que si $|i-j| \geq 2$, alors les ouverts $O_i$ et $O_j$
sont disjoints:
ceci est d\^u \`a la d\'efinition de $p=d_h(x,y)$. Par
cons\'equent, les arcs  $\gamma_i$ et $\gamma_j$ sont \'egalement
disjoints.
 Soit $\gamma= \gamma_1 * \cdots * \gamma_{p}$. Si
$\gamma$ est un arc, on a trouv\'e un arc de $h$--longueur inf\'erieure ou
\'egale \`a $p$ reliant $x$ et $y$.
Dans le cas contraire, on remplace
$\gamma$ par un arc $\gamma'$ en suivant la proc\'edure sugg\'er\'ee
par la
\figref{fig.simplification}:
on parcourt $\gamma_1$ en partant de $x_0$ jusqu'au
premier point de $\gamma_2$ qu'on appelle $x'_1$, et on remplace
l'arc $\gamma_1$ par le sous-arc $\gamma'_1=[x_0 x'_1]_{\gamma_1}$;
puis on parcourt $\gamma_2$ en partant de $x'_1$ jusqu'au
premier point de $\gamma_3$ qu'on appelle $x'_2$, et on remplace
l'arc $\gamma_2$ par le sous-arc $\gamma'_2=[x'_1 x'_2]_{\gamma_2}$;
\textit{etc.}.
On v\'erifie enfin que $\gamma'= \gamma'_1 * \cdots * \gamma'_{p}$ est
un arc. Sa $h$--longueur est clairement  inf\'erieure ou \'egale \`a
$p=d_h(x,y)$.
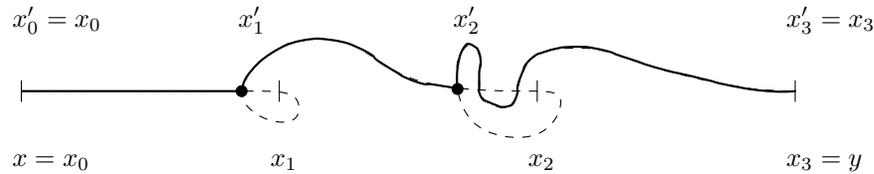
\begin{figure}[htp]\anchor{fig.simplification}
\centerline{\hbox{\input{fig-simplification.pstex_t}}}
\caption{\label{fig.simplification}Construction d'un arc g\'eod\'esique}
\end{figure}

R\'eciproquement, soit $\gamma$ un arc de $h$--longueur $k$ reliant deux
points $x$ et $y$. D'apr\`es la premi\`ere
partie de l'Affirmation~\ref{affi.arcs-et-domaines}, chaque sous-arc
libre de $\gamma$ est inclus dans un domaine de translation; on en
d\'eduit que $d_h(x,y) \leq k$, ce qui termine la preuve.
\end{proof}

\subsection{Quelques propri\'et\'es}
Rappelons qu'un hom\'eomorphisme de Brouwer n'a pas de point
p\'eriodique (Corollaire~\ref{coro.erre}). Par suite, $h^n$ est
encore un hom\'eomorphisme de Brouwer (pour tout entier $n$ non nul),
et on peut consid\'erer la distance $d_{h^n}$.
\begin{lemm}[Distance de translation pour $h^n$]\hspace{2cm}
\label{lemm.propriete-h-distance}
\begin{itemize}
\item Pour tous $x,y$ qui ne sont pas dans la m\^eme orbite de $h$,
et tout entier $n$,
$d_h(x,h^n(y))=d_h(x,y)$;
\item pour tout entier $n$ non nul, les distances $d_{h^n}$ et $d_h$
sont \'egales;
\item  tout arc $\gamma$ qui est g\'eod\'esique pour
$h$ est encore g\'eod\'esique pour $h^n$, et toute d\'ecomposition
minimale de
$\gamma$ pour $h$ est encore une d\'ecomposition minimale  de $\gamma$
pour
$h^n$.
\end{itemize}
\end{lemm}

\begin{proof}
La premi\`ere propri\'et\'e vient simplement de la d\'efinition de la
distance de translation, puisque les
domaines de translation sont invariants par $h$.

Pour la deuxi\`eme, il est clair que tout domaine de translation pour $h$
est encore un domaine de translation pour $h^n$ (car $\tau^n$ est
conjugu\'ee \`a $\tau$). Ceci prouve que $d_{h^n} \leq d_h$. D'autre part,
si $d_{h^n}(x,y)=1$, alors ou bien $x$ et $y$ sont dans la m\^eme
orbite de $h^n$ (et donc aussi dans la m\^eme orbite de $h$), ou bien
il existe un arc libre pour $h^n$ joignant
$x$ et $y$, et cet arc est inclus dans un domaine de translation de $h$
d'apr\`es le Lemme~\ref{lemm.ETDS}; dans les deux cas, on a
$d_h(x,y)=1$. On en d\'eduit
facilement que  $d_{h} \leq d_{h^n}$ (en utilisant par exemple une suite
g\'eod\'esique pour $h^n$).

Passons \`a la troisi\`eme propri\'et\'e. Soit   $\gamma=\gamma_1 *
\dots *
\gamma_d$ une d\'ecomposition minimale d'un arc g\'eod\'esique pour $h$.
Un arc libre pour $h$ est encore libre pour $h^n$
(Corollaire~\ref{coro.liberte}), donc chaque
sous-arc $\gamma_i$ est aussi libre pour $h^n$:  pour $h^n$, $\gamma_1 *
\dots *
\gamma_d$  est encore une d\'ecomposition de $\gamma$.
D'autre part, on vient de montrer que $d_{h^n}=d_h$, donc la
$h^n$--distance
entre les deux extr\'emit\'es de $\gamma$ vaut encore $d$, et la
$h^n$--longueur de $\gamma$ ne peut pas \^etre strictement plus petite
que $d$ (Corollaire~\ref{coro.definition-arcs}).
Ceci prouve que, pour $h^n$,  la d\'ecomposition est minimale, et que
$\gamma$ est
un arc g\'eod\'esique.
\end{proof}

\begin{lemm}
\label{lemm.iteres-arc-geo}
On se donne une d\'e\-com\-pos\-i\-tion min\-i\-male $\gamma=\gamma_1 *
\dots * \gamma_k$
d'un arc g\'eod\'esique ($k \geq 2$).
\begin{itemize}
\item Pour tous $i,j$ entre $1$ et $k$, pour tout
entier $n$, si $|i-j|>1$, alors
$h^n(\gamma_i) \cap \gamma_j = \emptyset$;
\item pour tout $i$ entre $1$ et $k-1$, pour tout entier $n$,
l'arc $\gamma_i * \gamma_{i+1}$ rencontre son image par $h^n$.
\end{itemize}
\end{lemm}
\begin{proof}
Prouvons le premier point. On note $(x_0, \dots, x_k)$ les sommets de
la d\'ecomposition. On peut supposer que $i<j$.
 Supposons que $h^n(\gamma_i)$ rencontre
$\gamma_j$, et soit $z$ un point de l'intersection; on a alors
$$
\begin{array}{rcl}
d_h(x_{i-1},x_{j})	& =	&  d_h(h^n(x_{i-1}),x_{j}) \\
		& \leq	&  d_h(h^n(x_{i-1}),z) + d_h(z, x_{j})	\\
		& \leq	&  1+1=2.
\end{array}
$$
Comme $\gamma$ est un arc g\'eod\'esique, ceci n'est possible que si
$|i-j|\leq1$. Le premier point est donc prouv\'e.

Montrons le second point. Si l'arc $\gamma_i * \gamma_{i+1}$ est libre
pour $h^n$, alors d'apr\`es le Lemme~\ref{lemm.propriete-h-distance},
$d_h(x_{i-1},x_{i+1})=d_{h^n}(x_{i-1},x_{i+1})=1$,
ce qui contredit le fait que $\gamma$ est un arc g\'eod\'esique.
\end{proof}

\subsubsection*{Les $h$--boules}
\begin{defi}
Pour tout point $x$ du plan et tout entier strictement positif $r$,
on notera
$B_r(x)$ la boule de centre $x$ et de rayon $r$ pour la distance
$d_h$:
$$B_r(x)=\{y	:   d_h(x,y) \leq r \}$$
On dira que $B_r(x)$ est la
\emph{$h$--boule de centre $x$ et de rayon $r$}.
\end{defi}
Cet ensemble est un ouvert connexe du plan (pour la topologie
usuelle), et est invariant par $h$: en effet,
d'apr\`es la d\'efinition de $d_h$, $B_1(x)$ est l'union des domaines de
translation contenant $x$, et $B_{r}(x)$ est l'union des domaines
de translation rencontrant $B_{r-1}(x)$.

\subsection{Fl\`eches horizontales}
\label{ssec.fleches}
Dans ce paragraphe, nous utilisons le lemme de Franks pour d\'efinir les
fl\`eches horizontales.
\subsubsection*{D\'efinition}
On se donne une d\'ecomposition minimale  $\gamma=\gamma_1 * \dots *
\gamma_d$ d'un arc g\'eod\'esique ($d \geq 1$).
\begin{affi}\label{affi.deux-possibilites}
 Soit $i$ un entier entre $1$ et $d-1$.
Une et une seule des deux possibilit\'es suivantes est
v\'erifi\'ee:
\begin{enumerate}
\item $h(\gamma_i) \cap \gamma_{i+1}\neq \emptyset$,
\item $\gamma_i \cap h( \gamma_{i+1}) \neq \emptyset$.
\end{enumerate}
\end{affi}

\begin{figure}[htp]\anchor{fig.22choseslune}
\centerline{\hbox{\input{fig-22choseslune.pstex_t}}}
\caption{\label{fig.22choseslune}D\'efinition d'une fl\`eche, \textit{via}
le lemme de Franks}
\end{figure}
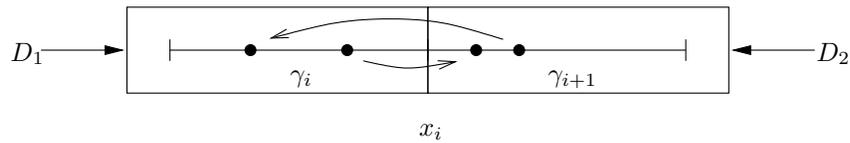

\begin{proof}[\proofname\ de l'Affirmation~\ref{affi.deux-possibilites}]
Puisque la d\'ecomposition est minimale, l'arc $\gamma_i *
\gamma_{i+1}$ ne peut pas \^etre libre. Puisque, par contre,
$\gamma_i$ et $\gamma_{i+1}$ sont libres, l'une des deux possibilit\'es
de la remarque doit \^etre v\'erifi\'ee.

Montrons que les deux possibilit\'es ne peuvent pas \^etre v\'erifi\'ees
simultan\'ement. On note $x_i$ le sommet commun aux deux arcs. En
utilisant le th\'eor\`eme de Schoenflies, on peut supposer que
$\gamma$ est un segment euclidien (\cf \figref{fig.22choseslune}). On
peut alors \'epaissir les deux arcs pour trouver deux disques topologiques
ferm\'es $D_1$, $D_2$ tels que
\begin{itemize}
\item $\gamma_i \setminus\{x_i\} \subset \inte(D_1)$ et   $\gamma_{i+1}
\setminus\{x_i\} \subset \inte(D_2)$;
\item les int\'erieurs de $D_1$ et $D_2$ sont disjoints;
\item $D_1$ et $D_2$ sont libres.
\end{itemize}
Supposons maintenant que la premi\`ere possibilit\'e de l'affirmation est
v\'erifi\'ee, c'est-\`a-dire qu'il existe un point $x$	de $\gamma_i$
dont l'image $h(x)$ est dans
$\gamma_{i+1}$. Alors $x \neq x_i$ car $\gamma_{i+1}$ est libre, et
$h(x) \neq x_i$ car $\gamma_{i}$ est libre. Par cons\'equent $x$ et
$h(x)$ appartiennent respectivement aux int\'erieurs de $D_1$ et
$D_2$, autrement dit $h(\inte(D_1))$ rencontre $\inte(D_2)$. De
m\^eme, si la seconde possibilit\'e de l'aff\-ir\-ma\-tion  est
v\'erifi\'ee, alors $h(\inte(D_2))$ rencontre $\inte(D_1)$.
Si ces deux rencontres avaient	lieu simultan\'ement,
la suite $(D_1,D_2)$ serait une cha\^\i ne de disques p\'eriodique,
ce qu'interdit	le lemme de Franks~\ref{lemm.franks}.
\end{proof}

Rappelons qu'un \emph{mot \`a $p$ lettres dans l'alphabet $A$} est une
suite de $p$ \'el\'ements de $A$, autrement dit une application de
l'ensemble $\{1, \dots, p\}$ dans $A$. On notera toujours $M_i$
la $i$\`eme
lettre du mot $M$.
\begin{defi}\label{defi.fleches-horizontales}
On appelle \emph{mot horizontal de $\gamma_1 * \dots *
\gamma_d$}, et	 on note
$$M^h_\leftrightarrow(\gamma_1 * \dots *
\gamma_d),$$
 le mot de $d-1$ lettres dans l'alphabet
$\{\leftarrow, \rightarrow \}$, dont la $i$\`eme lettre est:
\begin{enumerate}
\item $M^h_\leftrightarrow(\gamma_1 * \dots *
\gamma_d)_i:=\rightarrow$ si c'est la premi\`ere possibilit\'e qui est
satisfaite dans l'Affirmation~\ref{affi.deux-possibilites} ci-dessus;
\item $M^h_\leftrightarrow(\gamma_1 * \dots *
\gamma_d)_i:=\leftarrow$ si c'est la deuxi\`eme.
\end{enumerate}
S'il
n'y a pas d'ambigu\"{\i}t\'e sur l'hom\'eomorphisme, on \'ecrira
simplement
$$M_\leftrightarrow(\gamma_1 * \dots * \gamma_d).$$
\end{defi}

\subsection{\'Enonc\'e pr\'ecis du Th\'eor\`eme~\ref{theo.fleches-hori}}
\textit{A priori}, le mot horizontal d\'epend de l'arc g\'eod\'esique
$\gamma$, et
m\^eme du choix d'une d\'ecomposition de $\gamma$. Le
Th\'eor\`eme~\ref{theo.fleches-hori}
dit qu'en r\'ealit\'e, il ne d\'epend que des extr\'emit\'es de $\gamma$.
\setcounter{theobis}{1}
\begin{theobis}\label{theoBbis}
Soient $x$ et $y$ deux points du plan, et
 $\gamma=\gamma_1 * \dots * \gamma_d$ et $\gamma'=\gamma'_1 * \dots *
\gamma'_d$ deux d\'ecompositions minimales d'arcs g\'eod\'esiques joignant
$x$ et $y$.
Alors $M_\leftrightarrow(\gamma_1 * \dots *
 \gamma_d)=M_\leftrightarrow(\gamma'_1 * \dots *
\gamma'_d)$.
\end{theobis}

La preuve du th\'eor\`eme occupera les Sections~\ref{sec.crg},
\ref{sec.topologie-disque} et
\ref{sec.preuve-crg}. En admettant provisoirement le th\'eor\`eme,
on peut donner la
d\'efinition suivante.
\begin{defi}\label{def.mot-couple}
Soient $x$ et $y$  deux points distincts du plan.
On appellera \emph{mot horizontal du couple $(x,y)$}, et on notera
$M^h_\leftrightarrow(x,y)$, le mot $M_\leftrightarrow(\gamma_1 * \dots *
\gamma_d)$, o\`u
$\gamma_1 * \dots * \gamma_d$ est n'importe quelle d\'ecomposition
minimale de n'importe quel arc g\'eod\'esique joignant $x$ \`a $y$. S'il
n'y a pas d'ambigu\"{\i}t\'e sur l'hom\'eomorphisme, on \'ecrira
simplement
$M_\leftrightarrow(x,y)$.
\end{defi}
Remarquons que le mot horizontal entre deux points \`a $h$--distance $1$
est vide.

\subsection{Quelques propri\'et\'es \'el\'ementaires du mot horizontal}
Sur l'alphabet $\{\leftarrow, \rightarrow \}$, on d\'efinit l'op\'eration
``oppos\'ee'', not\'ee $-$, par
$$-\leftarrow \ := \ \rightarrow\quad\mbox{et}\quad -\rightarrow\ := \
\leftarrow.$$
 Ceci induit naturellement une op\'eration sur
les mots sur cet alphabet, encore not\'ee $-$, qui consiste \`a changer
chaque lettre du mot en son oppos\'ee. Le mot $-M$ est appel\'e
\emph{oppos\'e} du mot $M$.
Si $M$ est un mot de longueur $l$ sur cet alphabet, le \emph{palindrome
de $M$}
est le mot $\pali(M)$ d\'efini par $\pali(M)_i=M_{l+1-i}$.
Enfin, pour tout arc $\gamma$, on note $-\gamma$ l'arc ayant
l'orientation oppos\'ee \`a celle de $\gamma$: en termes de param\'etrage,
on a
$(-\gamma)(t)=\gamma(1-t)$.

On se donne une d\'ecomposition minimale $\gamma_1 * \dots *
\gamma_d$  d'un arc g\'eod\'esique $\gamma$ pour
$h$. Par abus, on notera parfois $\gamma$ la d\'ecomposition elle-m\^eme
(et $-\gamma$ la d\'ecomposition induite de mani\`ere \'evidente sur
l'arc $-\gamma$).
\begin{lemm}~\label{lemm.proprietes-mot-horizontal}
\begin{enumerate}
\item	 $M^{h^n}_\leftrightarrow(\gamma) =
M^h_\leftrightarrow(\gamma) =  -M^{h^{-n}}_\leftrightarrow(\gamma)$
(pour tout $n \geq 1$);
\item $M_\leftrightarrow(-\gamma)=-\pali(M_\leftrightarrow(\gamma))$.
\end{enumerate}
\end{lemm}

\begin{proof}
Tout d'abord, les liens entre $M^h_\leftrightarrow(\gamma)$,
$M^{h^{-1}}_\leftrightarrow(\gamma)$ et $M^h_\leftrightarrow(-\gamma)$
suivent imm\'ediatement de la d\'efinition du mot
horizontal. Il reste donc uniquement \`a \'etablir le lien entre
$M^{h^n}_\leftrightarrow(\gamma)$ et
$M^h_\leftrightarrow(\gamma)$ pour $n \geq 1$.

Soit  $i$ un entier entre $1$ et $d-1$.
D'apr\`es le Lemme~\ref{lemm.propriete-h-distance}, pour $h^n$, l'arc
$\gamma$ est encore un arc g\'eod\'esique dont $\gamma_1 * \dots *
\gamma_d$  est une d\'ecomposition minimale; notamment, le mot
$M^{h^n}_\leftrightarrow(\gamma)$ est bien d\'efini.
Montrons que les mots horizontaux associ\'es \`a $\gamma$ pour $h$ et
$h^n$ co\"{\i}ncident; pour cela, on s'int\'eresse \`a la $i$\`eme lettre
de ces deux mots. Pour fixer les id\'ees, supposons que
$M^h_\leftrightarrow(\gamma)_i=\leftarrow$, c'est-\`a-dire que
$h(\gamma_i) \cap \gamma_{i+1}= \emptyset$. \`A l'aide du
th\'eor\`eme de Schoenflies, on trouve deux disques topologiques
ferm\'es libres  $D_i$ et $D_{i+1}$, contenant respectivement les arcs
$\gamma_i$ et $\gamma_{i+1}$ dans leurs int\'erieurs, et tels que
$h(D_i) \cap D_{i+1}= \emptyset$ (\cf \figref{fig.fleches-hn}).
\begin{figure}[htp]\anchor{fig.fleches-hn}
\centerline{\hbox{\input{fig-fleches-hn.pstex_t}}}
\caption{\label{fig.fleches-hn}Si
$M^h_\leftrightarrow(\gamma)_i=\leftarrow$\ldots{}}
\end{figure}
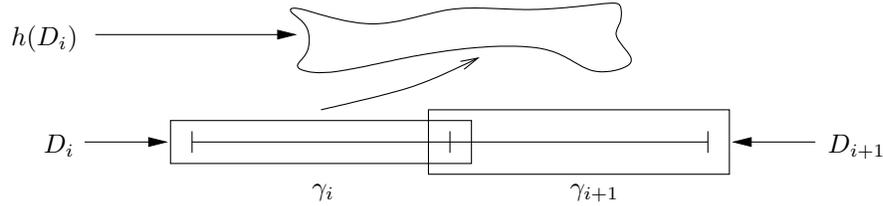
D'apr\`es le lemme de Franks
(Lemme~\ref{lemm.intervalle}), l'ensemble des entiers $n$ tels que
$h^n(\inte(D_{i}))$ rencontre $\inte(D_{i+1})$ est un intervalle de
$\bbZ$. Comme
 $\inte(D_{i+1})$ rencontre $\inte(D_{i})$  mais pas
$h(\inte(D_{i}))$, on en d\'eduit que pour tout entier $n
\geq 1$, $h^n(\inte(D_i)) \cap \inte(D_{i+1})  = \emptyset$, donc
aussi que
$ h^n(\gamma_i) \cap \gamma_{i+1} = \emptyset$.
 Puisque la d\'ecomposition est
minimale pour $h^n$, on doit alors avoir $\gamma_i \cap
h^n(\gamma_{i+1}) \neq \emptyset$
(par Affirmation~\ref{affi.deux-possibilites});
autrement dit,	$M^{h^n}_\leftrightarrow(\gamma)_i =\leftarrow=
M^h_\leftrightarrow(\gamma)_i$, ce que l'on voulait.
\end{proof}

\section{Composantes de Reeb}
\label{sec.crg}
Dans cette partie, nous g\'en\'eralisons le concept de composantes
de Reeb des feuilletages du plan (voir l'introduction), nous
\'enon\c{c}ons un
th\'eor\`eme d'ex\-ist\-ence et d'unicit\'e, et nous en d\'eduisons le
Th\'eor\`eme~\ref{theo.fleches-hori} (le mot horizontal d'un arc
g\'eod\'esique ne d\'epend que de ses extr\'emit\'es).
Le choix de la d\'efinition des composantes de Reeb est assez naturel,
il utilise le concept classique d'\emph{ensemble singulier} (d\^u \`a
Ker\'ekj\'art\'o). Un des points importants est qu'il faut
maintenant d\'efinir une composante de Reeb comme \'etant associ\'ee
\`a un couple de points $(x,y)$. Le Th\'eor\`eme~\ref{theoAbis}
affirme alors
l'existence d'exactement $d-1$ composantes de Reeb \emph{minimales}
associ\'ees \`a $(x,y)$, o\`u $d$ est la distance de translation
entre $x$ et $y$, et ces composantes sont donn\'ees par une
``formule'' topologique explicite.
Notons que la minimalit\'e est essentielle pour obtenir un nombre
fini de composantes associ\'ees \`a $(x,y)$ (voir
les exemples de l'Appendice~\ref{sec.exemples}). Elle jouera
\'egalement un r\^ole central dans la d\'efinition des fl\`eches
verticales (Section~\ref{sec.geodesiques-bout},
\ref{sec.fleches-verticales} et \ref{sec.continus-minimaux}). Par
contre, elle n'intervient pas dans la preuve du
Th\'eor\`eme~\ref{theo.fleches-hori}.

\subsection{Ensemble singulier}
\label{ssec.singulier}
\begin{defi}
Soient $x$ et $y$ deux points du plan. On dira que le couple $(x,y)$
 est \emph{singulier} si $x$ et $y$ n'appartiennent pas \`a la m\^eme
 orbite de $h$, et si
 pour tous voisinages $V_x$ et $V_y$ de $x$ et $y$ respectivement, il
 existe un entier
\emph{positif} $n$ tel que $h^n(V_x) \cap V_y \neq \emptyset$.
 On notera $\sing(h)$, et on
appellera \emph{ensemble singulier de $h$}, l'ensemble des couples
singuliers de $h$.
\end{defi}

Voici quelques exemples. La translation $\tau\co (x,y) \mapsto (x+1,y)$
n'admet aucun couple de points singuliers. Pour l'hom\'eomorphisme de
Reeb (repr\'esent\'e sur la \figref{fig.Reeb}), l'ensemble singulier
est l'ensemble des couples $(x,y)$ o\`u $x$ et $y$ appartiennent
respectivement aux bords sup\'erieur et inf\'erieur de la bande $B$.
En fait, on peut voir que si $h$ est le temps $1$ d'un flot, alors
l'ensemble singulier est la r\'eunion des produits $\Delta_0 \times
\Delta_1$ des bords de composantes de Reeb $(\Delta_0,\Delta_1)$ (au
sens de la d\'efinition donn\'ee en introduction de ce texte).

Il suit de la d\'efinition qu'un couple $(x,y)$ est singulier pour
$h^{-1}$ si et seulement si le couple $(y,x)$ est singulier pour
$h$. Nous verrons plus loin que $h^n$ a le m\^eme ensemble
singulier que $h$ si $n \geq 1$ (voir les commentaires qui suivent le
Lemme~\ref{lemm.singulier-limite}).

Rappelons que $B_1(x)$ d\'esigne la boule de rayon $1$ et de centre
$x$ pour la distance de translation de $h$, \cad\ la r\'eunion des
domaines de translation contenant $x$. Le lemme suivant
\'etablit un premier lien entre l'ensemble singulier et la
distance $d_h$.
\begin{lemm}[Couples singuliers et $h$--distance]
\label{lemm.distance-singulier}
Soit $(x,y)$ un couple singulier de $h$. Alors $d_h(x,y)=2$, et $y \in
\partial B_1(x)$.
\end{lemm}
\begin{proof}
Une translation n'admet pas de couple singulier, par suite il ne peut pas
exister de domaine de translation de $h$ contenant simultan\'ement
$x$ et $y$. Ceci prouve que $d_h(x,y) > 1$.

Soit $O$ un domaine de translation contenant $x$: c'est un voisinage
de $x$, et il est invariant par $h$. La d\'efinition des
couples singuliers entra\^\i ne que tout voisinage de $y$ rencontre
$O$. Autrement dit $y$ est dans l'adh\'erence de $O$, mais $O$ est
inclus dans $B_1(x)$, donc $y \in \adhe(B_1(x))$. Comme $B_1(x)$ ne
contient pas $y$, c'est que $y \in \partial B_1(x)$.
On en d\'eduit que tout domaine de translation qui contient $y$
rencontre $B_1(x)$, et donc que $d_h(x,y)=2$.
\end{proof}

\subsection{D\'efinition g\'en\'erale des composantes de Reeb}
\begin{figure}[htp]\anchor{fig.compo-Reeb}
\centerline{\hbox{\input{fig-compo-Reeb.pstex_t}}}
\caption{\label{fig.compo-Reeb}Une composante de Reeb non
d\'eg\'en\'er\'ee}
\end{figure}
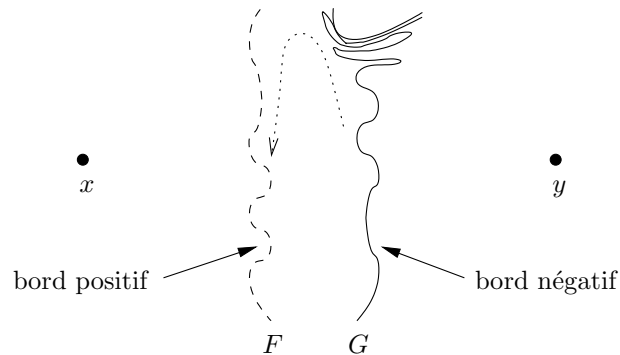

\begin{defi}[\rm(\figref{fig.compo-Reeb})]
Soient $x$ et $y$ deux points, et $F$ et $G$ deux sous-ensembles ferm\'es
connexes du plan. On dit que  $(F,G)$ est une \emph{composante de Reeb
pour $(x,y)$} si
\begin{enumerate}
\item $F\times G \subset \sing(h)$  ou $G \times F \subset \sing(h)$;
\item $F$ contient $x$, ou bien $F$  s\'epare $x$ et $G \cup \{y\}$;
\item $G$ contient $y$, ou bien $G$ s\'epare $y$ et $F \cup \{x\}$.
\end{enumerate}
Les ensembles $F$ et $G$ sont appel\'es \emph{bords} de la
composante: si $F\times G \subset \sing(h)$, $F$ en est le \emph{bord
n\'egatif}, $G$ le \emph{bord positif};  si $G \times F \subset
\sing(h)$, c'est le contraire.
\end{defi}
Il suit imm\'ediatement de la d\'efinition que les ensembles $F$ et
$G$ sont disjoints, puisqu'un couple $(x,x)$ n'est jamais singulier.

\begin{defi}
La composante de Reeb  est dite \emph{minimale} si,
 de plus,
\begin{itemize}
\item $F=\{x\}$, ou bien $F$ est minimal pour l'inclusion parmi les
ensembles ferm\'es connexes qui s\'eparent $x$ et $G \cup \{y\}$;
\item $G=\{y\}$, ou bien $G$ est minimal pour l'inclusion parmi les
ensembles ferm\'es connexes qui s\'eparent $y$ et $F \cup \{x\}$.
\end{itemize}
Autrement dit, $(F,G)$ est minimale si	pour
toute composante de Reeb  $(F',G')$ pour  $(x,y)$, si
$F' \subset F$ et $G' \subset G$, alors $F'=F$ et $G'=G$.
\end{defi}

\begin{defi}
La composante $(F,G)$ est dite \emph{d\'eg\'en\'er\'ee} si  $F$
contient $x$ ou  $G$ contient $y$.
\end{defi}
Voici l'exemple le plus simple de composante de Reeb: si $(x,y) \in
\sing(h)$, alors $(\{x\},\{y\})$ est une
composante de Reeb  d\'eg\'en\'er\'ee et minimale
pour $(x,y)$. L'Appendice~\ref{sec.exemples} contient des
exemples plus substantiels (on y montre en particulier qu'un couple
$(F,G)$ peut \^etre une composante de Reeb  minimale pour $(x,y)$,
et \^etre une composante de Reeb \emph{non} minimale pour un autre couple
$(x',y')$, m\^eme dans le cas non d\'eg\'en\'er\'e).

Remarquons que toute composante de Reeb pour $(x,y)$  ``contient'' une
 composante de Reeb  minimale pour $(x,y)$. En effet, le lemme de Zorn
 entra\^{\i}ne que
tout ferm\'e s\'eparant $x$ et $y$ contient un ferm\'e s\'eparant $x$
et $y$
et minimal pour l'inclusion parmi les ensembles de ce
type. En r\'ealit\'e, nous n'utiliserons pas cette remarque, puisque
nous allons
 obtenir les composantes de Reeb minimales de mani\`ere constructive.

\subsection{Construction des composantes de Reeb minimales}
\label{ssec.construction}
Soient $x$ et $y$ deux points du plan, tels que $d=d_h(x,y) \geq 2$.
On a
$$
\adhe(B_1(x)) \subset B_2(x) \subset \adhe(B_2(x)) \subset B_3(x)
\subset \cdots \subset B_{d-1}(x).
$$
Le point $y$ n'est dans aucun de ces ensembles. Par suite, pour $i=1,
\dots, d-2$, on peut consid\'erer l'ensemble
$O_y(\adhe(B_i(x)))$, d\'efini comme la composante connexe de l'ensemble
$\bbR^2 \setminus
\adhe(B_{i}(x))$ qui contient $y$. On pose alors
$$
F_1(x,y)=\partial O_y(\adhe(B_1(x))),
\dots,
F_{d-2}(x,y)=\partial O_y(\adhe(B_{d-2}(x))).
$$
Dans le cas o\`u l'ensemble $\partial B_{d-1}(x)$ ne contient pas non
plus le point $y$, on peut d\'efinir de la
m\^eme mani\`ere $O_y(\adhe(B_{d-1}(x)))$ et $F_{d-1}(x,y)$. Dans le cas
contraire, on pose $F_{d-1}(x,y)=\{y\}$.
Enfin, on d\'efinit sym\'etriquement les ensembles $O_x(\adhe(B_i(y)))$
et $F_i(y,x)$ pour
$i=1, \dots, d-1$.
D\'ecrivons imm\'ediatement quelques propri\'et\'es
\'el\'ementaires, concernant la topologie, puis la dynamique de ces
ensembles.
\begin{lemm}[Topologie des ensembles $F_i(x,y)$]
\label{lemm.topologie-crg}
  Pour tous $x,y$ et pour tout $i$ entre $1$ et $d-1$,
 l'ensemble $F_i(x,y)$ est connexe et  contenu dans  $\partial B_i(x)$.
Si  de plus $F_i(x,y) \neq \{y\}$, alors
cet ensemble s\'epare $x$ et $y$, et poss\`ede la propri\'et\'e de
 minimalit\'e suivante: tout ensemble $F$, contenu dans $\partial
B_i(x)$ et s\'eparant $x$ et $y$, contient $F_i(x,y)$.
\end{lemm}
\begin{proof}[\proofname\ du Lemme~\ref{lemm.topologie-crg}]
On se place dans le cas non trivial  ($F_i(x,y) \neq \{y\}$).
La connexit\'e de $F_i(x,y)$ r\'esulte des propri\'et\'es de la topologie
plane (appliquer le Th\'eor\`eme~\ref{theo.connexite} de l'appendice \`a
l'ouvert connexe $B_i(x)$). On a les inclusions
$$
F_i(x,y)=\partial  O_y(\adhe(B_{i}(x))) \subset \partial \adhe
(B_i(x)) \subset \partial B_i(x),
$$
l'\'egalit\'e a lieu par d\'efinition, la premi\`ere inclusion suit de la
Proposition~\ref{prop.bumpwire} de l'appendice, la derni\`ere est une
propri\'et\'e de topologie g\'en\'erale.
L'ensemble $F_i(x,y)$ s\'epare $x$ et $y$ puisqu'il est la fronti\`ere
d'un ensemble qui contient $y$ mais pas $x$.
Il reste \`a montrer la propri\'et\'e de minimalit\'e.
 Soit $z$  un point
quelconque de $F_i(x,y)$. On voit facilement que l'ensemble
$B_i(x) \cup \{z\} \cup O_y(\adhe(B_{i}(x)))$ est connexe. Un ensemble
$F$ inclus dans $\partial B_i(x)$ est disjoint de $B_i(x)$ et de
$O_y(\adhe(B_{i}(x)))$; par cons\'equent, si $F$ s\'epare $x$ et $y$, il
doit contenir le point $z$. On en d\'eduit que tout ensemble $F$ inclus
dans  $\partial B_i(x)$ et  s\'eparant $x$ et $y$ doit contenir
$F_i(x,y)$, ce que l'on voulait.
\end{proof}

Il n'est pas encore clair que les ensembles $F_i(x,y)$ et
$F_{d-i}(y,x)$ soient disjoints (comme sugg\'er\'e par la
\figref{fig-theoreme-crg} ci-dessous). Cependant, ceci sera une
cons\'equence du Th\'eor\`eme~\ref{theoAbis}.

\begin{lemm}~
\label{lemm.invariance-crg}
 Les ensembles $F_1(x,y), \dots, F_{d-2}(x,y)$ sont invariants
par $h$.
\end{lemm}

\begin{proof}[\proofname\ du Lemme~\ref{lemm.invariance-crg}]
Tout d'abord,  les ensembles $\adhe(B_i(x))$ sont invariants, donc les
ensembles $O_y(\adhe(B_i(x)))$ sont libres ou invariants (comme
composantes
connexes d'ensembles invariants).
 Si $1 \leq i \leq d-2$,  la $h$--boule $B_1(y)$ est disjointe de
$B_i(x)$ (in\'egalit\'e triangulaire), elle est donc incluse dans
l'ensemble
$O_y(\adhe(B_i(x)))$: celui-ci ne peut pas \^etre libre, puisqu'il
contient
l'ensemble invariant $B_1(y)$. Il est donc invariant.
L'ensemble $F_i(x,y)$ est aussi invariant comme fronti\`ere d'un ensemble
invariant.\footnote{Dans une premi\`ere version de ce texte, on affirmait
que l'ensemble $F_{d-1}(x,y)$ est invariant ou libre sous $h$. Suite \`a
une question du referee, l'auteur s'est aper\c cu de son erreur: on peut
construire des exemples ou cet ensemble n'est ni libre ni invariant.}
\end{proof}

\subsection{\'Enonc\'e pr\'ecis du Th\'eor\`eme~\ref{theo.crgm}}
\setcounter{theobis}{0}
\begin{theobis}\label{theoAbis}
Soient $x$ et $y$ deux points du plan, et $d=d_h(x,y)$.
Alors il existe exactement $d-1$  composantes de Reeb  minimales
pour $(x,y)$; ces composantes sont
 les  couples $(F_{d-i}(y,x),F_i(x,y))$ pour $i$ variant entre $1$
 et $d-1$.
\end{theobis}
\begin{figure}[htp]\anchor{fig-theoreme-crg}
\centerline{\hbox{\input{fig-theoreme-crg.pstex_t}}}
\caption{\label{fig-theoreme-crg}Illustration du
Th\'eor\`eme~\ref{theo.crgm}}
\end{figure}
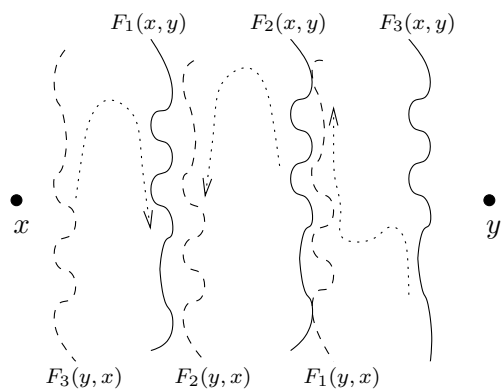
Le th\'eor\`eme est illustr\'e symboliquement par la
\figref{fig-theoreme-crg}, et plus concr\`etement par les exemples
de l'Appendice~\ref{sec.exemples}.
La preuve du th\'eor\`eme occupera la Section~\ref{sec.preuve-crg}.
Remarquons que parmi les $d-1$ composantes minimales, seules les deux
composantes ``extr\^emes''
$$(F_{d-1}(y,x),F_1(x,y)) \quad\text{et}\quad (F_{1}(y,x),F_{d-1}(x,y))$$
sont
susceptibles d'\^etre d\'eg\'en\'er\'ees. Ceci suit
imm\'ediatement de la d\'efinition des ensembles $F_i(x,y)$.
Notons au passage que les composantes d\'eg\'en\'er\'ees sont ``rares'':
en
effet, H. Nakayama a montr\'e que l'ensemble des points faisant partie
d'un
couple de points singuliers est topologiquement maigre (\cite{naka2});
g\'en\-\'e\-rique\-ment, les points $x$ et $y$ sont non singuliers, et aucune
des composantes de Reeb associ\'ees \`a $(x,y)$ n'est d\'eg\'en\'er\'ee.
 De
m\^eme, seuls les bords ``extr\^emes'' $F_{d-1}(y,x)$ et $F_{d-1}(x,y)$
peuvent ne pas \^etre invariants: ceci suit du
Lemme~\ref{lemm.invariance-crg}.

\subsection{Preuve du Th\'eor\`eme~\ref{theoBbis}: invariance
des fl\`eches
horizontales}
\label{ssec.preuve-1}
Dans cette section, nous admettons provisoirement le
Th\'eor\`eme~\ref{theoAbis} et nous en d\'eduisons le
Th\'eor\`eme~\ref{theoBbis}.
Commençons par donner l'id\'ee de la preuve.
Pour simplifier, consid\'erons d'abord le cas o\`u
 $h$ est  un hom\'eomorphisme de Brouwer qui
pr\'eserve chaque feuille d'un	feuilletage $\cal F$ (par exemple, $h$
peut \^etre obtenu en  int\'egrant un  champ  de vecteurs qui ne
s'annule  pas, \cf\ l'introduction).
 On voit facilement
que la $h$--distance entre les deux bords d'une
composante de Reeb du feuilletage est \'egale \`a $2$, et que le mot
horizontal
d'un arc g\'eod\'esique joignant le bord n\'egatif de la composante
\`a son bord positif vaut $(\rightarrow)$. Consid\'erons alors un
arc g\'eod\'esique $\gamma$, d'extr\'emit\'es $x$ et $y$, de longueur
$d \geq 2$, muni d'une
d\'ecomposition minimale de sommets $x_0=x, \dots, x_d=y$
(\cf\ \figref{fig.feuilletage-geo} de l'introduction). Cet arc
franchit
successivement les $d-1$ composantes de Reeb associ\'ees au couple
$(x,y)$, dont l'existence est assur\'ee par le
Th\'eor\`eme~\ref{theo.crgm}.
Puisque les deux  bords d'une composante de Reeb sont \`a $h$--distance
$2$, le sommet $x_i$ est situ\'e entre les deux bord de la
$i$--\`eme composante rencontr\'ee.
 Si $\gamma$ rencontre	le bord
n\'egatif  de cette composante avant le bord positif, alors la
$i$--\`eme fl\`eche du mot $M_\leftrightarrow(\gamma)$
 sera $\rightarrow$; elle sera $\leftarrow$
si $\gamma$ rencontre en premier lieu le bord positif. Il est
maintenant clair que tout autre arc g\'eod\'esique ayant les m\^emes
extr\'emit\'es que $\gamma$ donnera la m\^eme suite de fl\`eches.

La preuve dans le cas g\'en\'eral s'inspire fortement du cas des
 feuilletages. On montre d'abord qu'un arc g\'eod\'esique joignant $x$
 \`a $y$ doit
franchir successivement les $d-1$ composantes de Reeb minimales dans
 l'ordre attendu
(voir Lemme~\ref{lemm.ordre-crg}). Chaque franchissement d'une
composante de Reeb impose la pr\'esence d'un sommet,
 et impose aussi la fl\`eche correspondante
(Lemme~\ref{lemm.fleche-singulier}).

\begin{lemm}[Franchissement des composantes de Reeb]
\label{lemm.ordre-crg}
Soit $\gamma=\gamma_1 * \cdots * \gamma_d$ une d\'ecomposition
minimale d'un arc g\'eod\'esique joignant $x$ \`a $y$.
Alors pour tout entier $i$ entre $1$ et $d-1$, l'ensemble  $F_i(x,y)$
rencontre $\gamma$, et $F_i(x,y) \cap \gamma \subset \gamma_{i+1}$.
Si de plus $i \neq d-1$, alors $F_i(x,y) \cap \gamma \subset
\inte(\gamma_{i+1})$.
\end{lemm}
On a bien s\^ur un r\'esultat sym\'etrique en inversant les r\^oles de $x$
et $y$: notamment, $F_{d-i}(y,x)\cap \gamma  \subset \gamma_i$.
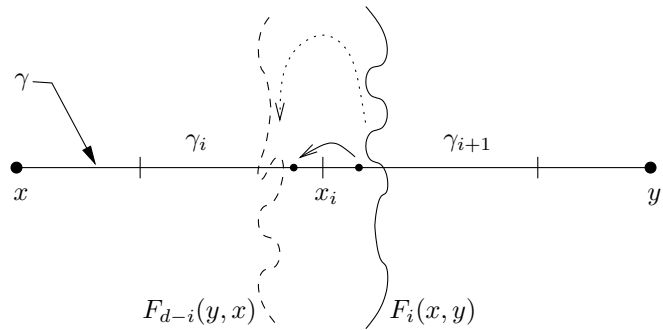
\begin{figure}[htp]\anchor{fig.rencontre-geodesique-crg}
\centerline{\hbox{\input{fig-rencontre-geodesique-crg.pstex_t}}}
\caption{\label{fig.rencontre-geodesique-crg}Positions relatives d'un arc
g\'eod\'esique et d'une composante de Reeb  minimale}
\end{figure}

Le second lemme relie la fl\`eche des arcs de longueur~2 \`a la
pr\'esence d'un couple singulier.
\begin{lemm}[Ensemble singulier et fl\`eches horizontales]
\label{lemm.fleche-singulier}
 Soit $\gamma=\gamma_1 * \gamma_2$ une d\'ecomposition minimale  d'un
 arc g\'eod\'esique de longueur~2. Alors
$$
\begin{array}{rcl} M_\leftrightarrow(\gamma_1 * \gamma_2) = (\to) \  \ \
&\Leftrightarrow & \ \ \ \gamma_1 \times
\gamma_2 \cap \sing(h) \neq \emptyset, \\
 M_\leftrightarrow(\gamma_1 * \gamma_2) = (\leftarrow) \  \ \ &
\Leftrightarrow & \ \ \ \gamma_2 \times
\gamma_1 \cap \sing(h) \neq \emptyset.
\end{array}
$$
\end{lemm}

\begin{proof}[\proofname\ du Lemme~\ref{lemm.ordre-crg}]
 L'ingr\'edient principal de la preuve est l'in\'egalit\'e
triangulaire.
Soit  $i$ un entier entre $1$ et $d-1$.
Notons $x_0=x, \dots, x_d=y$ les sommets de la d\'ecomposition.
Puisque l'ensemble $F_i(x,y)$ s\'epare $x$ et $y$ ou bien vaut $\{y\}$,
l'arc $\gamma$ doit rencontrer cet ensemble. Montrons que la rencontre
se fait n\'ecessairement entre $x_i$ et $x_{i+1}$. En effet:
\begin{itemize}
\item le sous arc $[x x_i]_\gamma$ est dans la $h$--boule $B_i(x)$ (car il
   contient $x$ et est r\'eunion de $i$ arcs libres), donc $\partial
B_i(x)$ ne rencontre pas $[x x_i]_\gamma$ (puisque les
$h$--boules sont des ouverts du plan);
\item  si $i < d-1$, le sous-arc $[x_{i+1} y]_\gamma$ est dans la
$h$--boule
   $B_{d-i-1}(y)$. Or d'apr\`es l'in\'egalit\'e triangulaire, on a
$B_{d-i-1}(y)\cap B_i(x) = \emptyset$,
   donc $\partial B_i(x)$ ne rencontre pas  $[x_{i+1} y]_\gamma$.
\end{itemize}
Puisque $F_i(x,y) \subset \partial B_i(x)$, on a bien
$$F_i(x,y) \cap \gamma \subset \ ]x_ix_{i+1}[_\gamma  \ =  \
\inte(\gamma_{i+1})
\quad\text{et}\quad
F_{d-1}(x,y) \cap \gamma \subset \ ]x_{d-1}x_{d}]_\gamma  \subset
\gamma_{d}$$
si $i<d-1$.
Ce que l'on voulait.
\end{proof}

\begin{proof}[\proofname\ du Lemme~\ref{lemm.fleche-singulier}]
On montre la premi\`ere \'equivalence, la seconde en d\'ecoule (en
changeant
$h$ en $h^{-1}$).

Montrons la premi\`ere implication. On suppose que
$M_\leftrightarrow(\gamma_1 * \gamma_2) = (\to)$. D'apr\`es les
Lemmes~\ref{lemm.propriete-h-distance} et
\ref{lemm.proprietes-mot-horizontal}, pour  $h^n$, $\gamma_1 *
\gamma_2$ est encore une d\'ecomposition minimale d'un arc g\'eod\'esique,
et le mot horizontal associ\'e est encore $(\to)$. Autrement dit, pour
tout entier $n >0$, on a $h^n(\gamma_1) \cap \gamma_2 \neq
\emptyset$. Par compacit\'e, on en d\'eduit facilement que $\gamma_1
\times \gamma_2$ contient un couple singulier.

Montrons l'implication r\'eciproque. Ceci revient \`a montrer
que si $M_\leftrightarrow(\gamma_1 * \gamma_2) = (\leftarrow)$, alors
  $\gamma_1 \times \gamma_2$ ne  contient pas de couple
singulier. Exactement comme dans la preuve du
Lemme~\ref{lemm.proprietes-mot-horizontal} (\cf
\figref{fig.fleches-hn}),
l'hypoth\`ese  $M_\leftrightarrow(\gamma_1 * \gamma_2) = (\leftarrow)$
permet de trouver deux disques topologiques ferm\'es $D_1$ et $D_2$,
dont les
int\'erieurs contiennent respectivement $\gamma_1$ et $\gamma_2$, et
tels que $h^n(\inte(D_1))$ soit disjoint de $\inte(D_2)$ pour tout
entier $n>0$. On en d\'eduit que $\gamma_1 \times \gamma_2$ ne	contient
pas de couple
singulier.
\end{proof}


\begin{proof}[\proofname\ du Th\'eor\`eme~\ref{theoBbis}]
Soient $x$ et $y$ deux points du plan, et $d=d_h(x,y)$. On suppose  $d
\geq 2$, sans quoi il n'y a rien  \`a montrer. Soit $i$ un entier
entre $1$ et $d-1$.
Le Th\'eor\`eme~\ref{theoAbis} nous dit que
$(F_{d-i}(y,x),F_i(x,y))$ est une composante de Reeb pour $(x,y)$, et
en particulier que  le produit $F_{d-i}(y,x) \times F_i(x,y)$
 est singulier pour $h$ ou bien pour $h^{-1}$.

Pour fixer les id\'ees, supposons  que $F_{d-i}(y,x) \times F_i(x,y)
\subset
\sing(h)$.  On consid\`ere
une d\'ecomposition minimale  $\gamma=\gamma_1 * \dots *
\gamma_d$ d'un arc g\'eod\'esique joignant $x$ \`a $y$.
 Le Lemme~\ref{lemm.ordre-crg} montre que dans ce cas, $\gamma_{i} \times
\gamma_{i+1} \cap \sing(h) \neq \emptyset$. D'autre part, il suit des
d\'efinitions que $M_\leftrightarrow(\gamma_1 * \dots *
\gamma_d)_i=M_\leftrightarrow(\gamma_{i} * \gamma_{i+1})$.
Le Lemme~\ref{lemm.fleche-singulier} dit
alors que $M_\leftrightarrow(\gamma_1 * \dots * \gamma_d)_i= \rightarrow$.
On a montr\'e que la $i$\`eme fl\`eche du mot ne d\'epend pas de
$\gamma$, ce qui  termine la preuve.
\end{proof}

Au passage, on a obtenu l'\'enonc\'e suivant (qui peut \^etre vu
comme une version ``technique'' du Th\'eor\`eme~\ref{theo.fleches-hori}).
\setcounter{theoter}{1}
\begin{theoter}\label{theoBter}
Pour tout $i=1, \dots, d-1$,
 \begin{alignat*}{3}
M_\leftrightarrow(\gamma_1 * \dots * \gamma_d)_i\,=\,\, \leftarrow  \  \ \
& \Leftrightarrow &\  \ \   &F_{d-i}(y,x) \times F_i(x,y)  \subset
\sing(h^{-1}) \\
M_\leftrightarrow(\gamma_1 * \dots * \gamma_d)_i\,=\,\, \rightarrow  \  \ \
& \Leftrightarrow &\  \ \   &F_{d-i}(y,x) \times F_i(x,y) \subset \sing(h).
\end{alignat*}
\end{theoter}

\section{Topologie de la suite des it\'er\'es d'un
disque libre}\label{sec.topologie-disque}
\emph{Cette partie est ind\'ependante des Sections~\ref{sec.distance} et
\ref{sec.crg}. Les d\'emonstrations de tous les \'enonc\'es se trouvent en
Section~\ref{ssec.preuves-disques}.}
Nous allons  \'etudier en d\'etail le comportement de la suite
des it\'er\'es positifs  d'un disque libre $D$. Le
premier r\'esultat est que cette suite converge
(Lemme~\ref{lemm.convergence}).
Bien s\^ur, ceci vaut \'egalement pour la suite des it\'er\'es n\'egatifs.
Les ensembles-limite positifs et n\'egatifs \'etant
disjoints et disjoints de $D$, leur r\'eunion d\'ecoupe le plan en trois
parties, une
partie ``n\'egative'', une partie ``centrale'' (contenant $D$) et une
partie ``positive'': on a ainsi associ\'e \`a $D$ une partition du plan
en cinq
morceaux (Lemme~\ref{lemm.topologie-disque2}). Un r\'esultat
technique essentiel compl\`ete la description; il affirme qu'on ne
peut pas aller contin\^ument de la partie centrale \`a l'un des deux
ensembles-limite sans rencontrer une infinit\'e d'it\'er\'es du
disque $D$
(Proposition~\ref{prop.non-accessibilite}, \cf
\figref{fig.iteres-disque}).
 Enfin, en faisant tendre le diam\`etre du disque $D$ vers $0$, on
localise ces constructions en un point $z$, ce qui permet d'associer
au point $z$ une nouvelle partition en cinq morceaux
(Lemme~\ref{lemm.partition}).
Au passage, le Lemme~\ref{lemm.singulier-limite} explique le lien avec
l'ensemble singulier de $h$: un couple $(x,y)$ est   singulier si
et seulement si $y$ appartient \`a l'ensemble-limite positif de
$x$. Ainsi, les r\'esultats de cette section permettent de construire
des ``grands'' ensembles de couples singuliers, et nous les utiliserons
dans la
Section~\ref{sec.preuve-crg} pour prouver l'existence des composantes
de Reeb
requises par le Th\'eor\`eme~\ref{theo.crgm}.

\subsection{Ensembles-limite}\label{sec.ensembles-limite}

On note $\cal K$ l'espace des compacts de la sph\`ere $\bbR^2 \cup
\{\infty\}$, muni de la topologie de Hausdorff (voir \cite{kura1},
chap. IV, paragraphes 38 et 42; une d\'efinition
possible de la convergence dans cet espace est donn\'ee dans la preuve
du lemme ci-dessous).
\begin{lemm}\label{lemm.convergence}
Soit $D$ un disque topologique ferm\'e libre; alors la suite des
it\'er\'es $(h^n(D))_{n \geq
0}$ converge dans $\cal K$.\footnote{Remarquons qu'on peut en d\'eduire
que ce r\'esultat est aussi vrai pour
un disque qui n'est pas libre: en effet, un tel disque est r\'eunion
d'un nombre fini de disques libres. Le r\'esultat n'est pas connu
pour un compact connexe quelconque (voir aussi le
Lemme~\ref{lemm.limitesup}).}
\end{lemm}
 Rappelons que l'espace $\cal K$ est m\'etrique,
compact, et que l'ensemble des \'el\'ements de $\cal K$ qui sont
connexes est encore compact.
La limite de cette suite est donc un compact connexe de la
sph\`ere $\bbR^2\cup\{\infty\}$, qui contient le point $\infty$ puisque
les orbites des points tendent vers $\infty$. On note $\lim^+D$
l'intersection de la limite avec le plan $\bbR^2$.
 On d\'efinit de m\^eme l'ensemble $\lim^-D$ comme la trace dans le
plan de la limite de la suite $(h^n(D))_{n \leq 0}$. Ces deux
ensembles sont des ferm\'es du plan, on les appellera
\emph{ensembles-limite} du disque $D$.
\begin{affi}\label{affi.disjointude}
 Les trois ensembles $D$, $\lim^+D$ et $\lim^-D$ sont deux \`a deux
 disjoints.
\end{affi}

\subsection{Partition du plan}
L'Affirmation~\ref{affi.disjointude} permet de d\'efinir les ensembles
suivants.
\begin{defi}[\rm(\figref{fig.iteres-disque})]
Soit $D$ un disque topologique ferm\'e libre.
\begin{itemize}
\item $U(D)$ est la composante connexe de  $\bbR^2 \setminus (\lim^+D
\cup \lim^-D)$ contenant $D$;
\item $U^-(D)$ est l'union des composantes connexes de $\bbR^2
\setminus \adhe(U(D))$ dont la fronti\`ere rencontre $\lim^-D$;
\item  $U^+(D)$ est l'union de celles dont la fronti\`ere rencontre
$\lim^+D$.
\end{itemize}
\end{defi}

\begin{lemm}[Partition associ\'ee \`a un disque libre]
\label{lemm.topologie-disque2}~
Soit $D$ un disque top\-o\-log\-ique ferm\'e libre. On a une partition
du plan
$$
\bbR^2 =  U^-(D) \sqcup \lim{}^-D \sqcup U(D)  \sqcup \lim{}^+D \sqcup
U^+(D),
$$
et chacun des ensembles de cette partition est invariant par $h$. De
plus,
$\partial U(D)=\lim{}^+D \cup \lim{}^-D.$
\end{lemm}

\subsection{Disposition des it\'er\'es de $D$}

\begin{figure}[htp]\anchor{fig.iteres-disque}
\centerline{\hbox{\input{fig-iteres-disque.pstex_t}}}
\caption{\label{fig.iteres-disque}Topologie de la suite des it\'er\'es
de $D$}
\end{figure}
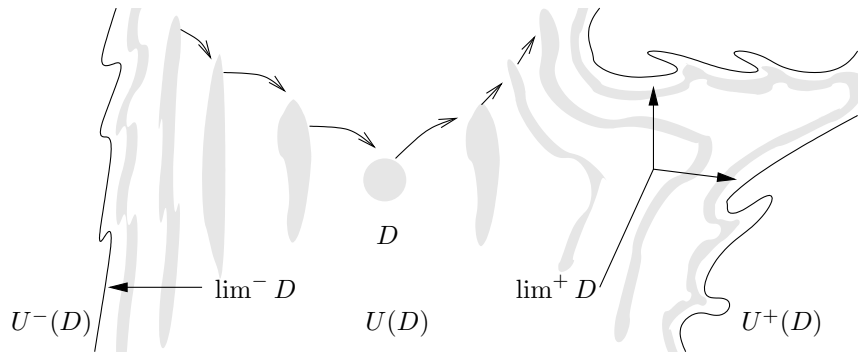

\begin{prop}[Non-accessibilit\'e]
\label{prop.non-accessibilite}
Soit $D$ un disque topologique ferm\'e libre.
Soit $\alpha$ un arc joignant un point de $U(D)$ \`a un point de
$\lim{}^+D \cup \lim{}^-D$.
 Alors il existe un entier $n$ tel que $\alpha$ rencontre  l'it\'er\'e
 $h^n(D)$.
\end{prop}
 Ce r\'esultat peut \^etre consid\'er\'e comme le lemme ``technique''
 fondamental
de notre \'etude. Il  est d\'ej\`a pr\'esent chez
Homma--Terasaka et Andrea (voir \cite{homm1} lemma 7; \cite{andr1},
proposition 1.7); cependant, le texte \cite{homm1} est tr\`es
difficile, et le texte \cite{andr1}  concerne les
ensembles-limite de trajectoires.  Nous en
donnerons une preuve fortement inspir\'ee du texte d'Andrea.
Une autre formulation consiste \`a dire que les points de $\lim^+ D
\cup \lim^- D$ ne sont pas accessibles depuis l'ouvert $U(D) \setminus
(\cup_{n \in \bbZ} h^n(D))$.
On peut aussi en d\'eduire un \'enonc\'e plus fort: sous les m\^emes
hypoth\`eses, l'arc $\alpha$ rencontre en fait une infinit\'e
d'it\'er\'es de $D$. En effet, dans le cas contraire,  $\alpha$
contiendrait un sous-arc $\alpha'$, disjoint de tous les
it\'er\'es de $D$, et joignant toujours un point de $U(D)$ \`a un
point de $\lim{}^+D \cup \lim{}^-D$; ceci contredirait le lemme.
La situation est symbolis\'ee par le dessin de la
\figref{fig.iteres-disque}.

\subsection{Partition associ\'ee \`a un point}
 Pour tout point $z$, on pose:
$$
\lim{}^+ z=\bigcap\{ \lim{}^+ D : D \mbox{ disque topologique libre
tel que } z \in \inte(D)\},
$$
$$
\lim{}^- z=\bigcap\{ \lim{}^- D  : D \mbox{ disque topologique libre
tel que } z \in \inte(D)\}.
$$
Ces deux ensembles sont appel\'es \emph{ensembles-limite} du point $z$.
\begin{lemm}
\label{lemm.singulier-limite}
Soit $(x,y)$ un couple de points du plan. Les assertions suivantes
sont \'equivalentes.
\begin{enumerate}
\item Le couple  $(x,y)$ est singulier;
\item $y \in \lim^+ x$;
\item $x \in \lim^- y$.
\end{enumerate}
\end{lemm}
La preuve de ce lemme est imm\'ediate. Une remarque en passant:
on d\'eduit du Lemme~\ref{lemm.convergence} que les ensembles-limite
d\'efinis pour les hom\'eomorphismes $h$ et $h^n$ ($n$ entier non
nul) co\"\i ncident; le Lemme~\ref{lemm.singulier-limite} montre
alors que $h$ et $h^n$ ont les m\^emes ensembles singuliers.
\begin{defi}
Comme dans le cas des disques, on d\'efinit les ensembles suivants.
\begin{itemize}
\item $U(z)$ est la composante connexe de $\bbR^2 \setminus (\lim^+z \cup
\lim^-z)$ contenant $z$;
\item $U^-(z)$ est l'union des composantes connexes de $\bbR^2
\setminus \adhe(U(z))$ dont la fronti\`ere rencontre  $\lim^-z$;
\item $U^+(z)$ est l'union de celles dont la fronti\`ere rencontre
$\lim^+z$.
\end{itemize}
\end{defi}
\begin{lemm}[Partition associ\'ee \`a un point]
\label{lemm.partition}
Soit $z$ un point du plan. On a une partition du plan
$$
\bbR^2 =  U^-(z) \sqcup \lim{}^-z \sqcup U(z)  \sqcup \lim{}^+z \sqcup
U^+(z),
$$
o\`u chacun de ces ensembles est invariant par $h$.
De plus,
$\partial U(z)=\lim{}^+z \cup \lim{}^-z$,
$\partial U^-(z) \subset \lim{}^-z$,
$\partial U^+(z) \subset \lim{}^+z$.
\end{lemm}

\begin{lemm}[Liens point/disques]
\label{lemm.liens}
On a
\begin{gather*}
U(z)=\bigcup\{U(D)\}\\
U^+(z)=\bigcap\{U^+(D)\}\\
U^-(z)=\bigcap\{U^-(D)\}
\end{gather*}
o\`u les unions et intersections ont lieu sur tous les disques
topologiques ferm\'es libres $D$ dont l'int\'erieur contient $z$.
\end{lemm}

\subsection{Preuves}\label{ssec.preuves-disques}

\begin{proof}[\proofname\ du Lemme~\ref{lemm.convergence}]
Soit $D$ un disque topologique ferm\'e libre.
Notons $S$ et $I$ les \emph{limites sup\'erieures et inf\'erieures}
de la suite $(h^n(D))_{n \geq 0}$:
\begin{align*}
S & =  \{x \in \bbR^2:\forall V \mbox{ voisinage de } x, \forall n_0 \geq
0, \exists
n \geq n_0 \mbox{ tel que } h^n(D) \cap V \neq \emptyset\}\\
I & =  \{x \in \bbR^2:\forall V \mbox{ voisinage de } x, \exists n_0 \geq
0, \forall
n \geq n_0, h^n(D) \cap V \neq \emptyset\}
\end{align*}
Prouver la convergence revient \`a montrer que $S=I$. Bien s\^ur, $I
\subset S$; montrons l'inclusion r\'eciproque. Soit $x$ un point de
$S$, et $V$ un disque ferm\'e libre voisinage de $x$. Soit $n_0$
un entier tel que $h^{n_0}(D) \cap \inte(V) \neq \emptyset$ (et donc
aussi  $h^{n_0}(\inte(D)) \cap \inte(V) \neq \emptyset$). D'apr\`es
le lemme de Franks (dans la version du Lemme~\ref{lemm.intervalle}),
l'ensemble des entiers $n$ tels
que $h^n(\inte(D))$ rencontre $\inte(V)$ est un intervalle de $\bbZ$;
mais il
contient $n_0$ et une infinit\'e d'entiers sup\'erieurs \`a $n_0$ (par
d\'efinition de $S$), donc il contient tous les entiers sup\'erieurs
\`a $n_0$. Autrement dit, $x$ est un point de $I$.
\end{proof}

\begin{proof}[\proofname\ de l'Affirmation~\ref{affi.disjointude}]
Comme $D$ est ferm\'e et libre, on peut trouver un autre disque
topologique ferm\'e libre $V$ qui contient $D$ dans son int\'erieur
(en utilisant le th\'eor\`eme de Schoenflies).
 D'apr\`es le Corollaire~\ref{coro.liberte} sur les disques libres,
ce disque $V$ est disjoint
de tous les it\'er\'es $h^n(\inte(V))$ pour $n \neq 0$, donc aussi des
it\'er\'es $h^n(D)$. Par suite, les
ensembles $\lim^+D$ et $\lim^-D$ ne rencontrent pas l'int\'erieur de
$V$, et sont donc disjoints de $D$.

On en d\'eduit que tout point $x$ de $\lim^+(D)$ poss\`ede un voisinage
$V'$
disque topologique ferm\'e libre disjoint de $D$, et un raisonnement
analogue montre que l'ensemble $\lim^-(D)$ ne rencontre pas $V'$. Par
cons\'equent, les ensembles $\lim^+(D)$  et  $\lim^-(D)$  sont disjoints.
\end{proof}

\begin{proof}[\proofname\ du Lemme~\ref{lemm.topologie-disque2}]
On se donne un disque topologique ferm\'e libre $D$.
\paragraph{Invariance}
Il est clair que les ensembles $\lim^+D$ et $\lim^-D$ sont invariants
par $h$. Montrons qu'il en est de m\^eme pour l'ouvert $U(D)$. En tant
que composante connexe d'un ensemble invariant, il est invariant ou libre.
D'apr\`es le Lemme~\ref{lemm.ETDS} (compacts connexes libres), il
existe un domaine de translation $\Delta$ contenant $D$. Comme $h$ est
conjugu\'e \`a une translation sur $\Delta$, les ensembles $\lim^+D$
et $\lim^-D$ sont disjoints de $\Delta$. Puisque $\Delta$ est connexe,
on en d\'eduit qu'il est inclus dans $U(D)$. D'autre part $\Delta$
est invariant par $h$, ce qui prouve que $U(D)$ n'est pas libre, donc
est invariant.
On en d\'eduit imm\'ediatement l'invariance des ensembles $U^+(D)$
et $U^-(D)$.

\paragraph{Fronti\`ere}
La fronti\`ere de l'ensemble $U(D)$ est
 incluse dans $\lim^+D\cup \lim^ -D$: ceci suit d'un r\'esultat
de topologie g\'en\'erale
 (Proposition~\ref{prop.bumpwire} de
l'Appendice~\ref{appe.topo}).
Puisque l'ensemble $U(D)$ est invariant, il contient tous les it\'er\'es
de $D$, ce
qui entra\^\i ne l'inclusion des deux ensembles $\lim^+D$ et
$\lim^-D$ dans $\adhe(U(D))$. Par ailleurs, par d\'efinition, $U(D)$ est
disjoint de  $\lim^+D$ et $\lim^-D$; donc la fronti\`ere de
$U(D)$ contient ces deux ensembles. On en d\'eduit l'\'egalit\'e $\partial
U(D)=\lim{}^+D \cup \lim{}^-D$.

\paragraph{Partition}
D'apr\`es la Proposition~\ref{prop.bumpwire}, la fronti\`ere de
toute composante connexe de $\bbR^2 \setminus \adhe(U(D))$ est non vide
et contenue dans
$\partial U(D)=\lim{}^+D \cup \lim{}^-D$. Par cons\'equent,  la
d\'efinition des
ensembles $U^+(D)$ et $U^-(D)$ implique que leur union recouvre le
compl\'ementaire de $\adhe(U(D))$. D'apr\`es le paragraphe sur la
fronti\`ere,  $\adhe(U(D))=U(D) \cup
\lim^+D \cup \lim^-D$, et par cons\'equent  les cinq
ensembles consid\'er\'es recouvrent le plan.
Il reste \`a voir que ces cinq ensembles sont disjoints deux \`a deux.
On a
\begin{itemize}
\item $U(D) \cap
\lim^+D =\emptyset$, par d\'efinition de $U(D)$;
\item	$U(D) \cap
U^+(D) =\emptyset$,  par d\'efinition de $U^+(D)$;
\item de m\^eme,  $U(D) \cap \lim^-D =\emptyset$ et  $U(D) \cap
U^-(D) =\emptyset$;
\item les ensembles $U^+(D)$ et
$U^-(D)$ sont chacun disjoints des ensembles $\lim^+(D)$ et
$\lim^-(D)$, d'apr\`es l'\'egalit\'e  $\partial
U(D)=\lim{}^+D \cup \lim{}^-D$;
\item  $\lim^+D \cap \lim^-D=\emptyset$, d'apr\`es
l'Affirmation~\ref{affi.disjointude}.
\end{itemize}
Il reste \`a voir que $U^+(D)$ et $U^-(D)$ sont disjoints.
Soit  $C$ une composante connexe du compl\'ementaire de
$\adhe(U(D))$; il s'agit de voir que la fronti\`ere de $C$ ne peut
pas rencontrer \`a la fois  $\lim^- D$ et $\lim^+
D$. D'apr\`es le Th\'eor\`eme~\ref{theo.connexite} de l'appendice,
$\partial C$ est connexe; et d'apr\`es la
Proposition~\ref{prop.bumpwire}, $\partial C$ est inclus dans $\partial
\adhe(U(D))$, \cad\ dans  $\lim^- D \cup \lim ^+ D$. Comme les deux
ensembles   $\lim^- D$ et $\lim^+ D$ sont des ferm\'es disjoints,
$\partial C$ est enti\`erement inclus dans l'un d'eux, ce que l'on
voulait.
\end{proof}

\begin{proof}[\proofname\ de la Proposition~\ref{prop.non-accessibilite}]
On se donne un disque topologique ferm\'e libre $D$.
\paragraph{Cas simple}
Commen\c{c}ons par montrer que \emph{tout arc $\alpha$,
\emph{libre}, joignant un point de $h^{n_0}(D)$ \`a un point de
$\lim^+D$, doit rencontrer tous
les it\'er\'es $h^n(D)$ pour $n \geq n_0$}. Supposons le contraire.
En utilisant le th\'eor\`eme de Schoenflies, on construit alors un disque
topologique ferm\'e $D'$, libre, contenant $\alpha$ dans son
int\'erieur, disjoint d'au moins un it\'er\'e $h^n(D)$ pour $n >
n_0$. Par d\'efinition de l'ensemble  $\lim^+D$, l'int\'erieur du disque
$D'$ rencontre des it\'er\'es $h^{n'}(D)$ pour $n'>n$;
 d'autre part, il rencontre aussi $h^{n_0}(D)$. L'ensemble des entiers
$k$ tels que $h^k(D)$ rencontre $\inte(D')$ n'est pas un intervalle,
ce qui contredit le lemme de Franks (\ref{lemm.intervalle}).
On a bien s\^ur un \'enonc\'e \'equivalent pour $\lim^-D$.

\paragraph{Type des points voisins}
Pour tout point $x$ du plan, on note  $B^\mathrm{eucl}_r(x)$ le disque
ouvert de centre $x$ et de rayon $r$ pour la distance euclidienne sur
le plan.
 On pose alors
$$
\phi(x)=\sup\{r : B^\mathrm{eucl}_r(x) \mbox{ est libre} \}.
$$
La fonction $\phi$ est strictement positive et continue sur le plan,
et pour tout $x$, le disque euclidien ouvert
$B^\mathrm{eucl}_{\phi(x)}(x)$ est libre.
On pose alors
$$
V=\{x \in U(D) : B^\mathrm{eucl}_{\phi(x)}(x) \cap \partial U(D) \neq
\emptyset\}.
$$
Cet ensemble est  appel\'e \emph{ensemble des points voisins de
$\partial U(D)$}. Il s'agit d'un ensemble ouvert, et $V \cup \partial
U(D)$ est un voisinage de $\partial U(D)$ dans
$\adhe(U(D))$.

Soit $x$ un point de $V$. On consid\`ere  un arc $\alpha$ joignant $x$
\`a un point
de $\partial U(D)$, inclus dans $B^\mathrm{eucl}_{\phi(x)}(x)$ et
d'int\'erieur inclus dans $U(D)$ (la d\'efinition de $V$ permet de
trouver un tel arc). On dira que $x$ est \emph{de
type $D$} si $\alpha$ rencontre un it\'er\'e de $D$, \emph{de type
$\partial U(D)$} sinon (cf \figref{fig.type1}). L'int\'er\^et de
cette d\'efinition est
qu'elle ne d\'epend pas du choix de $\alpha$. En effet, soit $\alpha'$
 un autre  arc comme ci-dessus, et supposons que $\alpha$ rencontre un
it\'er\'e de $D$ tandis que $\alpha'$ n'en rencontre pas (cf
\figref{fig.type2}). Alors
$\alpha \cup \alpha'$ est libre, et contient un arc $\beta$
rencontrant $\partial U(D)$ et un unique it\'er\'e de
$D$. L'existence de $\beta$ contredit le ``cas simple'' \'etudi\'e au
d\'ebut de cette d\'emonstration. Nous allons voir qu'en r\'ealit\'e,
tous les points voisins
sont de type $D$.

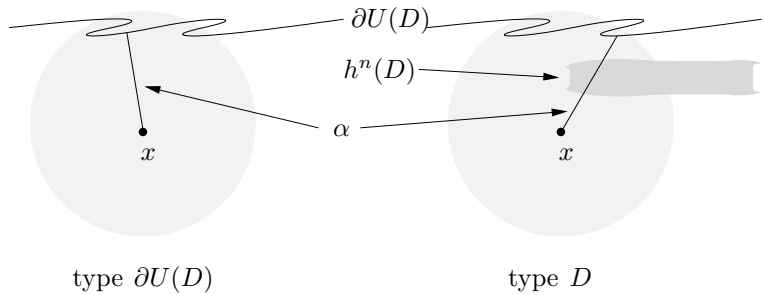
\begin{figure}[htp]\anchor{fig.type1}
\centerline{\hbox{\input{fig-type1.pstex_t}}}
\caption{\label{fig.type1}Type des points voisins}
\end{figure}
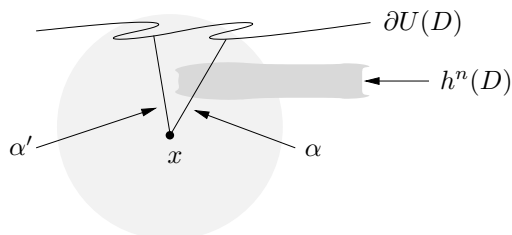
\begin{figure}[htp]\anchor{fig.type2}
\centerline{\hbox{\input{fig-type2.pstex_t}}}
\caption{\label{fig.type2}Le type est bien d\'efini}
\end{figure}

On voit facilement que le type des points de $V$ est localement
constant, donc constant sur les composantes connexes de
$V$. En particulier, si  une composante connexe de $V$ rencontre
un it\'er\'e de $D$, alors tous les points de cette composante connexe
 sont de type $D$.

\paragraph{Topologie plane}

Soit maintenant $\alpha$ un arc v\'erifiant les hypoth\`eses de la
Proposition~\ref{prop.non-accessibilite}, c'est-\`a-dire joignant un
point de $U(D)$ \`a
un point de $\partial U(D)$. Quitte \`a raccourcir $\alpha$, on peut
supposer que l'int\'erieur de $\alpha$ est inclus dans $U(D)$; on
appelle alors $z$ l'extr\'emit\'e de $\alpha$ sur $\partial U(D)$.
On applique l'affirmation suivante (dont la preuve est donn\'ee plus bas).
\begin{affi}\label{affi.ouvert-connexe}
Pour tout point $z$ de $\partial U(D)$,
il existe une composante connexe $V'$ de $V$ telle que $V' \cup \{z\}$
est un voisinage de $z$ dans $U(D)\cup\{z\}$.
\end{affi}

Puisque $z$ appartient \`a $\lim^+D \cup \lim^-D$, et que les it\'er\'es
de $D$ sont dans $U(D)$ (Lemme~\ref{lemm.topologie-disque2}),
il existe un it\'er\'e de $D$ qui rencontre $V'$.
Par cons\'equent,  les points de $V'$ sont tous de type $D$.
D'autre part, quitte \`a raccourcir $\alpha$, on peut supposer que
l'extr\'emit\'e $x$ de $\alpha$ dans $U(D)$ est un point de $V'$, et
que $\alpha$ est inclus dans le disque $B^\mathrm{eucl}_{\phi(x)}(x)$.
Le point $x$ est de type $D$, donc l'arc $\alpha$ rencontre un
it\'er\'e de $D$, ce que l'on voulait.
\end{proof}

\begin{proof}[\proofname\ de l'Affirmation~\ref{affi.ouvert-connexe}]
Il s'agit seulement de con\-sid\-\'e\-ra\-tions de top\-o\-log\-ie plane.
On pose $W=V \cup (\bbR^2 \setminus U(D))$; c'est une partie ouverte
du plan, contenant $\partial U(D)$,
et dont la fronti\`ere est incluse dans $U(D)$. Soit $W'$ la composante
connexe de $W$ contenant $z$; c'est un ensemble ouvert et connexe.
Soit  $V'=W' \cap U(D)$: c'est un ensemble inclus dans $V$, et $V'
\cup \{z\}$ est un voisinage de $z$ dans $U(D)\cup\{z\}$. Pour prouver
l'affirmation, il suffit donc de montrer que $V'$ est connexe.

 Ceci va d\'ecouler du lemme d'Alexander
(Th\'eor\`eme~\ref{theo.lemme-alex} de l'appendice).
On appelle $U_1$ la r\'eunion de $U(D)$ et de toutes les
composantes connexes de $W$ distinctes de $W'$. C'est un ensemble
ouvert, et comme toute	composante connexe de $W$ rencontre $U(D)$
(Proposition~\ref{prop.bumpwire}) qui
est connexe, c'est un ensemble connexe. L'ensemble $U_2=W'$ est
\'egalement un ouvert connexe, et $U_1 \cup U_2=U(D) \cup
W=\bbR^2$. Le lemme d'Alexander affirme alors que l'intersection $U_1
\cap U_2$ est connexe, et cette intersection n'est autre que
l'ensemble $V'$.
\end{proof}

\begin{proof}[\proofname\ du Lemme~\ref{lemm.partition}]
Dans toute la preuve, $D$ d\'esigne un disque topologique ferm\'e libre
dont l'int\'erieur
contient $z$.

Par d\'efinition,  $\lim^+D$ contient $\lim^+z$, donc
l'ouvert connexe $U(D)$ est disjoint de $\lim^+z$ (et de m\^eme il est
disjoint de $\lim^-z$), on en d\'eduit
l'inclusion $U(D) \subset U(z)$. Il s'ensuit que
$$
\lim{}^+z \subset \lim{}^+D \subset \adhe(U(D)) \subset \adhe(U(z)).
$$
Il est clair que  les ensembles $\lim^+z$ et $\lim^-z$	sont invariants
par $h$. Puisque l'ensemble $U(D)$ est invariant et inclus dans l'ensemble
$U(z)$,
celui-ci est aussi invariant. On en d\'eduit l'invariance de $U^+(z)$
et  $U^-(z)$.

D'apr\`es la Proposition~\ref{prop.bumpwire} et la d\'efinition de
 $U(z)$, on a $\partial U(z) \subset \lim^+z \cup \lim^-z$.
 D'autre part, l'inclusion $\lim^+z \subset \adhe(U(z))$ a \'et\'e
 montr\'ee dans le paragraphe pr\'eliminaire; on a bien s\^ur l'inclusion
 analogue
pour les limites n\'egatives.
L'\'egalit\'e $\partial U(z) = \lim^+z \cup \lim^-z$ s'ensuit.

Si $U$ est une composante connexe du compl\'ementaire de
$\adhe(U(z))$, on montre \`a l'aide de la
Proposition~\ref{prop.bumpwire} et du Th\'eor\`eme~\ref{theo.connexite}
que la
fronti\`ere de $U$ est incluse dans l'un des deux ensembles $\lim^-z$
et $\lim+z$ (exactement comme dans l'argument final de la preuve du
Lemme~\ref{lemm.topologie-disque2}). Ceci permet de voir que
$\partial U^-(z) \subset \lim^-z$ et $\partial U^+(z) \subset \lim^+z$.

De mani\`ere tout \`a fait analogue au cas des disques
(Lemme~\ref{lemm.topologie-disque2}), on
montre que les cinq ensembles recouvrent le plan, puis qu'ils sont
deux \`a deux disjoints.
\end{proof}

\begin{proof}[\proofname\ du Lemme~\ref{lemm.liens}]
Ici encore, $D$ d\'esigne un disque topologique ferm\'e libre dont
l'int\'erieur
contient $z$.
Montrons la premi\`ere des trois \'egalit\'es.
Soit $(D_k)_{k \geq 0}$ une suite d\'ecroissante de disques
topologiques ferm\'es libres, dont les int\'erieurs contiennent $z$,
et dont l'intersection est r\'eduite \`a $\{z\}$.
Il est clair que
l'application $D \mapsto \lim^+D$ est croissante (par rapport \`a
l'inclusion). D'apr\`es la d\'efinition de $U(D)$, l'application $D
\mapsto U(D)$ est alors d\'ecroissante.
 On en d\'eduit que
$$\lim{}^+z=\bigcap_{k \geq 0} \lim{}^+D_k \quad\text{et que}\quad
\bigcup\{U(D)\}= \bigcup_{k\geq 0}\{U(D_k)\}.$$
Soit maintenant un point  $y$ dans $U(z)$, c'est-\`a-dire tel que
$\lim^+z \cup \lim^-z$ ne s\'epare pas $y$ et $z$. Il existe alors un
arc $\gamma$
joignant $y$ \`a $z$ dans le compl\'ementaire de $\lim^+z \cup \lim^-z$;
puisque l'intersection des compacts $\cap_k\{(\lim^+D_k \cup \lim^-D_k)
\cap \gamma,
k \geq 0\}$ est vide, c'est
qu'il existe un disque $D_k$ tel que $(\lim^+D_k \cup \lim^-D_k) \cap
\gamma$ est
vide. On en d\'eduit que $\gamma$ est enti\`erement inclus dans
$U(D_k)$. Ceci prouve l'inclusion $U(z) \subset \cup \{U(D_k),k \geq 0\}$.
L'inclusion r\'eciproque vient de la propri\'et\'e
 $U(D) \subset U(z)$ vue dans le paragraphe pr\'eliminaire de la preuve
du Lemme~\ref{lemm.partition}.

Montrons maintenant que $U^+(z) \subset \cap\{U^+(D)\}$.
Pour cela, on consid\`ere une composante connexe $U$ de $\bbR^2
\setminus \adhe(U(z))$ dont l'adh\'erence rencontre $\lim^+z$. Pour
tout $D$, l'ensemble $\lim^+D$ contient $\lim^+z$, il rencontre donc
l'adh\'erence de $U$; mais par ailleurs l'ensemble  $\adhe(U(D))$ est
inclus dans
$\adhe(U(z))$ (voir le paragraphe pr\'eliminaire de la preuve
pr\'ec\'edente), il est donc disjoint de $U$. Ceci prouve que $U$
est inclus dans $U^+(D)$, ce que l'on cherchait.
De la m\^eme mani\`ere, on obtient l'inclusion	$U^-(z) \subset
\cap\{U^-(D)\}$.

Montrons les inclusions r\'eciproques. Soit $y$ un point de
$\cap\{U^+(D)\}$. Ce point n'est dans aucun des ensembles $U(D)$, donc
n'appartient pas \`a $U(z)$ (d'apr\`es la premi\`ere
\'egalit\'e). De m\^eme, $y$ n'appartient pas \`a $\lim^+(z)$ ni
\`a  $\lim^-(z)$. De m\^eme, puisque  $U^-(z) \subset \cap\{U^-(D)\}$, $y$
n'est pas dans $U^-(z)$. D'apr\`es le Lemme~\ref{lemm.partition}, $y$
appartient \`a $U^+(z)$. Ceci prouve que $ \cap\{U^+(D)\} \subset
U^+(z)$. L'inclusion pour les ensembles n\'egatifs est analogue.
\end{proof}

\section{Preuve du Th\'eor\`eme~\ref{theo.crgm}: existence des com\-pos\-antes
de Reeb }
\label{sec.preuve-crg}
Dans cette partie, nous utilisons les r\'esultats de la section
pr\'ec\'edente pour d\'emontrer le Th\'eor\`eme~\ref{theoAbis}.
Malheureusement, l'auteur n'a pas r\'eussi \`a trouver
une preuve qui ne brise pas la sym\'etrie des r\^oles de $x$ et de
$y$. Il s'en
excuse aupr\`es du lecteur, et incite celui-ci \`a chercher une
preuve sym\'etrique!
Pour montrer que les couples $(F_{d-i}(y,x),F_i(x,y))$ sont des
composantes de Reeb, la vraie difficult\'e
consistera \`a montrer que le produit de ces deux ensembles est
singulier. En effet, la d\'efinition de l'ensemble singulier
entra\^\i nera ensuite que les deux ensembles sont disjoints;  on en
d\'eduira facilement que les propri\'et\'es de s\'eparation (2 et 3) de la
d\'efinition des composantes de Reeb  sont
satisfaites.

Essayons de donner une id\'ee de la preuve du
Th\'eor\`eme~\ref{theoAbis}.
Pour cela, imaginons un point mobile $z$ qui s'\'eloigne d'un point $x$
immobile; le point $z$ sort successivement des $h$--boules de centre
$x$ et de rayon $1$, $2$, $3$, \textit{etc.}.
Dans la premi\`ere \'etape, nous montrons
 qu'au moment o\`u  la distance de translation de $z$ \`a
$x$ augmente d'une unit\'e (autrement dit si $z$ est sur la
 fronti\`ere d'une $h$--boule de centre $x$), ou bien l'union des deux
 ensembles-limite de $z$
contient $x$, ou bien elle s\'epare $x$ de $z$
(Affirmation~\ref{affi.quatre-possibilites}).

\begin{figure}[htp]\anchor{fig.idee-crg}
\centerline{\hbox{\input{fig-idee-crg.pstex_t}}}
\caption{\label{fig.idee-crg}Existence des
composantes de Reeb: id\'ee de la preuve}
\end{figure}
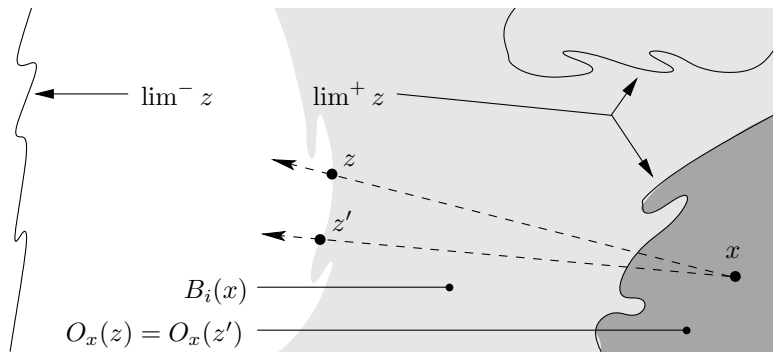

Pour simplifier, oublions la premi\`ere des deux possibilit\'es,
et plaçons-nous dans le cas (g\'en\'erique) o\`u l'ensemble $\lim^-z
\cup \lim^+z$
 s\'epare $x$ de $z$. Notons
$O_x(z)$ la composante connexe du compl\'ementaire de $\lim^-z \cup
\lim^+z$ qui contient $x$ (\figref{fig.idee-crg}). La fronti\`ere
de $O_x(z)$
rencontre un et un seul des deux ensembles-limite de $z$, et on dira
selon les cas
 que $z$ est de type positif ou n\'egatif par rapport \`a $x$.
La deuxi\`eme \'etape de la preuve consiste \`a montrer que ce que l'on
vient de
faire est ``localement constant''. De mani\`ere un peu plus
pr\'ecise, consid\'erons un autre point mobile $z'$ qui s'\'eloigne
\'egalement de $x$; on
montre que si $z$ et $z'$ sortent tous deux d'une m\^eme $h$--boule
de centre
$x$, et qu'ils sont assez proches l'un de l'autre ($d_h(z,z') \leq
1$), alors leur type par rapport \`a $x$ est le m\^eme. Mieux: les
ensembles $O_x(z)$ et $O_x(z')$ co\"\i ncident
(Affirmation~\ref{affi.invariance}).

Pour la troisi\`eme \'etape, on se donne deux point $x$ et $y$, et
on consid\`ere l'un des ensembles $F_i(x,y)$ candidats \`a former
les bords des composantes de Reeb. Cet ensemble est inclus dans la
fronti\`ere d'une $h$--boule de centre $x$, on peut donc appliquer le
r\'esultat de la premi\`ere \'etape; de plus, il est connexe, et
l'invariance locale obtenue \`a la deuxi\`eme \'etape  entra\^\i ne
que tous les points $z$ de $F_i(x,y)$ ont
le m\^eme type par rapport \`a $x$, et d\'efinissent le m\^eme ensemble
$O_x(z)$. On
consid\`ere alors l'ensemble $F$ fronti\`ere de $O_x(z)$, et on
montre que le couple $(F,F_i(x,y))$ forme une composante de Reeb
(Affirmation~\ref{affi.existence-crg}).

Ind\'ependamment de tout ce qui pr\'ec\`ede, on v\'erifie dans la
quatri\`eme \'etape que toute composante de Reeb pour $(x,y)$
``contient'' l'un des couples $(F_{d-i}(y,x),\allowbreak F_i(x,y))$
(Affirmation~\ref{affi.crg-minimales}). On en
d\'eduira d'abord que l'ensemble $F$ de la troisi\`eme \'etape
contient $F_{d-i}(y,x)$, et donc que les couples
$(F_{d-i}(y,x),F_i(x,y))$ sont bien des
composantes de Reeb (on retrouve ainsi la sym\'etrie perdue). On en
d\'eduira ensuite qu'il n'existe pas d'autre composante minimale.

\subsection{\'Enonc\'es des \'etapes de la preuve}
On donne ici les \'enonc\'es pr\'ecis des diff\'erentes \'etapes. Les
preuves sont fournies plus bas.
La premi\`ere affirmation est une cons\'equence des r\'esultats sur
 la disposition des it\'er\'es d'un
 disque, et  surtout du lemme de non-accessibilit\'e.
\begin{affi}\label{affi.quatre-possibilites}
Soient $x$ et $z$ deux points du plan. Si $z$ est sur la fronti\`ere
d'une $h$--boule de centre $x$, alors $x \not \in U(z)$.\footnote{On
pourrait \'egalement prouver la r\'eciproque.}
\end{affi}

Soient $x \in \bbR^2$, $i \geq 1$, $z \in \partial B_i(x)$. D'apr\`es
l'affirmation pr\'ec\'edente et la partition du plan donn\'ee par le
Lemme~\ref{lemm.partition}, une et une seule des quatre possibilit\'es
suivantes est v\'erifi\'ee:
$$
x \in	U^-(z), \ \ x \in  \lim{}^-z,  \ \ x \in  \lim{}^+z, \ \ x
\in U^+(z)
$$
Selon les cas, on dira que $z$ est \emph{de type $U^-$, $\lim{}^-$,
$\lim{}^+$
ou  $U^+$ par rapport \`a $x$}. Si $z$ est de type $U^-$ ou $U^+$,
on notera $O_x(z)$ la
composante connexe de $U^-(z)$ ou de $U^+(z)$ qui
contient $x$.
Les liens entre la
 partition associ\'ee \`a un point et celle associ\'ee \`a un disque
libre le contenant vont nous permettre de montrer la deuxi\`eme
affirmation.
\begin{affi}[Invariance locale du type]\label{affi.invariance}
Soient $x$ un point du plan,  $i \geq 1$, et $z$ et $z'$ deux points
de $\partial B_i(x)$.
On suppose que $d_h(z,z') \leq 1$. Alors
\begin{itemize}
\item $z$ et $z'$ sont du m\^eme type par rapport \`a $x$;
\item si ce type commun est $U^+$ ou $U^-$, on a $O_x(z)=O_x(z')$.
\end{itemize}
\end{affi}
L'invariance locale du type et la connexit\'e des ensembles
$F_i(x,y)$ vont  entra\^iner la troisi\`eme affirmation.
\begin{affi}[Existence des composantes]
\label{affi.existence-crg}
Soient $x$ et $y$ deux points du plan, et $1 \leq i \leq d_h(x,y)-1$.
Il existe alors  une partie $F$ du plan telle que
le couple  $(F, F_i(x,y))$ est une composante de Reeb  pour
$(x,y)$.
\end{affi}
La quatri\`eme affirmation, elle, est ind\'ependante des trois
premi\`eres. Elle r\'esultera surtout des propri\'et\'es de la
distance de translation (in\'egalit\'e triangulaire et existence
d'arcs g\'eod\'esiques).
\begin{affi}[Unicit\'e des composantes]
\label{affi.crg-minimales}
Soient $x$ et $y$ deux points du plan, et $d=d_h(x,y)$. Soit $(F,G)$
une composante de Reeb	pour $(x,y)$. Alors
 il existe  un entier $i$, $1 \leq i \leq d-1$, tel que
$$
 F_{d-i}(y,x) \subset F \subset \partial B_{d-i}(y) \mbox{ et }
 F_i(x,y) \subset G \subset \partial B_i(x).
$$
\end{affi}
Dans le cas non d\'eg\'en\'er\'e, l'entier $i$ sera \'egal \`a
la distance de translation entre le point $x$ et l'ensemble $F$.

\subsection{Premi\`ere \'etape: lien entre partition et distance de
translation}

\begin{proof}[\proofname\ de l'Affirmation~\ref{affi.quatre-possibilites}]
Puisque, d'apr\`es le Lemme~\ref{lemm.liens}, $U(z)$ est
 l'un\-i\-on des ensembles $U(D)$ pour les disques
 topologiques ferm\'es libres $D$ dont l'int\'erieur contient $z$,
 il suffit de prouver
que si $D$ est un tel disque, $x
\not \in U(D)$.
 On raisonne par l'absurde: \emph{on suppose que
$x \in U(D)$ pour un tel disque $D$, et on va montrer que tous
 les points de $D$ sont \`a la
m\^eme $h$--distance de $x$}, ce qui contredira le fait que $z \in
\partial B_i(x)$.

Puisque $D$ est libre, tout point de $D$ est \`a $h$--distance $1$ de
 tout autre point de $D$, donc aussi de tout point d'un
it\'er\'e de $D$ (invariance par $h$ de la distance $d_h$,
Lemme~\ref{lemm.propriete-h-distance}). Le r\'esultat est
donc vrai si $x$ appartient \`a un it\'er\'e de $D$.
On suppose maintenant que $x$ n'appartient \`a aucun it\'er\'e de
$D$. Soit
$y$ un point de $D$; on va montrer que \emph{pour tout point $y'$ de
$D$, on
a $d_h(x,y') \leq d_h(x,y)$}. Par sym\'etrie des r\^oles de $y$ et
 $y'$, ceci conclura la preuve en montrant que tous les points de $D$
 sont \`a m\^eme $h$--distance de $x$.
Pour montrer cette in\'egalit\'e, on consid\`ere un arc g\'eod\'esique
$\gamma$ de $x$ \`a
$y$ (\figref{fig.calculs-distance}). Soit
$p$ le premier point rencontr\'e, en parcourant $\gamma$ depuis
l'extr\'emit\'e  $x$, sur le ferm\'e
$$
\adhe\left( \bigcup_{n \in \bbZ}h^n(D) \right) =  \lim{}^-(D) \cup
\bigcup_{n \in
\bbZ}h^n(D) \cup \lim{}^+(D).
$$
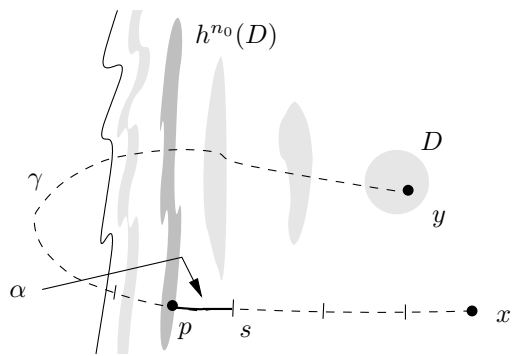
\begin{figure}[htp]\anchor{fig.calculs-distance}
\centerline{\hbox{\input{fig-calculs-distance2.pstex_t}}}
\caption{\label{fig.calculs-distance}Calculs de distance}
\end{figure}
Comme on a suppos\'e que $x$ appartient \`a $U(D)$, d'apr\`es la
proposition de non-accessibilit\'e
(\ref{prop.non-accessibilite}), le point $p$ appartient en fait  \`a un
it\'er\'e de $D$. Remarquons que $p$ est diff\'erent de $x$.
Soit $n_0$ tel que $p \in h^{n_0}(D)$. On se donne une
d\'ecomposition minimale de $\gamma$, et on note
 $s$ le dernier sommet de cette d\'ecomposition
situ\'e strictement avant $p$ lorsqu'on parcourt $\gamma$ en partant
du point $x$; on a donc:
\begin{itemize}
\item $d_h(x,s) \leq d_h(x,y)-1$;
\item l'arc $\alpha=[sp]_\gamma$ est libre.
\end{itemize}
D'autre part $h^{n_0}(D)$ est libre, et $\alpha$ ne rencontre aucun
autre it\'er\'e de $D$ (par d\'efinition de $p$). On en d\'eduit que
l'ensemble $\alpha \cup h^{n_0}(D)$ est encore libre. Soit $y'$ un
point quelconque de $h^{n_0}(D)$; il existe un arc $\beta$ reliant
$s$ \`a $y'$ et inclus dans $\alpha \cup h^{n_0}(D)$; en particulier,
$\beta$ est libre. D'apr\`es  l'Affirmation~\ref{affi.arcs-et-domaines},
 on a $d_h(s,y')=1$. Par invariance de la $h$--distance, ceci est
encore vrai si $y'$ est un point quelconque de $D$
(Lemme~\ref{lemm.propriete-h-distance}).
Un point $y'$ de $D$ v\'erifie donc
\begin{align*}
d_h(x,y') & \leq   d_h(x,s) +	d_h(s,y') \\
& \leq   d_h(x,y)-1 + 1 \\
& \leq  d_h(x,y).
\end{align*}
Ce que l'on voulait.
\end{proof}

\subsection{Deuxi\`eme \'etape: invariance locale du $x$--type}

\begin{proof}[\proofname\ de l'Affirmation~\ref{affi.invariance}]
 On se place sous les hypoth\`eses de l'aff\-ir\-ma\-tion. Puisque
 $d_h(z,z') \leq 1$, ou bien $z$ et $z'$ sont dans la m\^eme orbite de
 $h$, et la conclusion de l'affirmation est clairement satisfaite;
 ou bien il existe un arc libre joignant ces deux points, et (par
 \'epaississement), il existe aussi un disque topologique ferm\'e
 libre $D$
 dont l'int\'erieur contient $z$ et $z'$. Nous nous pla\c{c}ons dans
 ce cas.

\emph{Supposons d'abord que $z$ est de type $U^+$ par rapport \`a $x$},
autrement dit que $x$
appartient \`a $U^+(z)$. On va appliquer le Lemme~\ref{lemm.liens}
concernant les liens entre la partition associ\'ee au point $z$ et
celle associ\'ee au disque $D$.
On a $ O_x(z) \subset U^+(z) \subset U^+(D)$; or
 l'ensemble $U^+(D)$  est disjoint des ensembles $U^-(D)$,
$\lim^+D$ et $\lim^-D$,
qui contiennent respectivement $U^-(z')$, $\lim^+ z'$ et $\lim^- z'$, par
cons\'equent l'ensemble $O_x(z)$ est disjoint de ces trois ensembles.
D'autre part, l'ensemble $O_x(z)$  contient le point $x$ qui n'appartient
pas \`a $U(z')$
(Affirmation~\ref{affi.quatre-possibilites}), et
$\partial U(z')=(\lim^+ z')  \cup (\lim^- z')$
(Lemme~\ref{lemm.partition}): puisque $O_x(z)$ est un ensemble
connexe, qui contient un point hors de $U(z')$ et qui ne rencontre pas
 $\partial U(z')$, il est disjoint de $U(z')$. On
d\'eduit de tout ceci l'inclusion $O_x(z) \subset U^+(z')$.
Une premi\`ere cons\'equence est que $x \in U^+(z')$,
 donc $z'$ est bien du m\^eme type que $z$ par rapport \`a $x$.
Ensuite, comme $O_x(z')$ est la composante connexe de $U^+(z')$
contenant $x$, et
 que $O_x(z)$ est connexe, on a $O_x(z) \subset O_x(z')$. Par sym\'etrie
 des
 r\^oles de $z$ et $z'$, on a  $O_x(z')=O_x(z)$.

\emph{Le cas o\`u $z$ est de type $U^-$} se traite bien s\^ur en
appliquant le cas $U^+$ \`a $h^{-1}$. D'autre part, en \'echangeant les
r\^oles de $z$ et $z'$, on obtient
$$x \in U^+(z) \Leftrightarrow x \in U^+(z')
\quad\text{et}\quad
x \in U^-(z) \Leftrightarrow x \in U^-(z').$$
Consid\'erons maintenant \emph{le cas o\`u
$z$ est de  type $\lim^+$}; autrement dit, $x$ appartient \`a
$\lim^+(z)$. D'apr\`es ce qui pr\'ec\`ede, $z'$ ne peut pas
\^etre de type $U^-$ ou $U^+$. D'autre part,
dans ce cas, $x$ est dans $\lim^+(D)$, donc n'est pas dans
$\lim^-(D)$, qui contient $\lim^-(z')$: par cons\'equent $z'$ n'est pas
de type $\lim^-$, donc $z'$ est de type $\lim^+$. \emph{Le cas o\`u
$z$ est de type $\lim^-$} s'en d\'eduit.
\end{proof}

\subsection{Troisi\`eme \'etape: existence de composantes de Reeb}

\begin{proof}[\proofname\ de l'Affirmation~\ref{affi.existence-crg}]
Commen\c{c}ons par rappeler que,
s'il n'est pas
\'egal \`a $\{y\}$,
l'ensemble $F_i(x,y)$
 s\'epare
$x$ et $y$ (Lemme~\ref{lemm.topologie-crg}).

L'ensemble   $F_i(x,y)$ est inclus dans $\partial B_i(x)$, donc le
type des point de  $F_i(x,y)$ par rapport \`a $x$ est bien d\'efini
(Affirmation~\ref{affi.quatre-possibilites}).
Puisque $F_i(x,y)$ est connexe (Lemme~\ref{lemm.topologie-crg}), et
que ce type est localement constant
(Affirmation~\ref{affi.invariance}), tous les points de  $F_i(x,y)$
sont du m\^eme type par rapport \`a $x$.

\emph{Consid\'erons le cas o\`u tous les points de  $F_i(x,y)$ sont
de type
$\lim^+$}.
 Dans ce cas, pour tout $z$ dans $F_i(x,y)$, le couple
$(z,x)$ est singulier. L'ensemble $F=\{x\}$ v\'erifie alors la
conclusion de l'Affirmation~\ref{affi.existence-crg}.

\emph{On se place maintenant dans le cas o\`u tous les points de
$F_i(x,y)$ sont de type
 $U^+$}. Par connexit\'e et invariance locale,
l'ensemble $O_x=O_x(z)$ ne d\'epend pas du point $z$ de $F_i(x,y)$. On
pose
alors $F= \partial O_x.$

 D'apr\`es la Proposition~\ref{prop.bumpwire}, pour tout
point $z$ de $F_i(x,y)$, $F$ est inclus dans $\partial U^+(z)$, qui est
inclus dans $\lim^+z$ (Lemme~\ref{lemm.partition}): autrement dit,
$F_i(x,y) \times F  \subset
\sing(h)$ (Lemme~\ref{lemm.singulier-limite}). D'autre part, on a d\'ej\`a
vu que $F_i(x,y)$ est
connexe, et $F$ est \'egalement connexe  d'apr\`es le
Th\'eor\`eme~\ref{theo.connexite} (car $F$ est la fronti\`ere de $O_x(z)$,
qui est une composante connexe de $\bbR^2 \setminus \adhe(U(z))$,
et $U(z)$
est un ouvert connexe).
Il reste \`a prouver les propri\'et\'es de
s\'eparation de la d\'efinition des composantes de Reeb:
\begin{enumerate}
\item  $F_i(x,y)$, s'il n'est pas \'egal \`a $\{y\}$, s\'epare
$F\cup\{x\}$ et
$y$
\item  $F$ s\'epare $x$ et $F_i(x,y) \cup\{y\}$
\end{enumerate}
On remarque tout d'abord  que les ensembles $F$ et $F_i(x,y)$ sont
disjoints, puisque  $F_i(x,y) \times F \subset \sing(h)$.

Montrons alors le premier point. Pour chaque $z$ de $F_i(x,y)$, $z$
appartient \`a	$U(z)$, qui est
disjoint de $\adhe(U^+(z))$, qui contient $\adhe(O_x)$. Donc l'ensemble
$F_i(x,y)$ est disjoint de $\adhe(O_x)$, or cet ensemble est connexe
et contient
$x$. Si $F_i(x,y)$ n'est pas \'egal \`a $\{y\}$, il s\'epare $x$ et
$y$, par cons\'equent il s\'epare aussi $\adhe(O_x)$ et $y$, et
l'ensemble  $F\cup\{x\}$ est inclus dans $\adhe(O_x)$. Ceci prouve le
premier point.

Montrons le second point. Ce qui pr\'ec\`ede montre que $y \not \in
\adhe(O_x)$, et comme $x \in \adhe(O_x)$, l'ensemble $F=\partial O_x$
s\'epare $x$ et $y$. Il reste \`a voir que $F_i(x,y)$ est inclus dans
la composante connexe du compl\'ementaire de $F$ qui contient $y$,
not\'ee $O_y(F)$. Or dans le cas contraire, l'ensemble $O_y(F) \cup F$
est disjoint
de $F_i(x,y)$, par ailleurs il est  connexe
(Proposition~\ref{prop.bumpwire}),
mais ceci contredit le fait que $F_i(x,y)$ s\'epare $F$ et $y$
(premier point).

Il d\'ecoule de ces deux propri\'et\'es que le couple $(F,F_i(x,y))$ est
 une composante de Reeb pour $(x,y)$.
 La proposition est donc prouv\'ee dans
les cas $\lim^+$ et $U^+$. Les cas $\lim^-$ et $U^-$ s'en d\'eduisent.
\end{proof}

\subsection{Quatri\`eme \'etape: composantes de Reeb  minimales}

\begin{proof}[\proofname\ de l'Affirmation~\ref{affi.crg-minimales}]
On d\'efinit (pour tout point $z$ et toute partie $E$ du plan)
$$
d_h(z,E)=\min\{d_h(z,z') : z' \in E\}.
$$
\paragraph{Cas non d\'eg\'en\'er\'e}
\emph{On se place d'abord dans le cas o\`u la composante $(F,G)$ est non
d\'eg\'en\'er\'ee}: dans ce cas, $F$ s\'epare $x$ de $G\cup\{y\}$,
 et $G$ s\'epare $y$ de $F \cup\{x\}$.
L'entier  $i$ sera d\'efini comme $i=d_h(x,F)$.
Nous allons prouver successivement:
\begin{enumerate}
\item $F \subset \partial B_{d_h(y,G)}(y)$
\item $F_{d_h(y,G)}(y,x) \subset F$
\item $d_h(x,F)=d-d_h(y,G)$ (rappelons que $d=d_{h}(x,y)$)
\end{enumerate}

\subparagraph{Premier point}
Puisque la composante est non d\'eg\'en\'er\'ee, le point $y$
n'app\-ar\-tient
pas \`a $G$, et on a $d_h(y,G) \geq 1$.
Soit $z$ un point de $G$
tel que $d_h(y,z)=d_h(y,G)$. L'ensemble $B_{d_h(y,G)}(y)$ est ouvert,
c'est donc
un voisinage de $z$; il est invariant par $h$. Si $z'$ est un point
quelconque de $F$,
comme $(F,G)$ est une composante de Reeb, l'un des
deux couples $(z,z')$ et $(z',z)$ est singulier. La d\'efinition de
l'ensemble singulier entra\^\i ne alors que  $z'$ appartient \`a
$\adhe(B_{d_h(y,G)}(y))$.

Il reste \`a voir que $F$ ne rencontre pas $B_{d_h(y,G)}(y)$. Pour cela
on consid\`ere un arc g\'eod\'esique $\gamma$ joignant $y$ \`a un
point $z'$ de $F$ (\cf \figref{fig.distances}). Puisque $G$
s\'epare $F$ de $y$, l'arc $\gamma$ rencontre
$G$ en un point $z$. Puisque $d_h(z,z')=2$
(Lemme~\ref{lemm.distance-singulier}), toute d\'ecomposition de
$\gamma$ a un sommet entre $z$ et $z'$, on en d\'eduit que $d_h(y,z')-1
\geq d_h(y,z) \geq d_h(y,G)$. Par cons\'equent $z'$ n'appartient pas \`a
la $h$--boule $B_{d_h(y,G)}(y)$, ce que l'on voulait.
\begin{figure}[htp]\anchor{fig.distances}
\centerline{\hbox{\input{fig-distances.pstex_t}}}
\caption{\label{fig.distances}$F$ n'est pas trop pr\`es de $y$}
\end{figure}
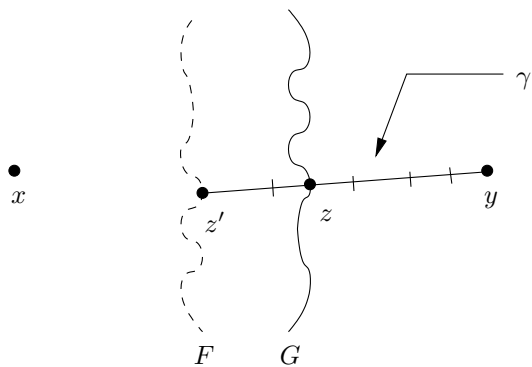

\subparagraph{Deuxi\`eme point}
L'ensemble $F$	s\'epare $x$ et $y$, il est inclus dans
$\partial B_{d_h(y,G)}(y)$, donc il contient $F_{d_h(y,G)}(y,x)$
d'apr\`es la
propri\'et\'e de minimalit\'e de cet ensemble
(Lemme~\ref{lemm.topologie-crg}).

\subparagraph{Troisi\`eme point}
On prouve d'abord l'in\'egalit\'e $d_h(x,F) \geq d-d_h(y,G)$. Soit $z$
un point
quelconque  de $F$. Comme les $h$--boules sont ouvertes, tous les
points $z'$ assez proches de $z$ v\'erifient $d_h(x,z')  \leq
d_h(x,z)$. Puisque $F \subset \partial B_{d_h(y,G)}(y)$, il existe un
tel point
$z'$ v\'erifiant \'egalement $d_h(z',y) \leq d_h(y,G)$. On a alors
\begin{eqnarray*}
d_h(x,z) & \geq & d_h(x,z') \\
	 & \geq & d_h(x,y) - d_h(z',y) \\
	& \geq	&  d - d_h(y,G).
\end{eqnarray*}
Prouvons l'in\'egalit\'e inverse. Soit $\gamma=[x_0 x_1]* \cdots *
[x_{d-1} x_d]$ une d\'ecomposition minimale d'un arc g\'eod\'esique de
$x=x_0$ \`a
$y=x_d$.
Comme $F$ s\'epare $x$ et $y$, l'arc $\gamma$ rencontre
$F$. Comme le sous-arc $[x_{d-d_h(y,G)}x_d]_\gamma$ est inclus dans
$B_{d_h(y,G)}(y)$,  ce sous-arc est disjoint de
$F$ (par le premier point), c'est donc que
$[x_0 x_{d-d_h(y,G)}]_\gamma$ rencontre $F$. Comme $[x_0
x_{d-d_h(y,G)}]_\gamma \subset
B_{d-d_h(y,G)}(x)$, on a bien  $d_h(x,F) \leq d-d_h(y,G)$.

\paragraph{Fin de la preuve dans le cas non d\'eg\'en\'er\'e}
On  recopie maintenant la preuve des premier et deuxi\`eme points
en \'echangeant les r\^oles de $x$ et de $y$ (et donc ceux de $F$ et de
$G$). On obtient ainsi que $F_{d_h(x,F)}(x,y) \subset G \subset \partial
B_{d_h(x,F)}(x)$.
Le troisi\`eme point montre alors que l'entier $i=d_h(x,F)$ convient.
Ceci termine la preuve dans le cas non d\'eg\'en\'er\'e.

\paragraph{Cas d\'eg\'en\'er\'e}
\emph{Il reste \`a traiter les cas o\`u la composante $(F,G)$ est
d\'eg\'en\'er\'ee}.
Dans le cas o\`u $G$ contient $y$ et $F$ contient $x$, le couple
$(x,y)$ est singulier pour $h$ ou $h^{-1}$, et la preuve est
tr\`es simple: on utilise le
Lemme~\ref{lemm.distance-singulier}, et on voit que
$d_h(x,y)=2$, $F_1(y,x)=\{x\}$ et $F_1(x,y)=\{y\}$ (par d\'efinition).

Traitons rapidement le cas o\`u $G$ contient $y$ et $F$ ne contient
pas $x$ en suivant la m\^eme d\'emarche que dans le cas non
d\'eg\'en\'er\'e.
 Le produit $\{y\} \times F$ est
singulier pour $h$ ou $h^{-1}$, donc $F \subset \partial B_1(y)$
(Lemme~\ref{lemm.distance-singulier}).
Avec les m\^emes arguments que dans la situation non d\'eg\'en\'er\'ee
(deuxi\`eme et troisi\`eme points), on en d\'eduit
$F_1(y,x) \subset F$, puis $d_h(x,F)=d-1$.
 On prouve ensuite que $G \subset \partial B_{d-1}(x)$ comme pour le
premier point du cas non
d\'eg\'en\'er\'e. On en d\'eduit que $F_{d-1}(x,y)=\{y\}$ par
d\'efinition. Tout ceci  prouve  que l'entier $i=d-1$ convient.
Le cas o\`u $F$ contient $x$ et $G$ ne contient
pas $y$ est sym\'etrique.
\end{proof}

\subsection{Derni\`ere \'etape: preuve du
Th\'eor\`eme~\ref{theoAbis}}

\begin{proof}[\proofname\ du Th\'eor\`eme~\ref{theoAbis}]
Soient $x$ et $y$ deux points du plan, et $d=d_h(x,y)$.
Montrons tout d'abord que pour chaque entier $j$ entre $1$ et $d-1$,
 $(F_{d-j}(y,x), F_j(x,y))$ est une composante de Reeb
 pour $(x,y)$.
D'apr\`es l'Affirmation~\ref{affi.existence-crg}, il existe un ensemble
$F$ tel que $(F, F_j(x,y))$ est une composante de Reeb
 pour $(x,y)$.
 On applique alors l'Affirmation~\ref{affi.crg-minimales} \`a cette
 composante.
 L'entier $i$ fourni par cette
affirmation v\'erifie notamment $F_i(x,y) \subset F_j(x,y)$; il
co\"{\i}ncide alors avec l'entier $j$, sans quoi ces deux ensembles
seraient respectivement inclus dans les ``$h$--sph\`eres'' $\partial
B_i(x)$ et $\partial B_j(x)$, mais celles-ci sont sont disjointes si $i
\neq j$.
On a donc, d'apr\`es l'Affirmation~\ref{affi.crg-minimales},
$F_{d-j}(y,x) \subset F$, par cons\'equent
$$F_{d-j}(y,x) \times F_j(x,y)
\subset F \times F_j(x,y) \subset \sing(h) \mbox{ ou }
\sing(h^{-1}).$$
Il reste \`a montrer les propri\'et\'es de s\'eparation requises par
la d\'efinition des composantes de Reeb.
Or le couple $(F, F_j(x,y))$ v\'erifie ces propri\'et\'es
(Affirmation~\ref{affi.existence-crg}), et comme $F_{d-j}(y,x) \subset
F$, on voit facilement que le couple $(F_{d-j}(y,x), F_j(x,y))$ les
v\'erifient \'egalement: notamment,  $F_j(x,y)$ s\'epare $y$ et
$F_{d-j}(y,x) \cup\{x\}$, ce qui constitue la propri\'et\'e 2 de la
d\'efinition des composantes, et
on en d\'eduit la propri\'et\'e 3 par sym\'etrie des r\^oles de
$x$ et $y$.

On a montr\'e que le couple $(F_{d-j}(y,x),
F_j(x,y))$ est une composante de Reeb
 pour $(x,y)$.
La minimalit\'e de cette  composante de Reeb
 est \'evidente au vu du
lemme topologique~\ref{lemm.topologie-crg}.
L'Affirmation~\ref{affi.crg-minimales} montre maintenant qu'il n'y a
pas plus de
$d-1$ composantes de Reeb  minimales pour $(x,y)$,
ce qui ach\`eve la d\'emonstration du th\'eor\`eme.
\end{proof}

\section{G\'eod\'esiques infinies et points Nord et Sud}

\label{sec.geodesiques-bout}
Dans cette partie, nous montrons quelques propri\'et\'es simples des
suites g\'e\-o\-d\'es\-iques infinies. Rappelons que celles-ci
apparaissent
naturellement dans
l'\'etude de la dynamique locale (voir l'introduction).
Dans la Section~\ref{ssec.geo-infinies}, on d\'efinit les g\'eod\'esiques
infinies et la relation d'\'equivalence
 ``d\'efinir les m\^emes  bouts'', et on d\'eduit du
 Th\'eor\`eme~\ref{theoAbis}
que deux g\'eod\'esiques infinies d\'efinissant les m\^emes bouts
induisent
la m\^eme suite de composantes de Reeb et le m\^eme mot horizontal.
Dans la Section~\ref{ssec.NS}, on associe \`a toute classe d'\'equivalence
$[\Gamma]$ de g\'eod\'esiques infinies une topologie sur le plan
augment\'e de deux
points $N$ et $S$, situ\'es ``\`a l'infini de part et d'autre de
$\Gamma$''\footnote{Dans la situation d\'ecrite en introduction  o\`u $h$
est le
relev\'e canoniquement associ\'e \`a un point fixe $x_0$ d'un
hom\'eomorphisme de surface, il faut penser \`a $N$ comme \'etant
un ``relev\'e'' du point fixe $x_0$: notamment, la projection sur la
surface d'une suite qui tend vers $N$ tend vers $x_0$.}.
 En un certain sens, les points $N$ et $S$ repr\'esentent aussi les
 ``deux façons d'aller \`a
l'infini'' sous l'it\'eration de $h$ (relativement \`a $[\Gamma]$);
ceci est pr\'ecis\'e
dans la Section~\ref{ssec.partition-dynamique}. La
Section~\ref{ssec.circulation} \'enonce une propri\'et\'e qui restreint
fortement les dynamiques possibles dans un bord de composante de
Reeb: ou bien presque toute orbite va de $N$ \`a $S$, ou bien presque
toute orbite va de $S$ \`a $N$.

\subsection{Les g\'eod\'esiques infinies et leurs bouts}
\label{ssec.geo-infinies}
\subsubsection*{G\'eod\'esiques infinies}
\begin{defi}
Une \emph{suite g\'eod\'esique infinie} est un plongement
isom\'etrique de $\bbZ$ dans le plan muni de la distance $d_h$,
autrement dit une suite $(x_k)_{k \in
\bbZ}$ de points du plan, telle que pour tous entiers $k,k'$, on a
$$
d_h(x_k,x_{k'})=|k-k'|.
$$
Une \emph{droite g\'eod\'esique} est une droite topologique
orient\'ee $\Gamma$, munie d'une d\'e\-com\-pos\-i\-tion
$$
\Gamma= \cdots * \gamma_{-1} * \gamma_0 * \gamma_{1} * \cdots
$$
 o\`u les sous-arcs $\gamma_i$ sont libres, et dont les sommets
 forment une suite g\'eod\'esique infinie $(x_k)$. Les indices sont
 choisis dans le sens croissant de l'orientation de $\Gamma$.
\end{defi}

\begin{defi} Le \emph{mot horizontal associ\'e} \`a la suite
g\'eod\'esique infinie	$(x_k)_{k \in \bbZ}$ est la suite
$(M_\leftrightarrow(x_{k-1},x_{k+1}))_{k \in  \bbZ}$.
\end{defi}

\begin{defi}
La \emph{suite de composantes de Reeb associ\'ee} \`a  la suite
g\'eod\'esique infinie	$(x_k)_{k \in \bbZ}$ est la suite
$(F_k,G_k)_{k \in \bbZ}$, o\`u
$(F_k,G_k)$ est l'unique composante de Reeb
minimale pour le couple $(x_{k-1},x_{k+1})$; autrement dit,
avec les notations du Th\'eor\`eme~\ref{theoAbis},
$(F_k,G_k)=(F_1(x_{k+1},x_{k-1}),F_1(x_{k-1},x_{k+1}))$.
\end{defi}
On associe \'egalement \`a toute droite g\'eod\'esique un
mot horizontal et une suite de composantes, \textit{via} la suite de
ses sommets (la Proposition~\ref{prop.memes-bouts} ci-dessous va
entra\^\i ner
que ceci ne d\'epend pas de la
d\'ecomposition choisie).

\begin{rema}\label{rema.non-degeneree}
Une composante de Reeb associ\'ee \`a une suite g\'eod\'esique infinie
n'est jamais d\'eg\'en\'er\'ee. En effet, supposons (pour fixer les
id\'ees) que la composante $(F_{0},G_{0})$ est d\'eg\'en\'er\'ee,  avec
$x_{1} \in G_{0}$.
La d\'efinition des composantes de Reeb montre qu'alors $x_{1} \in
\partial B_{1}(x_{-1})$. Mais d'autre part $x_{1} \in B_{1}(x_{2})$,
et les boules sont des ouverts du plan, par cons\'equent les boules
$B_{1}(x_{-1})$ et $B_{1}(x_{2})$ se rencontrent, et l'in\'egalit\'e
triangulaire donne $d_{h}(x_{-1}, x_{2})\!\leq 2$, ce qui contredit le
caract\`ere g\'eod\'esique.
\end{rema}

\subsubsection*{Bouts}
\begin{defi}
On dira que les suites g\'eod\'esiques $(x_k)_{k \in
\bbZ}$ et  $(x_k')_{k \in \bbZ}$ \emph{d\'e\-fin\-iss\-ent les m\^emes
bouts} s'il existe un entier $p$ tel que pour tout entier $k$,
$d_h(x_k,x'_k) \leq p$.\footnote{Cette d\'efinition peut para\^\i tre
un peu restrictive;
en r\'ealit\'e, c'est simplement qu'elle est adapt\'ee \`a notre
espace m\'etrique $(\bbR^2, d_h)$, qui se trouve avoir une
g\'eom\'etrie proche de celle d'un arbre. Par exemple, on pourrait
montrer que si
$\exists p $ tel que $\forall k$, $\exists k'$ tel que
$d_h(x_k,x_{k'}')\leq p$, alors la suite g\'eod\'esique infinie
$(x_k')$ d\'efinit les m\^emes bouts que l'une des deux suites $(x_k)$ et
$(x_{-k})$.
\`a l'oppos\'e, on pourrait \'egalement montrer que les deux suites
d\'efinissent les m\^emes bouts si et seulement si $\exists k_0$ tel
que $\forall k$, $d_h(x_k,x_{k+k_0}') \leq 2$.}
\end{defi}
Cette relation est clairement une relation d'\'equivalence.
Elle induit une relation d'\'equivalence sur les droites
g\'eod\'esiques, \textit{via} les suites g\'eod\'esiques de
sommets. On notera $[(x_k)_{k \in \bbZ}]$ et $[\Gamma]$ les classes
d'\'equivalences.

\begin{rema}\label{rema.suite-courbe}
Pour toute  suite g\'eod\'esique $(x_k)_{k \in
\bbZ}$, il existe une droite g\'e\-o\-d\'es\-ique $\Gamma'$
 dont les sommets forment une suite
g\'e\-o\-d\'es\-ique d\'efinissant les m\^emes bouts que  $(x_k)_{k \in
\bbZ}$.
\end{rema}
Montrons rapidement cette remarque. Il est clair qu'il existe une courbe
$\Gamma$, qui passe successivement par tous les sommets de la suite
g\'eod\'esique donn\'ee, telle que chaque sous-courbe $[x_{k} x_{k+1}
]_{\Gamma}$ est un arc libre.
En effa\c cant les boucles de $\Gamma$, on obtient une courbe $\Gamma'$
sans point double  (comme dans la
preuve du Corollaire~\ref{coro.definition-arcs}, \cf
\figref{fig.simplification}); cette courbe $\Gamma'$ est naturellement
d\'ecompos\'ee en sous-arcs $[x'_{k} x'_{k+1} ]_{\Gamma'} \subset
[x_{k} x_{k+1} ]_{\Gamma}$. En utilisant cette derni\`ere inclusion
et le caract\`ere g\'eod\'esique de la suite $(x_k)_{k \in \bbZ}$,
on v\'erifie que la nouvelle suite $(x'_k)_{k \in
\bbZ}$ est encore une suite g\'eod\'esique.
Il reste \`a voir que $\Gamma'$ est bien d'une droite topologique. Toute
partie	compacte du plan est incluse dans une boule $B_{k}(x'_{0})$;
comme la suite des sommets de $\Gamma'$ est g\'eod\'esique, seule une
partie born\'ee de $\Gamma'$ peut rester dans cette boule, ce qui montre
que $\Gamma'$ est proprement plong\'ee dans le plan.

Le r\'esultat principal concernant les bouts est le suivant.
\begin{prop}\label{prop.memes-bouts}
Deux suites g\'eod\'esiques infinies
d\'efinissant les m\^emes bouts ont les m\^emes mots horizontaux et
les m\^emes suites de composantes de Reeb associ\'ees
 (\`a d\'ecalage d'indice
pr\`es).
\end{prop}
Cet \'enonc\'e est avant tout une cons\'equence de ``l'unicit\'e'' dans le
Th\'eor\`eme~\ref{theoAbis}. La fin de cette section est consacr\'ee
\`a
la preuve de la Proposition~\ref{prop.memes-bouts}, qui sera organis\'ee
comme suit.
On remarque d'abord qu'une composante de Reeb pour un
couple de points donn\'e est encore une composante de Reeb pour
beaucoup d'autres couples de points
(Affirmation~\ref{affi.crg-reutilisee}). Pour une suite g\'eod\'esique
(finie ou infinie), cette propri\'et\'e se traduit
par une compatibilit\'e entre les composantes de Reeb associ\'ees aux
diff\'erents couples de sommets (Lemme~\ref{lemm.compatibilite-crg}).
Ceci permet de clarifier la topologie des
composantes de Reeb associ\'ees aux suites g\'eod\'esiques infinies
et de leurs composantes compl\'ementaires
(Lemme~\ref{lemm.pasdenom}). Nous en d\'eduirons la proposition.

Si $A$ est une	partie du plan et $x$ un point hors de $A$,  nous noterons
syst\'ema\-ti\-que\-ment $O_x(A)$ la composante connexe
du compl\'ementaire de $A$ qui contient $x$. La preuve de
l'affirmation suivante est laiss\'ee au lecteur.
\begin{affi}\label{affi.crg-reutilisee}
Soient $x$ et $y$ deux points du plan, et $(F,G)$ une composante de
Reeb  pour $(x,y)$, qui est non
d\'eg\'en\'er\'ee: $F$ ne contient pas $x$ et $G$ ne contient pas
$y$.  Si $x'$ est un point de  $O_x(F)$,
  et si $y'$ est un point de $O_y(G)$, alors
$(F,G)$ est encore une composante de Reeb  pour
$(x',y')$. De plus, si $(F,G)$ est minimale pour $(x,y)$, elle est aussi
minimale pour $(x',y')$.
\end{affi}
Nous reprenons les notations de la
Section~\ref{ssec.construction} concernant les composantes de Reeb
minimales.
\begin{lemm}[Compatibilit\'e des composantes]\label{lemm.compatibilite-crg}
Consid\'erons une suite g\'eo\-d\'e\-sique (finie) $(x=x_0, \dots,
y=x_d)$, et $i$ un
entier entre $1$ et $d-1$. Alors
$$F_{d-i}(y,x) = F_1(x_{i+1},x_{i-1}) \mbox{ et }
F_i(x,y)=F_1(x_{i-1},x_{i+1}).
$$
\end{lemm}
On en d\'eduit notamment que pour tous $0 \leq i < j \leq d$,
les composantes de  Reeb  minimales pour
$(x_i,x_j)$ sont des  composantes de  Reeb  minimales pour
$(x_0,x_d)$.

\begin{proof}[\proofname\ du Lemme~\ref{lemm.compatibilite-crg}]
\emph{Supposons tout d'abord que $1 < i < d-1$}. Dans ce cas,
d'apr\`es les d\'efinitions, la composante $(F_{d-i}(y,x), F_i(x,y))$
pour $(x,y)$ n'est pas d\'eg\'en\'er\'ee.
L'ensemble $F_{i}(x,y)$ est inclus dans la fronti\`ere de
la $h$--boule $B_i(x)$ (Lemme~\ref{lemm.topologie-crg}).
L'in\'egalit\'e triangulaire permet de voir que
$F_{i}(x,y)$ ne rencontre pas la $h$--boule $B_{d-i-1}(y)$.
Or cette $h$--boule est connexe, elle contient le point $x_{i+1}$,
on en d\'eduit
que $x_{i+1}$ est dans $O_y(F_{i}(x,y))$. De mani\`ere sym\'etrique, on
montre que $x_{i-1}$ est dans
$O_x(F_{d-i}(y,x))$. L'Affirmation~\ref{affi.crg-reutilisee} nous dit
maintenant que	$(F_{d-i}(y,x), F_i(x,y))$ est une composante de Reeb
 minimale pour le couple  $(x_{i-1},x_{i+1})$. Or
d'apr\`es le Th\'eor\`eme~\ref{theo.crgm}, ce couple admet une unique
composante de Reeb  minimale, d'o\`u les deux \'egalit\'es.

Si $i=1$ ou $d-1$, il faut \'etudier \`a part les cas o\`u la composante
est d\'eg\'en\'er\'ee, ce qui ne pose pas de probl\`eme. Les
d\'etails sont laiss\'es au lecteur.
\end{proof}

Soit $(x_k)_{k \in \bbZ}$ une suite g\'eod\'esique infinie,
$(F_k,G_k)$ la suite de composantes associ\'ees. Comme dans la preuve
du lemme ci-dessus, l'in\'egalit\'e triangulaire permet de voir que
 \emph{tous les sommets $x_k$ avec $k \leq -1$ appartiennent \`a
la m\^eme composante connexe de $\bbR^2 \setminus F_0$}; on la notera
$O_-(F_0)$ (\cf \figref{fig.geodesique-infinie}).
On voit \'egalement (par exemple \`a l'aide du lemme)
 que \emph{tous les  sommets $x_k$ avec $k
\geq 0$ appartiennent eux aussi \`a
une m\^eme composante connexe de $\bbR^2 \setminus F_0$}, que l'on notera
$O_+(F_0)$. On d\'efinit de m\^eme des ensembles $O_-(F_k)$  et
$O_+(F_k)$ pour tout entier $k$, et aussi des ensembles
$O_-(G_k)$  et $O_+(G_k)$ (\cf
\figref{fig.geodesique-infinie}).\footnote{Pour compl\'eter la
description, on pourrait \'egalement
d\'efinir
l'ouvert $O_k=O_{+}(F_k) \cap O_{-}(G_k)$.
\`A l'aide du lemme d'Alexander
(Th\'eor\`eme~\ref{theo.lemme-alex}), on peut voir que $O_k$ est
connexe, et qu'il s'agit de l'unique composante connexe de $\bbR^2
\setminus (F_k \cup G_k)$
dont l'adh\'erence rencontre \`a la fois $F_k$ et $G_k$.}
\begin{defi}\label{defi.ensembles-complementaires}
Les ouverts $O_-(F_k)$, $O_+(F_k)$, $O_-(G_k)$	et $O_+(G_k)$ sont
appel\'es respectivement \emph{ensembles
compl\'ementaires n\'egatif et positif} de $F_k$ et de $G_k$.
\end{defi}
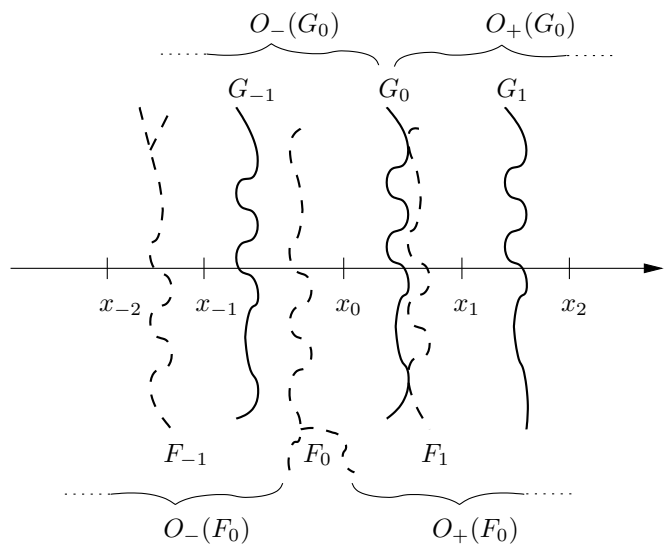
\begin{figure}[htp]\anchor{fig.geodesique-infinie}
\centerline{\hbox{\input{fig-geodesique-infinie2.pstex_t}}}
\caption{\label{fig.geodesique-infinie}Composantes de Reeb associ\'ees \`a
une droite g\'eod\'esique infinie}
\end{figure}

Le lemme suivant d\'ecrit la disposition topologique des
diff\'erents ensembles.
\begin{lemm}\label{lemm.pasdenom}
Soit $(x_k)_{k \in \bbZ}$ une suite g\'eod\'esique infinie.
\begin{enumerate}
\item Les deux ensembles compl\'ementaires positif et n\'egatif de
$F_k$ sont disjoints;
\item on a $k<l \Leftrightarrow O_-(F_k) \not \subseteq O_-(F_l)$;
\item les bords des composantes associ\'ees \`a $(x_k)$ et leurs
ensembles compl\'e\-men\-taires sont tous invariants par $h$;
\item les ensembles compl\'ementaires de $F_k$ ne d\'ependent que de
la classe d'\'equ\-i\-val\-ence $[(x_k)]$.
\end{enumerate}
\end{lemm}
Bien s\^ur, on a des r\'esultats sym\'etriques sur les
ensembles compl\'ementaires de $G_k$.

\begin{proof}[\proofname\ du Lemme~\ref{lemm.pasdenom}]
Sans perte de g\'en\'eralit\'e, on suppose que $k=0$.
\paragraph{Premier point}
 La composante de Reeb
$(F_0,G_0)=(F_1(x_{1},x_{-1}),F_1(x_{-1},x_{1}))$ du couple
$(x_{-1},x_1)$ est non d\'eg\'en\'er\'ee
(Remarque~\ref{rema.non-degeneree}). Par suite $F_0$ s\'epare $x_1$ et
$x_{-1}$, ce qui prouve le premier point, \'etant donn\'ee la d\'efinition
des ensembles compl\'ementaires de $F_0$.

\paragraph{Deuxi\`eme point}
Soit $0<l$. D'apr\`es le premier point, l'ensemble $F_l$ s\'epare
$x_{l-1}$ et $x_l$; Or l'ensemble $O_+(F_0)$ est connexe et contient
ces deux points, c'est donc que $F_l$ rencontre $O_+(F_0)$.
D'autre part, d'apr\`es le Lemme~\ref{lemm.compatibilite-crg}, $F_0$ et
$F_l$ sont respectivement inclus dans les fronti\`eres des $h$--boules
$B_{l+1}(x_{l+1})$ et	$B_{1}(x_{l+1})$; les $h$--boules \'etant
ouvertes, les fronti\`eres des boules de m\^eme centre et de rayon
distincts sont disjointes, et par cons\'equent $F_0$ et $F_l$ sont
disjoints. On en d\'eduit que $F_l$ est enti\`erement inclus dans
$O_+(F_0)$. D'apr\`es le premier point, $F_l$ est disjoint de
$O_-(F_0)$. Ce dernier ensemble est
connexe, disjoint de $F_l$, et il contient $x_{-1}$, il est donc
inclus dans $O_-(F_l)$. Mieux: comme $O_-(F_0) \cup F_0$ est connexe
(Proposition~\ref{prop.bumpwire}), $F_0$ est aussi inclus dans
$O_-(F_l)$. Ceci montre que l'inclusion $O_-(F_0) \subset O_-(F_l)$
est stricte, et entra\^\i ne  l'\'equivalence voulue.

\paragraph{Troisi\`eme point}
Le Lemme~\ref{lemm.compatibilite-crg} montre que $F_0=F_1(x_1,x_{-2})$, on
en d\'eduit l'invariance de $F_0$ par le
Lemme~\ref{lemm.invariance-crg}.
Les ensembles compl\'ementaires sont donc invariants ou libres.
Mais d'apr\`es l'in\'egalit\'e triangulaire, $O_-(F_0)$ contient
par exemple la $h$--boule $B_1(x_{-2})$, qui est invariante. On en
d\'eduit que $O_-(F_0)$ est invariante; on montre de m\^eme que
$O_+(F_0)$ l'est aussi.

Le dernier point sera prouv\'e en m\^eme temps que la
Proposition~\ref{prop.memes-bouts} (sans raisonnement circulaire\ldots{}).
\end{proof}

\begin{proof}[\proofname\ de la Proposition~\ref{prop.memes-bouts},
et du dernier point du Lemme~\ref{lemm.pasdenom}]\nl
Soient $(x_k)$ et $(x_k')$ deux suites g\'eod\'esiques infinies
d\'efinissant les m\^emes bouts, et
$p$ un entier tel que pour tout $k$, $d_h(x_k,x_k') \leq p$.
Nous commen\c{c}ons par prouver que \emph{toute composante de Reeb
associ\'ee
\`a la suite $(x_k)$ est aussi une composante de Reeb associ\'ee \`a la
suite $(x_k')$}. Il suffit bien s\^ur de prouver ceci pour la
composante $(F_0,G_0)=(F_1(x_1,x_{-1}),F_1(x_{-1},x_{1}))$.

Rappelons que (d'apr\`es la d\'efinition), l'ensemble $F_0$ est
inclus dans la $h$--boule $B_2(x_1)$.
D'apr\`es l'in\'egalit\'e triangulaire, l'ouvert $O_-(F_0)$ contient donc
la $h$--boule $B_p(x_{-k})$ pour tout $k \geq p+100$ (\cf
\figref{fig.compo-bouts}); en particulier, il
contient le point $x_{-k}'$. De m\^eme, l'ouvert $O_+(G_0)$ contient
$x_{k}'$.
D'apr\`es l'Affirmation~\ref{affi.crg-reutilisee}, ceci force
$(F_0,G_0)$ \`a \^etre une composante de Reeb minimale pour le couple
$(x_{-k}',x_{k}')$. On applique maintenant le
Lemme~\ref{lemm.compatibilite-crg} sur la compatibilit\'e des
composantes  pour voir que $(F_0,G_0)$ est une
composante de Reeb associ\'ee \`a la suite $(x_k')$.
\begin{figure}[htp]\anchor{fig.compo-bouts}
\centerline{\hbox{\input{fig-compo-bouts2.pstex_t}}}
\caption{\label{fig.compo-bouts}Toute composante de Reeb pour $(x_k)$
est une composante de Reeb pour $({x_k}')$}
\end{figure}
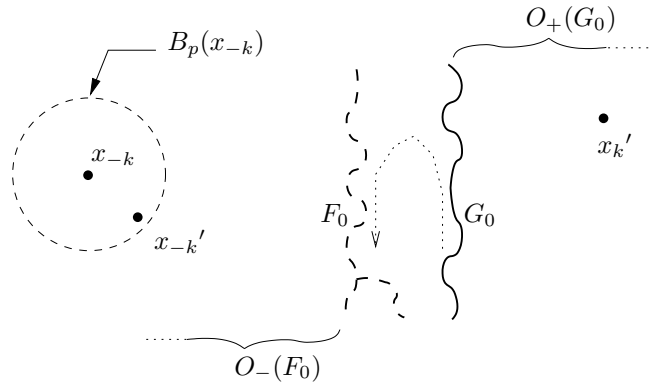

Par sym\'etrie, on obtient que l'ensemble des composantes de Reeb
associ\'ees \`a la suite $(x_k)$ co\"\i ncide avec l'ensemble des
composantes associ\'ees \`a $(x_k')$.
Notons $(F_k,G_k)$ et $(F_k',G_k')$ les deux suites de composantes
associ\'ees.
Pour tout entier $k$, on d\'efinit alors $\phi(k)$ comme l'unique
entier tel que $(F_k,G_k)=(F_{\phi(k)}', G_{\phi(k)}')$; il reste \`a
voir que l'application $\phi$ est un d\'ecalage (\ie $\phi(k)=k+k_0$ pour
un certain $k_0$). On a d\'ej\`a vu que $\phi$ est surjective. Il
suffit alors de montrer  qu'elle est (strictement) croissante.
On a d\'ej\`a vu que $O_-(F_0)$ contient le point $x_{-p-100}'$, on
en d\'eduit que $O_-(F_0)=O_-(F_{\phi(0)}')$; et plus
g\'en\'eralement que tous les ensembles compl\'ementaires
associ\'es \`a $(x_k)$ et $(x_k')$ co\"\i ncident. Ceci
prouve le dernier point du Lemme~\ref{lemm.pasdenom}. D'autre part,
 le Lemme~\ref{lemm.pasdenom} caract\'erise le fait que $k<l$ par la
propri\'et\'e topologique $O_-(F_k) \subset O_-(F_l)$; en appliquant ceci
successivement aux deux suites g\'eod\'esiques, on obtient
$$
\begin{array}{rcl}
k < l & \Leftrightarrow  &   O_-(F_k) \not \subseteq O_-(F_l)  \\
	 & \Leftrightarrow  &	O_-(F_{\phi(k)}') \not \subseteq
	 O_-(F_{\phi(l)}') \\
	 & \Leftrightarrow  &	\phi(k) < \phi(l).
\end{array}
$$
Ceci montre que $\phi$ est strictement croissante, et
termine la preuve en ce qui concerne les composantes de Reeb
g\'en\'eralis\'ees. Puisque les mots horizontaux d\'ependent
uniquement des composantes de Reeb
(Th\'eor\`eme~\ref{theoBter},
Section~\ref{ssec.preuve-1}), on en d\'eduit que les suites
g\'eod\'esiques ont aussi le m\^eme mot horizontal.
\end{proof}

\subsection{Les points $N$ et $S$}
\label{ssec.NS}
\begin{figure}[htp]\anchor{fig.compactification}
\centerline{\hbox{\input{fig-compactification.pstex_t}}}
\caption{\label{fig.compactification}Compactification Nord-Sud}
\end{figure}
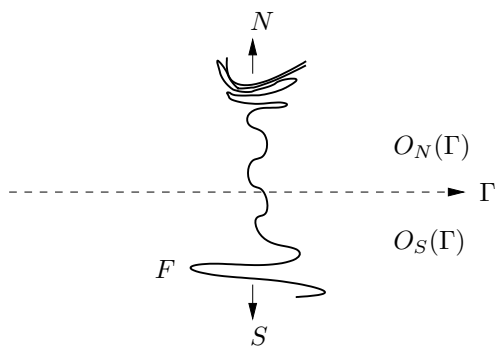

On se donne une droite g\'eod\'esique $\Gamma$.  On appelle $O_N(\Gamma)$
et $O_S(\Gamma)$ les deux composantes connexes de $\bbR^2 \setminus
\Gamma$,
avec la convention que $O_N(\Gamma)$ est situ\'ee \`a gauche de $\Gamma$
quand on
parcourt $\Gamma$ dans le sens positif (\cf
\figref{fig.compactification}).

Soit $F$ une partie ferm\'ee du plan telle
que $F \cap \Gamma$ est compacte; on
peut alors compactifier $F$ en un espace $\widehat F=F \cup \{N,S\}$
de la mani\`ere suivante: on ajoute un
point $N$ par la compactification d'Alexandroff de $F \cap
\adhe(O_N(\Gamma))$, et
un point $S$  par la compactification d'Alexandroff de $F \cap \adhe(
O_S(\Gamma))$.\footnote{Autrement dit, par exemple, une base de voisinages
de $N$ est form\'ee des ensembles $(O_N(\Gamma) \setminus K) \cap F$
o\`u $K$ d\'ecrit l'ensemble des compacts du plan.}
\begin{defi}\label{defi.compactification}
L'espace topologique $\widehat F=F \cup \{N,S\}$ est appel\'e
\emph{compactification Nord-Sud de $F$ (associ\'ee \`a $\Gamma$)}.
\end{defi}
Pr\'ecisons imm\'ediatement qu'en g\'en\'eral cette topologie d\'epend
de $\Gamma$ (et pas seulement de $[\Gamma]$), voir la remarque \`a la
fin de cette partie.

\subsubsection*{Compactification des boules}
Soit $F$ une partie du plan. Supposons que l'intersection $F
\cap \Gamma$ soit  non vide et d'adh\'erence compacte.

\begin{defi}\label{defi.separe-bouts}
On dira alors que \emph{$F$ s\'epare les deux bouts de $\Gamma$} si $F$
s\'epare les deux composantes connexes non born\'ees de $\Gamma
\setminus F$.\footnote{L\`a encore, cette notion d\'epend du choix
de $\Gamma$, mais elle ne d\'epend que de $[\Gamma]$ pour les ensembles
born\'es pour la distance de translation.}
\end{defi}

Par exemple, le premier point du Lemme~\ref{lemm.pasdenom} montre que
chaque bord de
composante de Reeb associ\'ee \`a $\Gamma$ s\'epare les deux
bouts de $\Gamma$.
Soit $y$ un point du plan. Pour tout entier positif $d$,
l'intersection de la boule $B_d(y)$ avec la g\'eod\'esique $\Gamma$ est
incluse dans une partie compacte de $\Gamma$.
Soit $(x_k)$ la suite des sommets de la d\'ecomposition de $\Gamma$; alors
l'ensemble $\adhe(B_1(x_1))$ s\'epare les deux bouts de $\Gamma$,
puisqu'il contient le bord $F_0$ d'une composante de Reeb
associ\'ee.
 Si on choisit $d$ assez grand pour que $B_d(y)$ contienne
 la boule $B_1(x_0)$, alors
 l'adh\'erence de la boule $B_d(y)$  s\'epare
les deux bouts de $\Gamma$.
\begin{affi}\label{affi.independance}
Soit $B_d(y)$ une boule du plan s\'eparant les deux bouts de
$\Gamma$.
Alors la topologie obtenue sur l'ensemble
$\widehat{\adhe}(B_d(y))=\adhe(B_d(y))\cup\{N,S\}$ par la
compactification Nord-Sud  associ\'ee \`a
$\Gamma$  ne d\'epend que de $[\Gamma]$.
\end{affi}

\begin{proof}[\proofname\ de l'Affirmation~\ref{affi.independance}]
Notons $O_-$ et $O_+$ les composantes
connexes du compl\'ementaire de $\adhe(B_{d}(y))$ qui contiennent les
parties non
born\'ees respectivement n\'egative et positive de $\Gamma$.
On d\'efinit de m\^eme les ensembles $O'_-$ et $O'_+$ comme  les
composantes
connexes du compl\'ementaire de $\adhe(B_{d+1}(y))$ qui contiennent les
parties non
born\'ees  de $\Gamma$. On montre facilement que les quatre ensembles
 $O_-$, $O_+$,	$O'_-$ et $O'_+$ ne
d\'ependent que du choix de  $[\Gamma]$ (comme pour le dernier point
du Lemme~\ref{lemm.pasdenom}).
En particulier, puisque $B_d(y)$ s\'epare les deux bouts de $\Gamma$,
les ensembles $O_-$ et $O_+$ sont disjoints, on en d\'eduit que
$B_d(y)$ s\'epare les deux bouts de toute droite g\'eod\'esique
  $\Gamma'$ d\'efinissant les m\^emes bouts que $\Gamma$.
Notons que $O'_-$ et $O'_+$ sont inclus respectivement dans $O_-$
et $O_+$.

Le lemme suivant  est tr\`es classique, on peut le prouver \`a l'aide de
triangulations; les d\'etails sont laiss\'es au lecteur.
\begin{lemm}[R\'egularisation]
\label{lemm.regularisation}
Soit $F$ une partie ferm\'ee et connexe de $\bbR^2$, et $O$ une partie
ouverte contenant $F$.
Alors il existe une partie ferm\'ee et connexe $F_0$, qui est une
sous-vari\'et\'e \`a bord du plan, contenant $F$ et contenue dans $O$.
\end{lemm}
On applique le Lemme~\ref{lemm.regularisation} avec $F=\adhe(B_d(y))$
et $O=B_{d+1}(y)$ (\cf \figref{fig.bandeB}), pour obtenir une
sous-vari\'et\'e \`a bords $F_0$. Il est clair que $F_0$ est
disjointe de $O'_-$ et $O'_+$. De plus, la
composante connexe $O_0^-$ de $\bbR^2 \setminus F_0$ qui contient
$O'_-$ n'est pas born\'ee. Par cons\'equent, le bord de
$O_0^-$ est une droite topologique, on la note $\Delta^-$. On d\'efinit
sym\'etriquement $\Delta^+$. Les deux droites topologiques $\Delta^-$
et $\Delta^+$ sont
disjointes, le th\'eor\`eme de Schoenflies permet de trouver un
hom\'eo\-mor\-phisme $\Phi$ du plan, pr\'eservant l'orientation, tel  que
$\Delta^-=\Phi(\{0\}\times \bbR)$ et
$\Delta^+=\Phi(\{1\}\times \bbR)$ (\cf
\figref{fig.bandeB}). Appelons $B$ la bande $\Phi([0,1]
\times \bbR)$. Par connexit\'e, le ferm\'e
$F=\adhe(B_d(y))$ est inclus dans la bande $B$. On a donc montr\'e:
\begin{affi}\label{affi.choix-de-B}
Il existe une partie $B$ du plan, telle que
\begin{itemize}
\item $B$ contient $\adhe(B_d(y))$, et est disjointe de $O'_-$ et $O'_+$;
\item il existe un hom\'eomorphisme $\Phi$ qui envoie $[0,1] \times
\bbR$ sur $B$;
\item le bord $\Delta^-=\Phi(\{0\} \times \bbR)$ est inclus dans $O^-$;
\item le bord $\Delta^+=\Phi(\{1\} \times \bbR)$ est inclus dans $O^+$.
\end{itemize}
\end{affi}
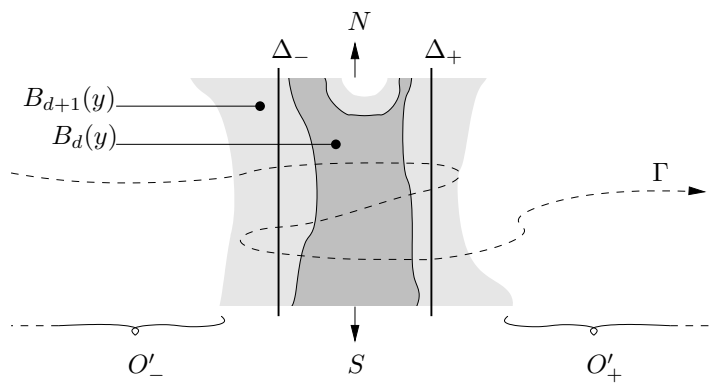
\begin{figure}[htp]\anchor{fig.bandeB}
\centerline{\hbox{\input{fig-bandeB.pstex_t}}}
\caption{\label{fig.bandeB}La compactification ne d\'epend que de
$[\Gamma]$}
\end{figure}
On  consid\`ere la compactification en bouts
$\widehat B=B
\cup \{N,S\}$, o\`u $N$ et $S$ sont choisis comme sur la
\figref{fig.bandeB}.\footnote{On compactifie $\bbR$ de mani\`ere
usuelle en $\widehat \bbR = \{-\infty\} \cup \bbR \cup \{+\infty\}$;
on peut alors construire $\widehat B$ comme le quotient de l'espace $[0,1]
\times \widehat \bbR$ par la relation d'\'equivalence qui identifie les
points de
$[0,1] \times \{+\infty\}$ entre eux (ce qui donne le point $N$), et
les points de
$[0,1] \times \{-\infty\}$ entre eux (ce qui donne le point $S$).}
Gr\^ace \`a la premi\`ere propri\'et\'e de l'affirmation, l'intersection
$B
\cap \Gamma$ est compacte et non vide, et $B$ s\'epare les deux bouts
de $\Gamma$. On en d\'eduit facilement que la topologie de la
compactification Nord-Sud
 associ\'ee \`a $\Gamma$ sur l'ensemble $\widehat{\adhe}(B_d(y))=
 \adhe(B_d(y))
\cup \{N,S\}$ co\"\i ncide avec  la topologie induite par la topologie
de $\widehat B$.
Mais comme les ensembles  $O_-$, $O_+$,  $O'_-$ et $O'_+$ ne
d\'ependent que du choix de  $[\Gamma]$,
tout ceci est encore vraie pour une g\'eod\'esique infinie
$\Gamma'$  d\'efinissant les m\^emes bouts que
$\Gamma$. En particulier, la topologie de la compactification Nord-Sud
sur $\widehat{\adhe}(B_d(y))$ ne d\'epend que de $[\Gamma]$.
\end{proof}

\subsubsection*{Topologie sur $\bbR^2 \cup \{N,S\}$}
Remarquons que si $B_{d'}(y')$ est une boule qui contient
$B_d(y)$, alors la topologie de la compactification Nord-Sud sur
$\widehat{\adhe}(B_{d'}(y'))$ induit
celle sur $\widehat{\adhe}(B_{d}(y))$.

On peut  maintenant d\'efinir une topologie sur $\widehat{\bbR^2}=\bbR^2
\cup
\{N,S\}$, en d\'ecidant qu'une partie $\widehat O$ de $\widehat{
\bbR^2}$ est ouverte si et
seulement si, pour toute boule $B_d(y)$
s\'eparant les deux bouts de $\Gamma$, l'ensemble $\widehat O \cap
\widehat{\adhe}(B_d(y))$ est une partie ouverte de
$\widehat{\adhe}(B_d(y))$.
Cette topologie sera appel\'ee \emph{topologie associ\'ee \`a $\Gamma$}.
Il est clair que, l\`a encore, cette topologie ne d\'epend que du
choix de $[\Gamma]$.

 La droite topologique orient\'ee $h(\Gamma)$ est une
droite g\'eod\'esique, et  elle d\'efinit les m\^emes bouts que
$\Gamma$ (puisque la distance de translation entre un point et son
image par $h$ vaut $1$). On en d\'eduit que l'hom\'eomorphisme $h$
s'\'etend en un
hom\'eomorphisme de $\widehat{\bbR^2}$ qui fixe $N$ et $S$. On notera
encore $h$ cette extension.

\begin{rema}
Soit $F$ une partie ferm\'ee du plan telle que	$F \cap
\Gamma$ est compacte;
la compactification Nord-Sud $\widehat F$ associ\'ee \`a $\Gamma$
co\"\i ncide-t-elle avec la topologie induite par la topologie
de $\widehat{ \bbR^2}$ associ\'ee \`a $\Gamma$? En g\'en\'eral, la
r\'eponse est n\'egative; cependant, les deux topologies co\"\i
ncident effectivement dans le cas des $h$--boules, donc aussi si $F$ est
born\'e pour la distance $d_h$; en fait, cette topologie n'aura
d'int\'er\^et que pour de tels ensembles. Le cas des bords de
composantes de Reeb  associ\'ees \`a $\Gamma$ nous
int\'eressera tout
particuli\`erement.
\end{rema}

\subsection{Partition dynamique du plan}\label{ssec.partition-dynamique}
Comme avant, on se donne une droite g\'eod\'esique $\Gamma$,
et on note $O_N(\Gamma)$ et $O_S(\Gamma)$ les deux composantes connexes du
compl\'ementaire.
\begin{coro}\label{coro.partition-dynamique}
Pour tout point $x$ du plan, de deux choses l'une:
\begin{enumerate}
\item ou bien il existe un entier $n_0$ tel que pour tout $n \geq n_0$,
$h^n(x) \in O_N(\Gamma)$;
 \item ou bien il existe un entier $n_0$ tel que pour tout $n \geq n_0$,
$h^n(x) \in O_S(\Gamma)$.
\end{enumerate}
De plus, cette alternative ne d\'epend que de $[\Gamma]$: si
$[\Gamma']=[\Gamma]$, alors on a l'\'equivalence
$$
\exists n_0, \forall n \geq n_0, h^n(x) \in O_N(\Gamma)
\Longleftrightarrow \exists n_0', \forall n \geq n_0',	h^n(x) \in
O_N(\Gamma').
$$
\end{coro}
\begin{proof}
Soit $y$ un point du plan, et $d$ un entier positif tel que la $h$--boule
$B_d(y)$ s\'epare les deux bouts de $\Gamma$. L'espace
$\widehat{\adhe}(B_d(y))$ est compact, et tous les points de
$\adhe(B_d(y))$ sont errants pour $h$. On en d\'eduit facilement
que l'orbite positive de tout point de $\adhe(B_d(y))$ doit tendre
vers $N$ ou $S$ (voir par exemple \cite{lero2001}, preuve du
corollaire 9.8). Puisque ceci est vrai pour toute $h$--boule s\'eparant
les deux bouts de $\Gamma$, l'orbite positive de tout point $x$ du
plan tend vers $N$ ou $S$ dans $\widehat{\bbR^2}$.
Puisque les ensembles $O_N(\Gamma)\cup\{N\}$ et
$O_S(\Gamma)\cup\{S\}$ sont respectivement des voisinages de
 $N$ et $S$, on a bien la dichotomie annonc\'ee. Comme la topologie sur
$\widehat{\bbR^2}$ ne d\'epend que de $[\Gamma]$, l'alternative ne
d\'epend elle aussi que de $[\Gamma]$.
\end{proof}

On dira que le point $x$ \emph{va vers le
Nord} dans le premier cas, \emph{va vers le Sud} dans le second (et on
utilisera de m\^eme les expressions \emph{venir du Nord} et
\emph{venir du Sud}).
 On d\'efinit  les ensembles
\begin{eqnarray*}
\wnf(\bbR^2)& = & \bigl\{x \in \bbR^2 : \lim_{n \rightarrow -\infty}
h^n(x)=N\bigl\}\\
\wfn(\bbR^2)& = & \bigl\{x \in \bbR^2 : \lim_{n \rightarrow +\infty}
h^n(x)=N \bigl\}
\end{eqnarray*}
et, de mani\`ere similaire, les ensembles $\wsf(\bbR^2)$ et
$\wfs(\bbR^2)$;
puis on pose
\begin{align*}
\wnn(\bbR^2) &=\wnf(\bbR^2) \cap \wfn(\bbR^2); &
\wss(\bbR^2) &=\wsf(\bbR^2) \cap \wfs(\bbR^2) \\
\wns(\bbR^2) &= \wnf(\bbR^2) \cap \wfs(\bbR^2); &
\wsn(\bbR^2) &= \wsf(\bbR^2) \cap \wfn(\bbR^2).
\end{align*}
On notera \'egalement $\wnf(F)$ l'ensemble $\wnf(\bbR^2) \cap F$,
\textit{etc.}.

On a ainsi associ\'e \`a toute classe d'\'equivalence de
droites g\'eod\'esiques  une partition dynamique du plan,
$$
\bbR^2=\wss \sqcup \wsn \sqcup \wns \sqcup \wnn.
$$

\begin{rema}
Si l'on part non plus d'une droite g\'eod\'esique mais
simplement d'une suite g\'eod\'esique infinie $(x_k)$, on peut
r\'ecup\'erer les constructions des Sections~\ref{ssec.NS}
et~\ref{ssec.partition-dynamique}. En effet, on commence par se donner
une droite g\'eod\'esique $\Gamma$ dont la suite des sommets
d\'efinit les m\^emes bouts que $(x_k)$
(Remarque~\ref{rema.suite-courbe}); on construit les topologies
associ\'ees \`a $\Gamma$, et on remarque que ces topologies ne
d\'ependent pas du choix de $\Gamma$ (d'apr\`es
l'Affirmation~\ref{affi.independance}). Nous parlerons donc dans ce
qui suit de la \emph{compactification Nord-Sud associ\'ee \`a une
suite g\'eod\'esique infinie $(x_k)$}.
\end{rema}

\subsection{Sens de circulation dans les bords de composantes de Reeb}
\label{ssec.circulation}
La preuve de la proposition suivante occupera la
Section~\ref{sec.continus-minimaux}.
L'id\'ee intuitive est que les bords de composantes de Reeb
 minimales sont trop \'etroits pour que la
circulation puisse se faire dans les deux sens. Notons que l'\'enonc\'e
est faux pour les composantes de Reeb non minimales, comme le montre
l'exemple de la \figref{fig.flot-perturbe} de
l'Appendice~\ref{sec.exemples}.
\begin{prop}\label{prop.continus-minimaux}~
Soit $(x_k)_{k \in \bbZ}$ une suite g\'eod\'esique infinie, et soit $F$
un bord
de composante de Reeb	associ\'ee \`a
$(x_k)$. Soit $\widehat F$ la compactification Nord-Sud de $F$
associ\'ee \`a $(x_k)$. Alors:
\begin{itemize}
\item l'un des deux ensembles
 $\wsn(F)$ et $\wns(F)$
 est vide;
\item l'autre est un ouvert dense de $F$, et il est connexe;
\item les ensembles $\wnn(F)$ et $\wss(F)$ sont ferm\'es.
\end{itemize}
\end{prop}
On dira que $F$ est \emph{de type dynamique Sud-Nord} si $\wsn(F) \neq
\emptyset$, et que  $F$ est \emph{de type dynamique Nord-Sud} si $\wns(F)
\neq
\emptyset$. D'apr\`es la Proposition~\ref{prop.continus-minimaux}, $F$
est d'un et d'un seul des deux types, et  $F$ est de type Sud-Nord si et
seulement si il existe un ouvert (non vide)  de points de $F$ dans
$\wfn(F)$.

\section{Fl\`eches verticales}
\label{sec.fleches-verticales}
Le r\'esultat principal de cette section est que deux bords de
composantes de Reeb  ``adjacents'' ont le m\^eme
type dynamique (d\'efini pr\'ecis\'ement  \`a la
Section~\ref{sec.geodesiques-bout}):  le ``sens de circulation'', du
Nord vers le Sud ou bien du Sud vers le Nord, est le m\^eme. Ceci
permet d'associer un mot infini \`a toute suite g\'eod\'esique infinie,
 en ins\'erant, dans le mot horizontal, une fl\`eche $\uparrow$
aux endroits correspondant \`a des bords de type Sud-Nord et une fl\`eche
$\downarrow$ pour les bords de type Nord-Sud. Nous utilisons ici tous
les r\'esultats de la section pr\'ec\'edente.

\subsection{\'Enonc\'e pr\'ecis du
Th\'eor\`eme~\ref{theo.fleches-verticales}}
\begin{figure}[htp]\anchor{fig.theo-verticales}
\centerline{\hbox{\input{fig-theo-verticales.pstex_t}}}
\caption{\label{fig.theo-verticales}Fl\`eches verticales et dynamique
des sommets
d'une droite g\'eod\'esique}
\end{figure}
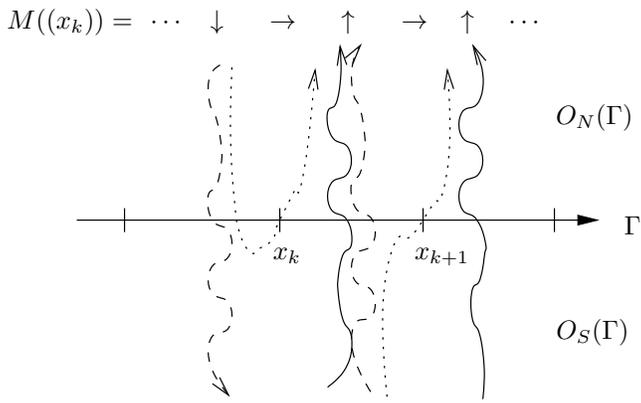

\setcounter{theobis}{3}
\begin{theobis}\label{theoDbis}
Soit $(x_k)_{k \in \bbZ}$ une suite g\'eod\'esique infinie et
$((F_k,G_k))_{k \in \bbZ}$ la
suite des composantes de Reeb  associ\'ee.
Alors pour tout entier $k$, $G_k$ et $F_{k+1}$ ont le m\^eme type
dynamique.

Supposons de plus que $(x_k)$ est la suite des sommets d'une droite
g\'eod\'esique	$\Gamma$. Alors,
\begin{itemize}
\item si l'unique lettre du mot $M_\leftrightarrow(x_{k-1},x_{k+1})$
est $\rightarrow$, alors
$G_k$ et $F_{k+1}$ sont de type dynamique Sud-Nord si et seulement si
le point
$x_k$ va vers le Nord;
\item  si cette lettre est $\leftarrow$, alors
$G_k$ et $F_{k+1}$ sont de type dynamique Sud-Nord si et seulement si
le point
$x_k$ vient du Sud.
\end{itemize}
\end{theobis}
On a bien s\^ur un \'enonc\'e analogue pour $x_{k+1}$ (obtenu \`a partir
du Th\'eor\`eme~\ref{theoDbis} en
inversant l'orientation de $\Gamma$):
\begin{itemize}
\item si $M_\leftrightarrow(x_{k},x_{k+2})=(\leftarrow)$, alors
$G_k$ et $F_{k+1}$ sont de type dynamique Sud-Nord si et seulement si
le point
$x_{k+1}$ va vers le Nord;
\item  si $M_\leftrightarrow(x_{k},x_{k+2})=(\rightarrow)$, alors
$G_k$ et $F_{k+1}$ sont de type dynamique Sud-Nord si et seulement si
le point
$x_{k+1}$ vient du Sud. Le th\'eor\`eme est illustr\'e par la
\figref{fig.theo-verticales}.
\end{itemize}
Notamment, l'\'enonc\'e implique que, \'etant donn\'ees les fl\`eches
horizontales associ\'ees \`a la droite g\'eod\'esique $\Gamma$, les
dynamiques de deux sommets cons\'ecutifs sont \'etroitement reli\'ees.

\subsection{D\'efinition des fl\`eches verticales}
\begin{defi}
Le \emph{mot vertical} associ\'e \`a $(x_k)$ est la suite
$({M_\updownarrow}_k)_{k \in \bbZ}$  d'\'el\'ements de l'alphabet
$\{\uparrow,\downarrow\}$ d\'efinie par
\begin{gather*}
{M_\updownarrow}_k=\ \uparrow \quad  \Leftrightarrow \quad \mbox{$G_k$
et $F_{k+1}$ sont de
type dynamique Sud-Nord}
\\
{M_\updownarrow}_k=\ \downarrow  \quad \Leftrightarrow \quad \mbox{$G_k$
et $F_{k+1}$ sont de
type dynamique Nord-Sud.}
\end{gather*}\end{defi}
On obtient ainsi un invariant de conjugaison, associ\'e \`a la
donn\'ee de $h$ et d'une classe d'\'equivalence de suites
g\'eod\'esiques infinies. Cet invariant compl\`ete le mot horizontal
d\'efini \`a la Section~\ref{ssec.fleches}; la d\'efinition suivante
regroupe
fl\`eches horizontales et verticales dans un m\^eme mot.
\begin{defi}
On associe
\`a toute suite g\'eod\'esique infinie $(x_k)$
le mot
\!$M((x_k))$ obtenu en
intercalant les lettres des mots horizontal et vertical,
$$
M((x_k)):=(\dots, {M_\updownarrow}_{k-1}, {M_\leftrightarrow}_k,
{M_\updownarrow}_k,{M_\leftrightarrow}_{k+1},
{M_\updownarrow}_{k+1}, \dots  ).
$$
\end{defi}
Notons que le mot $M((x_k))$ d\'etermine le type dynamique de chacun
des points $x_{k}$.

Ce mot n'intervient pas dans la suite du texte. Mentionnons cependant une
propri\'et\'e de ce mot, qui fournit \'egalement une interpr\'etation
des fl\`eches qui le composent, et  justifie ainsi le choix de ces
symboles.
 En utilisant les techniques de
d\'ecompositions en briques (voir \cite{leca3}, \cite{lero2001}), on
montrera dans  l'article~\cite{leroinf} les r\'esultats
suivants. Pour toute composante de Reeb
$(F_k,G_k)$ associ\'ee \`a $(x_k)$, il existe une droite de Brouwer
$\Delta_k$ qui s\'epare les deux bords	$F_k$ et $G_k$. De plus,
l'indice partiel de $h$ entre $\Delta_k$ et $\Delta_l$ (au sens de
\cite{lero2001}) est enti\`erement
d\'etermin\'e par le mot $M((x_k))$. Plus pr\'ecis\'ement, cet indice
est \'egal
au ``nombre de tours qu'effectue la fl\`eche'' dans le mot
$$
({M_\leftrightarrow}_k, {M_\updownarrow}_k, \dots
,{M_\updownarrow}_{l-1}, {M_\leftrightarrow}_{l} ).
$$
Par exemple,  le mot $(\rightarrow, \uparrow, \rightarrow, \downarrow,
\leftarrow)$ correspond \`a un indice \'egal \`a $-1/2$.

\subsection{Preuve dans le cas simple (mot attractif ou r\'epulsif)}
\subsubsection*{Notations}
Soit $(x_k)$ une suite g\'eod\'esique infinie.
D'apr\`es la Section~\ref{sec.geodesiques-bout},

\noindent -- on peut choisir une droite g\'eod\'esique $\Gamma$ dont la
   suite des sommets d\'efinit les m\^emes bouts que $(x_k)$
   (Remarque~\ref{rema.suite-courbe});

\noindent -- la suite de composantes  associ\'ee \`a $\Gamma$ est
la m\^eme
   que celle associ\'ee \`a $(x_k)$ (Proposition~\ref{prop.memes-bouts});

\noindent -- la d\'efinition des types dynamiques de $G_k$ et $F_{k+1}$ ne
   d\'epend pas du choix de cette courbe g\'eod\'esique $\Gamma$.

Par cons\'equent, on peut supposer dans toute la preuve du
Th\'eor\`eme~\ref{theoDbis} que
$(x_k)$ est la suite des sommets d'une droite g\'eod\'esique
$\Gamma$.
Comme avant, on notera $O_N(\Gamma)$ et $O_S(\Gamma)$ les composantes
connexes gauches et
droites du compl\'ementaire de $\Gamma$.
Afin d'all\'eger les notations,  on supposera
(sans perte de g\'en\'eralit\'e) que $k=0$, et on notera $G=G_0$ et
$F=F_1$. On note
$$
\Gamma= \cdots * \gamma_{-1} * \gamma_0 * \gamma_{1} * \cdots
$$
la d\'ecomposition de $\Gamma$. En appliquant le
Lemme~\ref{lemm.ordre-crg} sur le franchissement des composantes de Reeb,
 et le Lemme~\ref{lemm.compatibilite-crg} sur la compatibilit\'e des
 composantes de Reeb,
on voit que  $\Gamma \cap (G \cup
F)$ est non vide et inclus dans l'int\'erieur de l'arc $\gamma_0$.
\emph{On se place ici dans les cas ``attractif'' ou ``r\'epulsif'',
\ie\ on suppose que  le mot
$M_\leftrightarrow(x_{-1},x_2)$ vaut $(\rightarrow \leftarrow)$ ou
$(\leftarrow \rightarrow)$.} De plus,
il suffit de faire la preuve quand
$M_\leftrightarrow(x_{-1},x_2) = (\rightarrow \leftarrow)$, l'autre
cas s'en d\'eduisant en appliquant le r\'esultat \`a $h^{-1}$.
\subsubsection*{D\'emonstration}
Le fait que $\Gamma$ soit g\'eod\'esique va entra\^\i ner la
propri\'et\'e-cl\'e suivante.
\begin{affi}\label{affi.fleches-verticales-facile1}
Pour tout $n \geq 1$, $h^n(\gamma_0) \cap \Gamma=\emptyset$.
\end{affi}
\proof
Soit $n \geq 1$.
On doit montrer que pour tout entier $k$, $h^n(\gamma_0) \cap
\gamma_k=\emptyset$.
\begin{itemize}
\item Pour $k=0$, ceci vient de la libert\'e de $\gamma_0$
(Corollaire~\ref{coro.liberte}).
\item Pour $k=1$ ou $k=-1$, ceci vient de l'hypoth\`ese
$M_\leftrightarrow(x_{-1},x_2) = (\rightarrow \leftarrow)$ et de la
d\'efinition des fl\`eches horizontales
(Affirmation~\ref{affi.deux-possibilites} et
D\'efinition~\ref{defi.fleches-horizontales}).
\item Pour $|k|>1$, ceci provient de la g\'eod\'esicit\'e de
$\Gamma$
(\cf premier point du Lemme~\ref{lemm.iteres-arc-geo} sur les it\'er\'es
d'un arc g\'eod\'esique).\qed
\end{itemize}

Soit $D_0$ un disque topologique ferm\'e libre, v\'erifiant $\Gamma \cap
D_0 =
\gamma_0$ et $\inte(\gamma_0) \subset \inte(D_0)$. On obtient
facilement un tel disque en \'epaississant $\gamma_0$ (\`a l'aide du
th\'eor\`eme de Schoenflies).
On peut alors renforcer l'affirmation pr\'ec\'edente:
\begin{affi}\label{affi.fleches-verticales-facile2}
Sous les hypoth\`eses ci-dessus,
$h^n(D_0) \cap \Gamma=\emptyset$
pour tout $n \geq 1$.
\end{affi}
\begin{proof}
On pourrait faire la d\'emonstration en recopiant la preuve de
l'aff\-irm\-a\-tion pr\'ec\'edente, en utilisant une version adapt\'ee du
lemme de Franks. Voici une autre possibilit\'e, qui consiste \`a
``perturber'' $\Gamma$.
On raisonne par l'absurde. Supposons que $n$ est un entier strictement
positif tel que $h^n(D_0)$ rencontre $\Gamma$, et soit $y \in D_0 \cap
h^{-n}(\Gamma)$. Soit $\gamma_0'$ un arc inclus dans $D_0$, de m\^emes
extr\'emit\'es que $\gamma_0$, et passant par $y$. On obtient une
nouvelle droite g\'eod\'esique $\Gamma'$ \`a partir de
$\Gamma$ en rempla\c{c}ant $\gamma_0$ par $\gamma_0'$.
Puisque $D_0$ est libre, il est disjoint de tous ses it\'er\'es, et
$h^n(y)$ n'est pas sur $\gamma_0$; par cons\'equent  $y \in  \gamma'_0
\cap
h^{-n}(\Gamma')$. D'autre part, la nouvelle
droite g\'eod\'esique d\'efinit clairement les m\^emes bouts que
l'ancienne; par
cons\'equent, leurs mots horizontaux co\"\i ncident
(Proposition~\ref{prop.memes-bouts}), et on peut
appliquer l'Affirmation~\ref{affi.fleches-verticales-facile1} \`a
$\Gamma'$. Ceci contredit le fait que $y \in  \gamma'_0 \cap
h^{-n}(\Gamma')$.
\end{proof}

\begin{proof}[\proofname\ du
Th\'eor\`eme~\ref{theoDbis} dans le cas
simple, fin]
On reprend le dis\-que $D_0$ introduit ci-dessus.
Soit $x$ un point quelconque de $D_0 \cap G$ (qui n'est pas vide
puisque $G$ rencontre $\gamma_0$).
\emph{Supposons que $x$ soit dans $\wfn(G)$}: il existe un entier
$n_0$ tel que pour tout $n \geq n_0$, $h^n(x)$ est dans l'ouvert
$O_N(\Gamma)$
\`a gauche de $\Gamma$. Alors
l'Affirmation~\ref{affi.fleches-verticales-facile2} implique que
pour tout $n \geq n_0$, le disque $h^n(D_0)$ est enti\`erement inclus
dans l'ouvert $O_N(\Gamma)$; en particulier,
tous les points de $D_0$ vont vers le Nord.
Les ensembles $G \cap \inte(D_0)$ et $F \cap \inte(D_0)$ sont des
ouverts non vides de $\wfn(G)$ et $\wfn(F)$ respectivement. D'apr\`es la
Proposition~\ref{prop.continus-minimaux} et la d\'efinition des types
dynamiques, $G$ et $F$ sont tous les
deux de type dynamique Sud-Nord. De plus, les points $x_0$ et $x_1$
vont vers le Nord.
\emph{Si $x$ est dans $\wfs(G)$}, on montre de mani\`ere sym\'etrique que
$G$ et $F$ sont de type dynamique Nord-Sud, et que $x_0$ et $x_1$ vont
vers le Sud. Ceci termine la preuve du th\'eor\`eme dans les cas
``attractif'' et ``r\'epulsif''.
\end{proof}

\begin{affi}[Compl\'ement au cas attractif]
On suppose que le mot est attractif, c'est-\`a-dire que
$$M_\leftrightarrow(x_{-1},x_2)=(\rightarrow \leftarrow).$$
Dans ce cas,  si $h(\gamma_0)$ est inclus dans	$O_S(\Gamma)$, alors
$h^n(\gamma_0)$ est encore inclus dans $O_S(\Gamma)$ pour tout $n \geq 2$.
En particulier, on a les \'equivalences
$$
\begin{array}{rcl}
h(\gamma_0) \subset O_S(\Gamma)   & \Leftrightarrow   & x_{k} \mbox{
et } x_{k+1} \mbox{ vont vers le Sud}  \\
  & \Leftrightarrow  &	G_k \mbox{ et } F_{k+1} \mbox{ sont de type
  dynamique Nord-Sud} \\
  & \Leftrightarrow  &	{M_\updownarrow}_k=\ \downarrow .
\end{array}
$$
\end{affi}
\begin{proof}
On suppose que $h(\gamma_0)$
est inclus dans  $O_S(\Gamma)$.
L'arc $\gamma_0$ est libre, mais $\gamma_0 \cup \gamma_1$ ne l'est
pas, et la fl\`eche $M_\leftrightarrow(x_{0},x_2)=\leftarrow$ indique
que c'est l'image de $\gamma_{1}$ qui rencontre $\gamma_{0}$.
 Il existe donc un point $x$ sur $\gamma_1$
tel que $[x_0 x[_\Gamma$ soit libre, mais $[x_0 x]_\Gamma$ ne le soit
pas: le point $h(x)$ est donc sur $\gamma_0$. On consid\`ere l'arc
$\alpha=[h(x) x]_\Gamma$.
Le point-cl\'e, prouv\'e plus bas, consiste \`a voir que
$$
\forall n \geq 2, \ \  h^n(\alpha) \cap \Gamma = \emptyset.
$$
En effet, on en d\'eduit que l'ensemble
$$
\bigcup_{n \geq 2} h^n(\alpha)
$$
est disjoint de $\Gamma$, or il est connexe et il rencontre
$h(\gamma_0)$, c'est donc qu'il est inclus dans $O_S(\gamma)$;
mais pour tout $n \geq 2$, $h^n(\gamma_0)$ rencontre \`a son tour cet
ensemble, et est disjoint de $\Gamma$ d'apr\`es
l'Affirmation~\ref{affi.fleches-verticales-facile1}, c'est
donc que $h^n(\gamma_0) \subset O_S(\Gamma)$.

Montrons le point-cl\'e. Soit $y$ un point de $\alpha$. Si $y$ est dans
$[h(x) x_1]_\gamma$, alors  $h^n(y)$ appartient \`a
$h^n(\gamma_0)$, et l'Affirmation~\ref{affi.fleches-verticales-facile1}
 dit que $h^n(y) \not \in
\Gamma$ pour tout $n \geq 1$. Si $y=x$, on conclut de m\^eme, puisque
$h(x) \in \gamma_0$. Il reste le cas o\`u $y$ est dans $]x_1 x [_\Gamma$.
Dans ce cas, l'arc $[x_0 y]_\Gamma$ est libre: par suite, on obtient
une nouvelle d\'ecomposition de $\Gamma$ en posant
$$
{x_1}'=y, \quad {\gamma_0}'=[x_0 {x_1}']_\Gamma, \quad
{\gamma_1}'=[{x_1}' x_2]_\Gamma,
$$
$$
{\gamma_i}'=\gamma_{i} \quad \forall i \neq 0,1 \quad \mbox{ et } \quad
{x_{i}}'=x_{i} \quad\forall i \neq 1.
$$
 On applique alors l'Affirmation~\ref{affi.fleches-verticales-facile1}
  \`a cette nouvelle
d\'ecomposition. Ceci donne notamment, pour tout $n \geq 1$,
$h^n({x_0}') \not \in \Gamma$, comme voulu.
\end{proof}

\subsection{Preuve dans le cas difficile (mot indiff\'erent)}
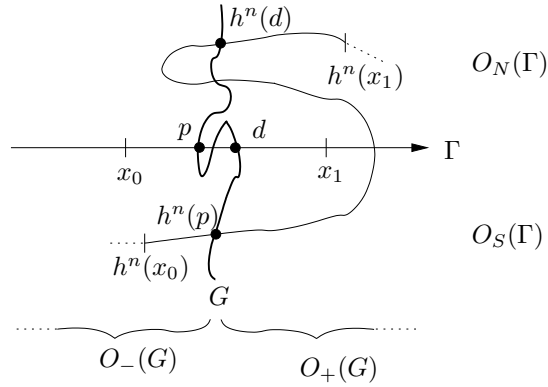
\begin{figure}[htp]\anchor{fig.image-gamma0-indif}
\centerline{\hbox{\input{fig-image-gamma0-indif2.pstex_t}}}
\caption{\label{fig.image-gamma0-indif}Disposition de l'image de
$\gamma_0$ (cas indiff\'erent)}
\end{figure}

\emph{On se place maintenant dans le cas ``indiff\'erent'', o\`u le mot
$M_\leftrightarrow(x_{-1},x_2)$ vaut $(\rightarrow \rightarrow)$ ou
$(\leftarrow \leftarrow)$.} On reprend les notations du cas simple.
Quitte \`a remplacer $h$ par $h^{-1}$, on peut supposer que
$M_\leftrightarrow(x_{-1},x_2)=(\rightarrow \rightarrow)$.
On note  $O_-(G)$  et $O_+(G)$ les ensembles compl\'ementaires
de $G$ (D\'efinition~\ref{defi.ensembles-complementaires});
rappelons que ces ensembles sont invariants (Lemme~\ref{lemm.pasdenom}).
On note $p$ le premier point et $d$ le dernier point sur $G$ quand on
parcourt la droite orient\'ee $\Gamma$ de $(-\infty)$ vers $(+\infty)$
(\cf
\figref{fig.image-gamma0-indif}). Les points $p$ et $d$ sont
situ\'es entre $x_0$ et $x_1$ sur $\Gamma$ (d'apr\`es le
Lemme~\ref{lemm.compatibilite-crg} et l'in\'egalit\'e triangulaire). Les
demi-droites
$]-\infty, p[_\Gamma$ et $]d, +\infty[_\Gamma$ sont incluses
respectivement dans  les ouverts $O_-(G)$ et $O_+(G)$.

Nous allons prouver successivement les r\'esultats suivants.
\begin{affi}\label{affi.fleches-verticales1}
Pour tout $n \geq 1$,
\begin{itemize}
\item  $h^n([x_0p]_\Gamma) \cap \Gamma=\emptyset$;
\item  $h^{-n}([d x_1]_\Gamma) \cap \Gamma=\emptyset$.
\end{itemize}
\end{affi}
\begin{affi}\label{affi.fleches-verticales2}
Soit $n \geq 1$. On a l'\'equivalence suivante:
$$
h^n([x_0 p]_\Gamma) \subset O_S(\Gamma) \quad \quad \Leftrightarrow \quad
\quad h^{-n}([d x_1]_\Gamma) \subset O_N(\Gamma).
$$
\end{affi}
\begin{coro}\label{coro.fleches-verticales2}
On a les \'equivalences:
$$
\begin{array}{rcl}
x_0  \mbox{ va vers le Sud } & \Leftrightarrow & x_1  \mbox{ vient
du Nord} \\
 & \Leftrightarrow & p \in \wfs(G) \\
 & \Leftrightarrow & d \in \wnf(G).
\end{array}
$$
\end{coro}

\begin{affi}\label{affi.fleches-verticales3}~
\begin{enumerate}
\item Si  $x_1$ vient du Nord, alors $G$ est de type dynamique
Nord-Sud;
\item Si  $x_0$ va vers le Sud, alors $F$ est de type dynamique
Nord-Sud.
\end{enumerate}
\end{affi}

\begin{proof}[\proofname\ de l'Affirmation~\ref{affi.fleches-verticales1}]
Soit $n \geq 1$, et montrons le premier point. Comme dans la preuve de
l'Affirmation~\ref{affi.fleches-verticales-facile1}, on doit montrer
que pour tout entier $k$,  $h^n([x_0p]_\Gamma) \cap \gamma_k=
\emptyset$ (\cf \figref{fig.image-gamma0-indif}).
Comme avant, ceci est vrai pour $k=0$ (libert\'e) et pour $|k|>1$
(g\'eod\'esicit\'e). Pour $k=-1$, ceci suit du fait que
$M_\leftrightarrow(x_{-1},x_1)=(\rightarrow)$.
Il reste \`a voir que $h^n([x_0 p]_\Gamma) \cap \gamma_{1}=\emptyset$.
Ceci vient du fait que l'arc $[x_0 p]_\Gamma$ est inclus dans l'ensemble
ferm\'e invariant $O_-(G)
\cup G$, tandis que  l'arc $\gamma_1=[x_1x_2]_\Gamma$ est inclus dans
$O_+(G)$.

Le deuxi\`eme point de l'affirmation se d\'emontre de fa\c{c}on
sym\'etrique.
\end{proof}

\begin{proof}[\proofname\ de l'Affirmation~\ref{affi.fleches-verticales2}]
Afin de simplifier les notations, on \'ecrit seulement la preuve pour $n=1$; la
preuve pour $n>1$ \'etant
identique. D'apr\`es l'Affirmation~\ref{affi.fleches-verticales1}, il
suffit de montrer que $h(p)$ est \`a droite de $\Gamma$  (\ie\ dans
$O_S(\Gamma)$)
 si et seulement si  $h^{-1}(d)$ est  \`a gauche de $\Gamma$ (\ie\
 dans $O_N(\Gamma)$).
Soit $\Gamma'$ la droite topologique orient\'ee image de $\Gamma$ par $h$.
Comme $h$
pr\'eserve l'orientation, le point $h^{-1}(d)$ est \`a gauche de $\Gamma$
si et seulement si le point $d$ est \`a gauche de
$\Gamma'$. Il s'agit
donc de montrer que \emph{le point $B=h(p)$ est
\`a droite de $\Gamma$ si et seulement si le point $A=d$ est \`a gauche de
$\Gamma'$.} Ceci va r\'esulter de consid\'erations de topologie plane.
On note (\cf \figref{fig.gamma-prime}):
\begin{align*}
\Gamma^-&=\ ]{-}\infty,  A]_\Gamma,& \Gamma^+&=[A,+\infty[_\Gamma,\\
{\Gamma'}^-&=\ ]{-}\infty,B]_{\Gamma'},& {\Gamma'}^+&=[B,+\infty[_{\Gamma'}.
\end{align*}
Puisque $\Gamma$ est une droite topologique, d'apr\`es le
th\'eor\`eme de Schoenflies, on peut supposer que c'est une droite
euclidienne horizontale orient\'ee de la gauche vers la droite, comme
sur la \figref{fig.gamma-prime}.
\begin{figure}[htp]\anchor{fig.gamma-prime}
\centerline{\hbox{\input{fig-gamma-prime.pstex_t}}}
\caption{\label{fig.gamma-prime}Disposition de la droite g\'eod\'esique
$\Gamma$ et de son image $\Gamma'$}
\end{figure}
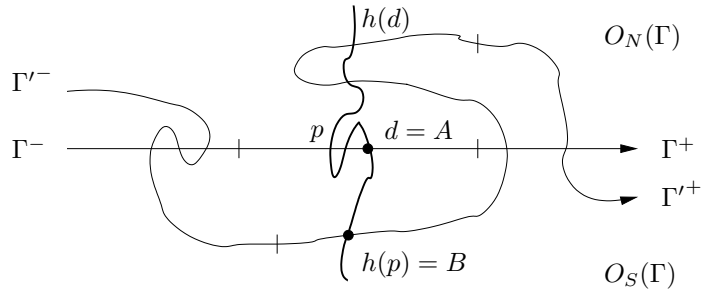

Montrons d'abord les  propri\'et\'es suivantes.
\begin{enumerate}
\item $\Gamma^+ \cap {\Gamma'}^+ \neq \emptyset$;
 \item $\Gamma^- \cap {\Gamma'}^- \neq \emptyset$;
\item $\Gamma^- \cap {\Gamma'}^+ = \emptyset$;
\item $\Gamma^+ \cap {\Gamma'}^- = \emptyset$.
\end{enumerate}
Les deux premi\`eres propri\'et\'es viennent du fait que $\Gamma$ est une
droite g\'eod\'esique (deuxi\`eme point du
Lemme~\ref{lemm.iteres-arc-geo}).
Pour la troisi\`eme, on peut m\^eme montrer que les demi-droites
$]{-}\infty,x_1]_\Gamma$ et $h([x_0,+\infty[_\Gamma)$ sont disjointes,
en  utilisant encore le caract\`ere g\'eod\'esique (premier point du
m\^eme lemme), ainsi que l'hypoth\`ese
$M_\leftrightarrow(x_{-1},x_2)=(\rightarrow \rightarrow)$ et la
d\'efinition des  fl\`eches. Enfin, la derni\`ere propri\'et\'e provient
simplement   des inclusions $\Gamma^+ \subset O_+(G) \cup \{d\}$ et
${\Gamma'}^- \subset O_-(G) \cup \{h(p)\}$ (rappelons que $h(p) \neq d$,
puisque $p$ et $d$ appartiennent \`a un m\^eme arc libre).

\emph{Supposons maintenant que $B$ est
\`a droite de ${\Gamma}$, et montrons que $A$ est \`a gauche de
${\Gamma'}$}. Remarquons tout d'abord que $A$ n'est pas sur
$\Gamma'$ (d'apr\`es les deux derni\`eres propri\'et\'es).
Soit $]C D[_{\Gamma'}$ la composante connexe de $\Gamma' \setminus
\Gamma$ contenant $B$, o\`u $C$ et $D$ sont choisis de fa\c{c}on \`a
ce que
$C$ soit situ\'e avant $D$ sur la droite topologique orient\'ee
$\Gamma'$; l'existence de ces deux
points vient des deux premi\`eres propri\'et\'es ci-dessus (\cf
\figref{fig.B-a-droite}).
D'apr\`es les deux derni\`eres propri\'et\'es, le point $C$ est sur
$\Gamma^-
\cap {\Gamma'}^-$, et le point $D$ sur	$\Gamma^+
\cap {\Gamma'}^+$ (et par cons\'equent $C$ est \'egalement situ\'e
avant $D$ sur $\Gamma$).
 Par d\'efinition, l'arc ouvert $]C D[_{\Gamma'}$ est disjoint de
$\Gamma$, et la courbe $[C D]_{\Gamma'} \cup [CD]_\Gamma$ est une courbe
de Jordan (\ie elle est hom\'eomorphe au cercle). Soit $U$ le disque
topologique ouvert bord\'e cette courbe.
 Remarquons  que puisque la droite
$\Gamma$ est disjointe de $]C D[_{\Gamma'}$, elle est aussi disjointe de
l'ouvert $U$.
D'apr\`es le th\'eor\`eme de Schoenflies, sachant que $B$ est \`a
droite de $\Gamma$, la situation est  hom\'eomorphe au
dessin de la \figref{fig.B-a-droite}: on peut supposer que
$U$ est un demi-disque euclidien, et $U$ est situ\'e \`a gauche de l'arc
orient\'e $[C D]_{\Gamma'}$ et \`a droite de l'arc orient\'e $[C
D]_\Gamma$.
\begin{figure}[htp]\anchor{fig.B-a-droite}
\centerline{\hbox{\input{fig-B-a-droite.pstex_t}}}
\caption{\label{fig.B-a-droite}Si $B$ est \`a droite de $\Gamma$\ldots{}}
\end{figure}
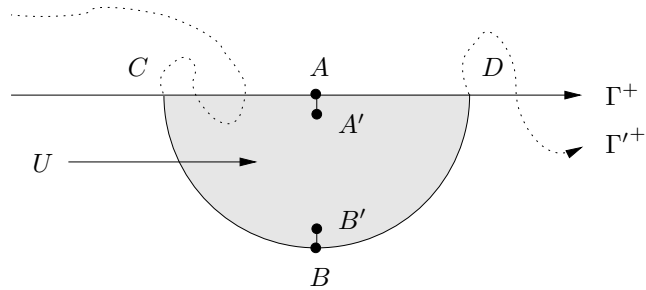

Sur cette figure, choisissons un point $A'$ dans $U$, suffisamment proche
de $A$
pour que le segment euclidien $[AA']$ soit disjoint de $\Gamma'$.
 De m\^eme, on choisit	un point $B'$ dans $U$,  suffisamment proche
 de $B$
pour que le segment euclidien $]BB']$ soit disjoint de $\Gamma'$.
 Puisque $B'$ est \`a gauche de
$\Gamma'$, pour que le point $A$ soit situ\'e \`a gauche de $\Gamma'$ il
suffit qu'il existe un arc reliant $A'$ \`a $B'$ en \'evitant
$\Gamma'$. Nous allons trouver un tel arc dans $U$.

Supposons que cet arc n'existe pas. Alors l'ensemble $\Gamma' \cap U$,
qui est une
partie ferm\'ee de $U$, s\'epare $A'$ et $B'$ dans $U$. L'ouvert $U$
\'etant hom\'eomorphe au plan, on peut appliquer
 le Th\'eor\`eme~\ref{theo.separation-plane} de l'appendice: il existe
une composante connexe de $\Gamma' \cap U$ qui s\'epare encore $A'$ et
$B'$. Cette composante connexe est un intervalle de ${\Gamma'}^-$ ou bien
de ${\Gamma'}^+$; elle est disjointe de $[C D]_{\Gamma'}$. Puisqu'elle
s\'epare
$A'$ et $B'$ dans $U$, et ne rencontre ni $[AA']$ ni $[BB']$, son
adh\'erence doit rencontrer chacun des deux arcs
$[CA]_\Gamma$ et $[AD]_\Gamma$. Ceci contredit l'une des deux
propri\'et\'es~3 et~4.

Nous avons donc montr\'e  que si $B$ est
\`a droite de ${\Gamma}$, alors $A$ est \`a gauche de
${\Gamma'}$. On montre de m\^eme que  si $B$ est
\`a gauche de ${\Gamma}$, alors $A$ est \`a droite de
${\Gamma'}$. Ce qui conclut.
\end{proof}

\begin{proof}[\proofname\ du Corollaire~\ref{coro.fleches-verticales2}]
Si $x_0$ va vers le Sud, alors il existe un entier $n_0$ tel que pour
tout $n \geq n_0$, $h^n(x_0) \in O_S(\Gamma)$. Mais alors d'apr\`es
l'Affirmation~\ref{affi.fleches-verticales1}, pour $n \geq n_0$, on a
$h^n([x_0p]_\Gamma) \subset O_S(\Gamma)$, et en particulier le point
$p$ va aussi vers le Sud. Et d'apr\`es
l'Affirmation~\ref{affi.fleches-verticales2}, toujours pour $n \geq
n_0$, on a $h^{-n}([d x_1]_\Gamma) \subset O_N(\Gamma)$, par
cons\'equent les points $d$ et $x_1$ viennent du Nord.
\end{proof}

\begin{proof}[\proofname\ de l'Affirmation~\ref{affi.fleches-verticales3}]
Pour prouver cette affirmation, nous com\-men\-\c{c}ons par remplacer
$\Gamma$
par une autre droite g\'eod\'esique $\Gamma'$, qui d\'efinira les
m\^emes bouts que $\Gamma$. Par cons\'equent, $G$ sera encore un bord
de composante de Reeb associ\'ee \`a $\Gamma'$
(Proposition~\ref{prop.memes-bouts}), et on pourra appliquer le
Corollaire~\ref{coro.fleches-verticales2} \`a cette nouvelle situation.

Puisque les ensembles $\wnn(G)$ et $\wss(G)$ sont ferm\'es
(Proposition~\ref{prop.continus-minimaux}) et disjoints, il existe un
voisinage $V$
du point $d$ qui ne rencontre pas simultan\'ement ces deux
ensembles.

En \'epaississant l\'eg\`erement l'arc $[d x_1]_\Gamma$ (\`a
l'aide du th\'eor\`eme de Schoenflies),
on choisit maintenant un disque topologique ferm\'e $D$ tel
que (\cf \figref{fig.gammaprime})
\begin{itemize}
\item $\Gamma \cap D$ est un sous-arc de $\Gamma$ dont $x_1$ est une
extr\'emit\'e;
\item $d \in \inte(D)$;
\item $D \cap G \subset V$;
\item $D$ est libre.
\end{itemize}
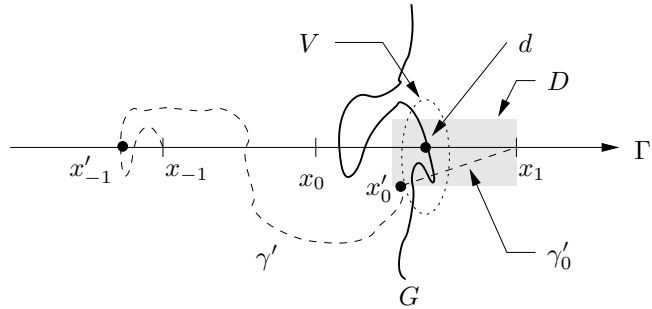
\begin{figure}[htp]\anchor{fig.gammaprime}
\centerline{\hbox{\input{fig-gammaprime.pstex_t}}}
\caption{\label{fig.gammaprime}Construction de la g\'eod\'esique
$\Gamma'$ (I)}
\end{figure}
Rappelons que l'ensemble $G$ est inclus dans l'adh\'erence de
$B_1(x_{-1})$ (par
d\'e\-fi\-ni\-tion). Par construction, $D$ contient donc des points \`a
$h$--distance 1 de $x_{-1}$, autrement dit
il existe un arc libre $\gamma'$ joignant $x_{-1}$ \`a un
point de $D$; et on peut supposer, quitte \`a raccourcir cet arc,
 que $\gamma' \cap D$ contient un
unique point qu'on appelle $x'_0$. Remarquons que $\gamma'$  est
disjoint de $[x_1, +\infty[_\Gamma$, car $(x_k)$ est une suite
g\'eod\'esique.
On pose $x'_1=x_1$, et on choisit
un arc $\gamma'_0$ joignant $x'_0$ \`a $x'_1$ dans $D$ (cet arc est
donc libre, et ne rencontre $\gamma'$ qu'au point $x'_0$).
 On appelle $x'_{-1}$ le dernier point d'intersection de
$\gamma'$ avec $]-\infty, x_{-1}]_\Gamma$ lorsqu'on parcourt $\gamma'$
depuis $x_{-1}$ vers $x'_0$. Remarquons encore que ce point est situ\'e
sur $]x_{-2}, x_{-1}]_\Gamma$, toujours \`a cause du
caract\`ere g\'eod\'esique de $(x_k)$.
On pose alors (\cf \figref{fig.gammaprime2})
\begin{align*}
\gamma'_{-2}&=[x_{-2} x'_{-1}]_\Gamma,\qquad\qquad
\gamma'_{-1}=[x'_{-1} x'_0]_{\gamma'} \\
\Gamma'&=\ ]{-}\infty, x'_{-2}]_\Gamma \ \cup \gamma'_{-2}  \ \cup \
\gamma'_{-1}  \ \cup  \ \gamma'_0
 \ \cup  \ [x_1,+\infty[_\Gamma.
\end{align*}
Remarquons que par d\'efinition du point $x'_{-1}$,  ce point est
l'unique point d'in\-ter\-sec\-tion des arcs $\gamma'_{-2}$ et $\gamma'_{-1}$.

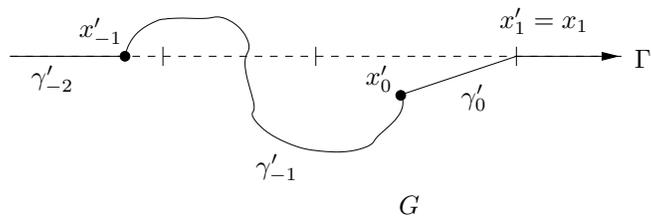
\begin{figure}[htp]\anchor{fig.gammaprime2}
\centerline{\hbox{\input{fig-gammaprime2.pstex_t}}}
\caption{\label{fig.gammaprime2}Construction de la g\'eod\'esique
$\Gamma'$ (II)}
\end{figure}

\emph{Dans cette situation, la droite $\Gamma'$, munie de la
d\'ecomposition
$$
\Gamma'=\cdots * \gamma_{-3} * \gamma'_{-2} * \gamma'_{-1}* \gamma'_0 *
\gamma_1 * \cdots
$$
 est une droite g\'eod\'esique, qui d\'efinit les m\^emes
 bouts que $\Gamma$.}
En effet,  on v\'erifie que la courbe est sans point double par
 construction; le fait	qu'il s'agit bien d'une droite g\'eod\'esique, et
 qu'elle d\'efinit les m\^emes bouts que $\Gamma$, est imm\'ediat.

En appliquant le Lemme~\ref{lemm.ordre-crg} sur le franchissement des
composantes de Reeb, on voit que $\Gamma' \cap
G \subset \gamma'_0$, donc $\Gamma' \cap G \subset D \cap G \subset   V$.
Comme pour $\Gamma$, on d\'efinit les points $p'$ et $d'$ comme
\'etant les premier et dernier points sur $G$ lorsqu'on parcourt
$\Gamma'$ de $(-\infty)$ vers $(+\infty)$. Ces deux
points sont dans $V$.

\emph{On suppose maintenant que le point $x_1=x'_1$ vient du Nord}. On
applique alors le Corollaire~\ref{coro.fleches-verticales2} \`a
la droite g\'eod\'esique $\Gamma'$; on obtient que le point $p'$
est dans $\wfs(G)$ et le point $d'$ dans $\wnf(G)$.
Le choix de $V$ emp\^eche alors qu'on ait simultan\'ement  $p' \in
\wss(G)$ et $d' \in \wnn(G)$. On en conclut
que l'un des deux points au moins est dans $\wns(G)$. Ceci prouve que
le type dynamique de $G$ est Nord-Sud.

La deuxi\`eme partie de l'Affirmation~\ref{affi.fleches-verticales3}
provient d'une construction enti\`e\-re\-ment sym\'etrique.\footnote{Plus
pr\'ecis\'ement, on consid\`ere
l'hom\'eomorphisme $h'=h^{-1}$, et la droite g\'eod\'esique $\Gamma'$
obtenue \`a partir de la droite topologique $\Gamma$ en renversant
l'orientation  et en choisissant la
d\'ecomposition dont la suite des sommets est $(x'_k:=x_{1-k})$. On a
alors $M^{h'}_\leftrightarrow(x'_{-1},x'_2)=(\rightarrow
\rightarrow)$, tandis que les r\^oles de $x_0$ et de $x_1$ sont
invers\'es,
ainsi que ceux de $N$ et de $S$, et ceux de $F$ et de $G$. En
appliquant le point 1 de
l'affirmation \`a cette situation sym\'etrique, on obtient alors le
point 2.}
\end{proof}

\begin{proof}[\proofname\ du
Th\'eor\`eme~\ref{theoDbis} dans le cas
difficile, fin]
 Le point
$x_0$ va vers le Sud ou vers le Nord
(Corollaire~\ref{coro.partition-dynamique}).
On suppose par exemple qu'il va vers le Sud, l'autre cas se traitant
de mani\`ere identique en inversant les r\^oles des points $N$ et
$S$.
Dans ce cas, d'apr\`es le Corollaire~\ref{coro.fleches-verticales2},
le point
$x_1$ vient du Nord.
D'apr\`es l'Affirmation~\ref{affi.fleches-verticales3}, $F$ et $G$
sont tous deux de type dynamique Nord-Sud.
 Ceci ach\`eve la preuve du th\'eor\`eme.
\end{proof}

\section{Dynamique dans les continus minimaux}
\label{sec.continus-minimaux}
Dans cette section, nous prouvons la
Proposition~\ref{prop.continus-minimaux}.
\emph{Les Sections~\ref{ssec.enon} et
\ref{ssec.demo-continus} sont ind\'ependantes du reste du texte.}
Rappelons qu'un \emph{continu} est un espace topologique (m\'etrisable)
compact et connexe.

Le r\'esultat principal de cette section dit qu'un bord $F$ de composante
de Reeb  (associ\'e \`a une suite
g\'eod\'esique infinie $(x_k)$) contient soit des points qui vont du
Nord au
Sud, soit des points qui vont du Sud au Nord, mais pas les deux \`a
la fois. Cette propri\'et\'e est avant tout la cons\'equence de la
minimalit\'e de $F$ (en tant qu'ensemble ferm\'e, connexe,
s\'eparant les deux bouts de $(x_k)$). Ceci va nous entra\^\i ner
dans des consid\'erations de dynamique topologique abstraite, o\`u
les seuls
souvenirs de la dynamique plane seront deux propri\'et\'es
techniques, cons\'equences du lemme de Franks, et prouv\'ees dans la
Section~\ref{ssec.continus-prelim}.
De ce point de vue,  cette partie appara\^\i t comme \'etant  le
prolongement de l'article~\cite{lero9}: dans un contexte plus
g\'en\'eral (qui correspond aux trois premi\`eres propri\'et\'es
de l'hypoth\`ese \textrm{(H)} ci-dessous, augment\'ees d'une
hypoth\`ese de connexit\'e locale qui n'\'etait pas essentielle),
on montrait l'existence d'au moins \emph{un} point dont l'orbite traverse
$F$ du Nord au Sud ou du Sud au Nord. L'hypoth\`ese de minimalit\'e et
les ``souvenirs du plan''
que nous ajoutons ici permettent donc de renforcer s\'erieusement les
r\'esultats pr\'eliminaires de~\cite{lero9}, r\'esultats que  nous
r\'eutilisons ici
 (sous la forme du Lemme~\ref{lemm.CRAS}). Je dois \`a Patrice Le
Calvez l'id\'ee d'utiliser des compacts minimaux dans le contexte de
l'article~\cite{lero9}.

\subsection{Pr\'eliminaires: deux cons\'equences du lemme de Franks}
\label{ssec.continus-prelim}
Rappelons que la \emph{limite sup\'erieure} d'une suite de parties
 d'un espace topologique $X$ est
l'ensemble des points de $X$ dont tout voisinage rencontre une infinit\'e
d'\'el\'ements de la suite. Nous aurons \'egalement besoin du concept
suivant, qui est une bonne alternative \`a la notion de \emph{composantes
connexes par arcs} dans les espaces o\`u celles-ci sont trop petites
pour \^etre utiles (par exemple, quand elles sont r\'eduites \`a des
points).
\begin{defi}\label{defi.ccc}
Soit $X$ un espace topologique. On d\'efinit une
 relation d'\'equ\-i\-val\-ence sur $X$, en d\'ecidant que
deux points sont \'equivalents si ils appartiennent \`a une m\^eme partie
compacte et connexe de $X$. On appellera
 \textit{composantes compactement connexes} les classes d'\'equivalence
pour cette relation.
\end{defi}
Ces d\'efinitions sont utilis\'ees par W. Daw dans l'article
\cite{daw1}. Elles apparaissent d\'ej\`a dans le livre de
K. Kuratowski, \cite{kura1} (paragraphe 42, VIII; les composantes
compactement connexes
y sont appel\'ees des \emph{constituants}.)
Dans les bons espaces, les composantes compactement
connexes sont ``assez grandes'',  elles ``vont \`a l'infini''
 (voir l'appendice, Corollaire~\ref{coro.ccc}).
Notons aussi qu'une composante compactement connexe n'est pas
n\'ecessairement  ferm\'ee.

On consid\`ere un hom\'eomorphisme de Brouwer $h$.
Nous d\'emontrons les deux \'enonc\'es suivants.
\begin{lemm}\label{lemm.limitesup}
Si $C$ est une partie  compacte du plan, alors la limite sup\'erieure
de la suite
$(h^n(C))_{n \geq 0}$ ne contient pas $C$.
\end{lemm}
\begin{lemm}\label{lemm.limitesup2}
Si $E$ est une partie compactement connexe du plan, libre par $h$,
et si $E \not \subset \adhe(h(E))$, alors  la limite sup\'erieure de
la suite
$(h^n(\adhe(E)))_{n \geq 0}$ ne contient pas $E$.
\end{lemm}

\begin{proof}[\proofname\ du Lemme~\ref{lemm.limitesup}]
Pour donner l'id\'ee de la preuve, commen\c{c}ons par le cas simple
o\`u $C$
est un disque topologique ferm\'e. On peut alors \'ecrire $C$ comme la
r\'eunion d'un nombre fini de disques topologiques ferm\'es $D_i$, libres,
d'int\'erieurs deux \`a deux disjoints. On sait que la suite
$(h^n(C))$ converge (Lemme~\ref{lemm.convergence}).
Si $C$ \'etait inclus dans la  limite
 de cette suite, alors pour tout disque
$D_i$, il existerait un disque $D_j$ et un entier $n_i>0$ tel que
$h^{n_i}(\inte(D_j))$ rencontre $\inte(D_i)$. Puisque les disques
$D_i$ sont en nombre finis, on en d\'eduirait l'existence d'une cha\^\i ne
de disques p\'eriodique, ce qu'interdit le lemme de Franks.

Voici la preuve dans le cas g\'en\'eral.\footnote{Par rapport au cas
simple, la complication vient du fait
qu'on ne sait pas recouvrir $C$ facilement par des disques
topologiques ferm\'es d'int\'erieurs deux \`a deux disjoints, et dont les
int\'erieurs rencontrent $C$.}
 Soit $C$ v\'erifiant les hypoth\`eses
du lemme. On consid\`ere un  quadrillage d'un carr\'e du plan contenant
$C$ (comme sur la
\figref{fig.quadrillage}). Pour tout point $x$ de $C$, on note
$D(x)$ la r\'eunion des petits carr\'es ferm\'es contenant $x$: $D(x)$
est un
disque topologique ferm\'e  r\'eunion
de $1$, $2$ ou $4$ petits carr\'es, et on peut choisir le quadrillage
pour que tous les ensembles $D(x)$ soient libres.
Remarquons que $x$ appartient  \`a l'int\'erieur de $D(x)$. D'autre part,
si $y$
est un point de $C$ qui n'est pas dans $D(x)$, alors les deux disques
$D(y)$ et $D(x)$ sont d'int\'erieurs disjoints.
\begin{figure}[htp]\anchor{fig.quadrillage}
\centerline{\hbox{\input{fig-quadrillage.pstex_t}}}
\caption{\label{fig.quadrillage}approximation de $C$ par des carr\'es
libres}
\end{figure}
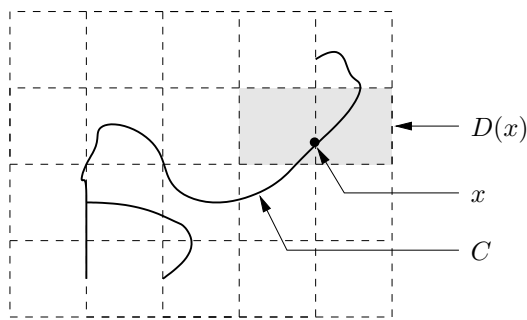

On suppose que $C$ est inclus dans la limite sup\'erieure $C_\infty$
de
$(h^n(C))_{n \geq 0}$, et on va construire par r\'ecurrence une cha\^\i
ne de disques (pour
$h^{-1}$). Soit $x_1$ un point quelconque de $C$. Puisque $x_1 \in
\inte(D(x_1)) \cap C_\infty$, il existe un entier $n_1>0$ et un point
$x_2$ de $C$
tels que $h^{n_1}(x_2)$ appartient \`a l'int\'erieur de $D(x_1)$. En
it\'erant
la construction, on trouve une suite $(x_i)_{i \geq 1}$ de points de
$C$ et une suite d'entiers $n_i >0$ telles que pour
tout $i \geq 1$, $h^{n_i}(x_{i+1})$ appartient \`a l'int\'erieur de
$D(x_i)$.
Puisque le nombre de petits carr\'es du quadrillage qui rencontrent $C$
est fini, il existe
un entier $i  \geq 1$ tel que $x_{i+1}$ appartient \`a l'un des
disques $D(x_j)$ pour $1 \leq j \leq i$. On appelle $i_2$ le plus
petit des tels entiers $i$, et $i_1 \leq i_2$ un entier tel que
$x_{i_2+1}$ appartient \`a $D(x_{i_1})$. On a
\begin{itemize}
\item les disques $D(x_1),\ldots, D(x_{i_2})$ sont d'int\'erieurs deux \`a
deux disjoints, gr\^ace au choix de $i_2$;
\item pour $1 \leq i \leq i_2 -1$, l'int\'erieur de
$h^{n_i}(D(x_{i+1}))$ rencontre l'int\'erieur de $D(x_i)$
 (puisque $h^{n_i}(x_{i+1}) \in \inte(D(x_i))$);
\item  l'int\'erieur de $h^{n_{i_2}}(D(x_{i_1}))$ rencontre l'int\'erieur
de $D(x_{i_2})$
(puisque $D(x_{i_1})$ contient	$x_{i_2+1}$ et $h^{n_{i_2}}(x_{i_2+1})
\in \inte(D(x_{i_2}))$).
\end{itemize}
Autrement dit, la suite $(D(x_{i_2}), D(x_{i_2-1}), \dots, D(x_{i_1}))$
 est une cha\^\i ne de disques p\'eriodique,
ce qui contredit le lemme de Franks.\footnote{Ceci est le seul endroit
du texte o\`u nous utilisons
l'absence de cha\^\i ne p\'eriodique de plus de deux \'el\'ements.}
\end{proof}

\begin{proof}[\proofname\ du Lemme~\ref{lemm.limitesup2}]
Soit $E$ v\'erifiant les hypoth\`eses du lemme.
Puisque $E \not \subset \adhe(h(E))$,  il
existe	un  disque topologique ferm\'e libre $D_1$, disjoint de
$h(E)$, dont l'int\'erieur rencontre $E$ (\cf \figref{fig.chaine2}).
On raisonne par l'absurde, en supposant que la limite sup\'erieure
$E_\infty$ de la suite $(h^n(\adhe(E)))_{n \geq 0}$ contient $E$. En
particulier, $E_\infty$ rencontre l'int\'erieur de $D_1$, et il existe
un entier
$n \geq 100$ tel que $h^n(E)$ rencontre l'int\'erieur de $D_1$.
Puisque $E$ est
compactement connexe et rencontre  $\inte(D_1)$ et $h^{-n}(\inte(D_1))$,
il existe une
partie compacte et connexe $C$ de $E$  qui rencontre
encore $\inte(D_1)$ et $h^{-n}(\inte(D_1))$. On a
\begin{enumerate}
\item $C$ est libre (car $E$ l'est);
\item $h(C)$ est disjoint de $D_1$ (car $h(E)$ l'est).
\item $h(\inte(D_1))$ rencontre $h(C)$;
\item $h^{n-1}(h(C))$ rencontre $\inte(D_1)$.
\end{enumerate}
\begin{figure}[htp]\anchor{fig.chaine2}
\centerline{\hbox{\input{fig-chaine2.pstex_t}}}
\caption{\label{fig.chaine2}Construction d'une cha\^\i ne p\'eriodique \`a
deux disques}
\end{figure}
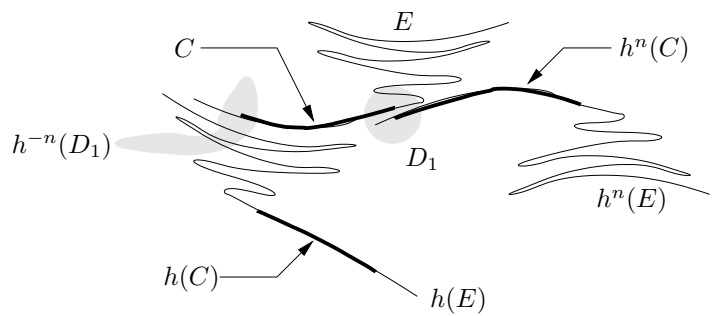

Si (par miracle!) $C$ est un disque topologique ferm\'e, alors la suite
$(D_1, h(C))$ forme une cha\^\i ne de disques p\'eriodique, contredisant
le lemme de Franks~\ref{lemm.franks}, et la preuve est termin\'ee.
Dans le cas g\'en\'eral, on remplace $C$ par un disque $D_2$, de la façon
suivante. On utilise les propri\'et\'es 3 et 4 ci-dessus pour trouver un
point $x_1$ de l'int\'erieur de $D_1$ tel que le point $y_1:=h(x_1)$ est
dans $h(C)$, et un point $x_2$ de $h(C)$ tel que le point
$y_2:=h^{n-1}(x_2)$ est dans l'int\'erieur de $D_1$. Puisque $h(C)$ est un
compact connexe libre et disjoint de $D_1$, il est inclus dans un
ouvert connexe libre et disjoint de $D_1$ (par exemple un
$\epsilon$--voisinage de $C$), et il suffit de choisir un disque
topologique ferm\'e $D_2$ inclus dans cet ouvert, et contenant les points
$y_1$ et $x_2$ dans son int\'erieur. La suite $(D_1, D_2)$ forme alors
une cha\^\i ne de disques p\'eriodique, et on aboutit \`a la m\^eme
contradiction que dans le cas miraculeux.
\end{proof}

\subsection{Hypoth\`eses}
\label{ssec.hypo}
On se place maintenant dans un espace topologique abstrait $\widehat
F$, et on consid\`ere  un hom\'eomorphisme $h$ de
$\widehat F$. On note \textrm{(H)} l'ensemble des hypoth\`eses suivantes.
\begin{enumerate}
\item \label{item.topologie} l'espace topologique $\widehat F$ est
compact, connexe, m\'etrisable;
\item \label{item.points-fixes} l'hom\'eomorphisme $h\co\widehat F
\rightarrow \widehat F$ a
exactement deux points fixes, not\'es $N$ et $S$; on pose alors
$F=\widehat F \setminus\{N,S\}$;
\item\label{item.errants} tous les  points de l'ensemble $F$ sont errants;
\item\label{item.minimalite} ({\sl propri\'et\'e de minimalit\'e\/})\qua
il n'existe pas de sous-ensemble strict $L$ de
$\widehat F$ compact, connexe, et contenant $N$ et $S$;
\item\label{item.souvenirs} ({\sl souvenirs du plan\/})\qua
\begin{enumerate}
\item si $C$ est une partie  compacte de $F$, alors la limite sup\'erieure
de la suite
$(h^n(C))_{n \geq 0}$ ne contient pas $C$;
\item si $E$ est une partie compactement connexe de $F$, libre par $h$,
et si $E \not \subset \adhe(h(E))$, alors  la limite sup\'erieure de
$(h^n(\adhe(E)))_{n \geq 0}$ ne contient pas $E$.
\end{enumerate}
\end{enumerate}

\begin{affi}\label{affi.hypotheses}
Soit $(x_k)_{k \in \bbZ}$ une suite g\'eod\'esique, et soit $F$ un bord
de composante de Reeb  associ\'ee \`a $(x_k)$. On
munit l'ensemble $\widehat F= F \cup \{N,S\}$ de la topologie de la
compactification Nord-Sud associ\'ee \`a $[(x_k)]$.
Alors l'espace $\widehat F$ et
l'hom\'eomorphisme $h$ prolong\'e \`a $\widehat F$
 v\'erifient les hypoth\`eses {\rm(H)}.
\end{affi}

\proof
On consid\`ere une droite g\'eod\'esique $\Gamma$ d\'efinissant
les m\^emes bouts que la suite $(x_k)$. Nous v\'erifions successivement
les diff\'erents points de l'hypoth\`ese \textrm{(H)}.
\begin{enumerate}
\item La compacit\'e de $\widehat F$ suit facilement de la
D\'efinition~\ref{defi.compactification}. Remarquons en particulier que
$\widehat F$ est un espace topologique s\'epar\'e.

Puisque $F$ est connexe, pour
voir que $\widehat F$ l'est
\'egalement, il suffit de montrer la densit\'e du premier dans le second.
Par l'absurde, si $F$ n'est pas dense dans $\widehat F$, alors il
existe un voisinage de $N$ ou de $S$  dans $\widehat F$ qui est
disjoint de
$F$; supposons qu'il s'agisse de $N$. Par d\'efinition de la topologie
sur $\widehat F$, ceci signifie que l'ensemble $F \cap \adhe(O_N(\Gamma))$
est
compact.  Ceci contredit le fait que $F$  s\'epare les deux bouts de
$\Gamma$ (Lemme~\ref{lemm.pasdenom}).

On peut v\'erifier que l'espace $\widehat F$  poss\`ede une base
 d\'enombrable d'ouverts. On a vu qu'il est aussi compact (et
 notamment s\'epar\'e).
 On en d\'eduit  la m\'etrisabilit\'e \`a l'aide des th\'eor\`eme
 g\'en\'eraux (voir par exemple le
corollaire 2.59 du chapitre 2 de \cite{hock1}).

\item Les points 2\ldots{}

\item \ldots{}et 3 de l'hypoth\`ese \textrm{(H)} sont clairement v\'erifi\'es
(l'errance vient du
Corollaire~\ref{coro.erre}).

\item Soient $x$ et $y$ deux sommets de $\Gamma$ tels que $F$ soit un
bord d'une des composantes de Reeb minimales associ\'ees \`a $(x,y)$.
 Soit $L$ un sous-ensemble strict de $\widehat F$, compact, contenant $N$
et $S$.  D'apr\`es le Lemme~\ref{lemm.topologie-crg} sur la topologie
des ensembles $F_i(x,y)$, $L$ ne s\'epare
pas les points $x$ et $y$. Comme $L$ est compact, les composantes
connexes du compl\'ementaire de $L$ dans le plan sont ouvertes, donc
connexes par arcs; il existe alors un
arc $\gamma'$ reliant $x$ \`a $y$ et \'evitant $L$.
On en d\'eduit facilement l'existence d'une droite topologique
$\Gamma'$, incluse dans $\Gamma \cup \gamma'$, et disjointe de
$L$. Comme $\Gamma'$ s\'epare les points $N$ et $S$ dans $\bbR^2
\cup\{N,S\}$, l'ensemble $L$ ne peut pas \^etre connexe. Ceci prouve
le quatri\`eme point.

\item Ces deux propri\'et\'es correspondent aux
Lemmes~\ref{lemm.limitesup} et~\ref{lemm.limitesup2}.\qed
\end{enumerate}

Sous les hypoth\`eses \textrm{(H)}, on d\'efinit les ensembles
$\wnf(F)$, \textit{etc.} comme indiqu\'e
\`a la Section~\ref{ssec.partition-dynamique}.
Rappelons que l'errance des points implique que tout point va vers le
Nord ou le Sud, et vient du Nord ou du Sud.

\subsection{\'Enonc\'es}
\label{ssec.enon}
Le r\'esultat principal de cette section est une version abstraite de
la Proposition~\ref{prop.continus-minimaux}; jointe \`a
l'Affirmation~\ref{affi.hypotheses}, cette version abstraite
 entra\^\i nera imm\-\'e\-di\-at\-e\-ment la
Proposition~\ref{prop.continus-minimaux}.
\renewcommand{\thepropbis}{\ref{prop.continus-minimaux}-bis}
\begin{propbis}
Sous les hypoth\`eses {\rm(H)},
\begin{itemize}
\item l'un des deux ensembles
 $\wsn(F)$ et $\wns(F)$
 est vide;
\item l'autre est un ouvert connexe dense dans $F$;
\item les ensembles $\wnn(F)$ et $\wss(F)$ sont ferm\'es.
\end{itemize}
\end{propbis}
Voici les \'enonc\'es des lemmes de la d\'emonstration. Les preuves se
trouvent dans la section suivante.
\begin{lemm}
\label{lemm.CRAS}
On se place sous les hypoth\`eses~\ref{item.topologie},
\ref{item.points-fixes} et \ref{item.errants} de~{\rm(H)}.
Supposons que pour tous voisinages $V_N$ de $N$ et $V_S$ de $S$, il
existe un  entier positif $n$
tel que $h^n(V_S)$ rencontre $V_N$.
Alors l'ensemble $\wsn(F)$ est non vide.
\end{lemm}

\begin{lemm}\label{lemm.monodynamie}
 Sous les hypoth\`eses~{\rm(H)},
toute composante compactement connexe de $F$ est enti\`erement
incluse dans l'un des quatre ensembles $\wns(F)$, $\wsn(F)$,
$\wnn(F)$, $\wss(F)$.
\end{lemm}

\begin{lemm}\label{lemm.composante-NS}
On se place sous les hypoth\`eses~{\rm(H)}.
Soit $E$ une composante compactement connexe de $F$ incluse dans
 $\wsn(F)$.
\begin{itemize}
\item Si $E$ est invariante par $h$, alors $E$ est dense dans $\widehat
F$;
\item dans le cas contraire, $E$ est libre, et $E$ est incluse dans
l'adh\'erence de $h(E)$ ou dans celle de $h^{-1}(E)$. De plus, l'union
de tous
les it\'er\'es de $E$ est un ensemble connexe, dense dans $\widehat F$,
et invariant par $h$.
\end{itemize}
\end{lemm}

\begin{coro}\label{coro.NS}
Si l'ensemble $\wsn(F)$ n'est pas vide, il est dense dans $F$ et connexe.
\end{coro}

\subsection{Preuves}
\label{ssec.demo-continus}
On note ${\cal K}(\widehat F)$ l'espace des parties compactes de
$\widehat F$, muni de la m\'etrique de Hausdorff. On rappelle que
c'est un espace m\'etrique compact.

\begin{proof}[\proofname\ du Lemme~\ref{lemm.CRAS}]
La preuve est essentiellement contenue dans la pr\-euve de la
Proposition~6 de l'article~\cite{lero9}. En voici une version
abr\'eg\'ee.
L'hy\-po\-th\-\`ese entra\^{\i}ne qu'il existe une suite $(x_p)_{p \geq 0}$
d'\'el\'ements de $F$ convergeant vers $S$, et une suite d'entiers
positifs $(n_p)$ telle que la suite $(h^{n_p}(x_p))_{p \geq 0}$
converge vers $N$. Soit
$$K_p=\{x_p, h(x_p), \dots, h^{n_p}(x_p)\}.$$
Quitte \`a extraire, on peut supposer que la suite $(K_p)$ converge dans
l'espace ${\cal K}(\widehat F)$.
On note $L$ sa limite. La fin de la preuve consiste \`a montrer que $L$
contient
une orbite allant de $S$ \`a $N$.
Pour cela, on montre d'abord que l'ensemble $L$ est invariant par $h$,
et qu'il est localement fini, au sens o\`u toute partie compacte de
$F=\widehat F \setminus\{N,S\}$ ne contient qu'un nombre fini de points
de $L$ (ceci vient de l'errance des points,
propri\'et\'e 3 de l'hypoth\`ese \textrm{(H)}).
Soit $d$ une distance sur $\widehat F$ (propri\'et\'e 1 de
l'hypoth\`ese \textrm{(H)}). Puisque $L$ est localement fini, il
existe un r\'eel positif $r$ tel qu'aucun point de $L$ n'est \`a
distance  $r$ de $N$. On appelle alors $V_N$ l'ensemble des points
$x$ de $\widehat F$ v\'erifiant $d(N,x)<r$, $V_S$ l'ensemble de ceux
v\'erifiant $d(N,x) > r$.  On pose maintenant
$\esfn(L)=L \cap V_S \cap h^{-1}(V_N)$ et $\enfs(L)=L \cap V_N \cap
h^{-1}(V_S)$.
Ce sont des ensembles finis.
Le fait que $L$ soit limite de morceaux d'orbites allant d'un point
proche de $S$ \`a un point proche de $N$ permet de montrer
 que la diff\'erence entre le nombre
d'\'el\'ements de $\esfn(L)$ et le nombre
d'\'el\'ements	de $\enfs(L)$ vaut $1$.
D'autre part, un petit	argument combinatoire montre que cette
 diff\'erence est \'egale \`a la diff\'erence entre le nombre
 d'orbites de $L$ dans $\wsn(F)$ et  le nombre
 d'orbites de $L$ dans $\wns(F)$.
En particulier, $L$ doit contenir au moins une orbite dans $\wsn(F)$,
ce que
l'on voulait.
\end{proof}

\begin{proof}[\proofname\ du Lemme~\ref{lemm.monodynamie}]
On raisonne par l'absurde, en supposant qu'il existe une composante
compactement connexe de $F$ contenant deux points $x$ et $y$ ayant
des destins
dynamiques diff\'erents; pour fixer les id\'ees, on suppose par exemple
que $x \in \wfn(F)$ et $y \in \wfs(F)$. Par d\'efinition des composantes
compactement connexes, il existe un compact connexe $C$, inclus dans
$F$,  contenant $x$ et $y$.
Soit $C_\infty$ une valeur d'adh\'erence, dans l'espace ${\cal
K}(\widehat F)$,
 de la suite  $(h^n(C)_{n \geq 0})$. Il s'agit d'un sous-ensemble
 compact et
connexe de $\widehat F$, et puisque $C$ contient le point $x$ qui va
vers le Nord et le point $y$ qui va vers le Sud, l'ensemble $C_\infty$
contient $N$ et $S$.
D'apr\`es l'hypoth\`ese de minimalit\'e
(propri\'et\'e~\ref{item.minimalite} de l'hypoth\`ese \textrm{(H)}),
on a alors $C_\infty=\widehat F$. En
particulier $C_\infty$ contient $C$, et \emph{a fortiori} la limite
sup\'erieure de la suite $(h^n(C))_{n \geq 0}$ contient $C$. Ceci
contredit le premier ``souvenir du plan''
(hypoth\`ese~\ref{item.souvenirs} (a)
de \textrm{(H)}).
\end{proof}

\begin{proof}[\proofname\ du Lemme~\ref{lemm.composante-NS}]
Tout d'abord,  il est clair que l'image par $h$ d'une composante
compactement connexe de $F$ est encore une composante compactement
connexe de $F$; par cons\'equent, chaque composante compactement connexe
est soit
invariante, soit libre par $h$.

\paragraph{Premier cas}
Soit $E$ une composante compactement connexe de $F$ incluse dans $\wsn$,
que l'on suppose invariante par $h$.
Comme $E$ contient une orbite allant de $S$ \`a $N$,
 l'adh\'erence de $E$ dans $\widehat F$, qui est compacte et connexe,
 contient $N$ et
 $S$. Par minimalit\'e (propri\'et\'e~\ref{item.minimalite}),
 on a donc $\adhe(E)=\widehat F$.

\paragraph{Deuxi\`eme cas}
Soit $E$ une composante compactement connexe de $F$ incluse dans $\wsn$,
que l'on suppose libre par $h$. D'apr\`es le Corollaire~\ref{coro.ccc}
de l'appendice, $E$ ``va \`a l'infini'', \ie
 l'adh\'erence de $E$ contient l'un au moins des deux points
$N$ et $S$. \emph{On suppose dans un premier temps que $\adhe(E)$
contient le point $S$}. Nous allons montrer que dans ce cas, $E
\subset \adhe(h(E))$.  Soit $E_\infty$ une valeur d'adh\'erence, dans
${\cal K}(\widehat F)$, de la
suite des it\'er\'es positifs de $\adhe(E)$:  l'ensemble $E_\infty$
est compact et
connexe. Comme le point $S$ est dans $\adhe(E)$, il
appartient aussi \`a $E_\infty$. Comme
$E$ est incluse dans $\wsn$, l'ensemble $E_\infty$ contient aussi le
point $N$. Par minimalit\'e (hypoth\`ese~\ref{item.minimalite}),
 on a alors $E_\infty=\widehat F$.
En particulier, $E_\infty$ contient $E$. On conclut que  $E
\subset \adhe(h(E))$ gr\^ace au deuxi\`eme
``souvenir du plan'' (hypoth\`ese~\ref{item.souvenirs} (b) de
\textrm{(H)}).

On a donc prouv\'e que si $S \in \adhe(E)$, alors $E
\subset \adhe(h(E))$. On en d\'eduit d'abord que l'union des it\'er\'es de
$E$ est connexe. Ensuite, la suite $(h^n(\adhe(E)))$ est
une suite croissante d'\'el\'ements de ${\cal K}(\widehat F)$, donc
convergente. Par minimalit\'e, sa limite est $\widehat F$
et par cons\'equent, l'union des it\'er\'es de $E$ est dense dans
$\widehat
F$. Il est clair qu'il s'agit d'un ensemble invariant par $h$.

\emph{Il reste \`a traiter le cas o\`u $\adhe(E)$
contient le point $N$}. En \'etudiant cette fois-ci les it\'er\'es
de $E$ par $h^{-1}$, on montre alors que
$E \subset \adhe(h^{-1}(E))$, et on conclut de mani\`ere analogue.
\end{proof}

\begin{proof}[\proofname\ du Corollaire~\ref{coro.NS}]
Si l'ensemble $\wsn(F)$ n'est pas vide, d'apr\`es le
Lemme~\ref{lemm.monodynamie}, il contient une composante compactement
connexe $E$ de $F$. Le Lemme~\ref{lemm.composante-NS} entra\^\i ne alors
la densit\'e de l'union $E'$  des it\'er\'es de $E$, donc celle de
$\wsn(F)$. Ce lemme montre aussi que $E'$ est connexe; $\wsn(F)$ est
connexe car il contient un ensemble connexe dense.
\end{proof}

\begin{proof}[\proofname\ de la Proposition~\ref{prop.continus-minimaux}]
\emph{Montrons tout d'abord que l'un au moins des deux ensembles
$\wsn(F)$ et
$\wns(F)$ n'est pas vide.}\footnote{Ce raisonnement fournit une
alternative \`a la d\'emonstration de l'article~\cite{lero9}, alternative
qui
permet de se passer du lemme de Birkhoff, et donc de l'hypoth\`ese de
connexit\'e locale.}
On raisonne par l'absurde: dans le cas contraire, les deux ensembles
$\wss(F) \cup \{ S\}$ et $\wnn(F) \cup \{ N\}$ forment une partition de
$\hat F$. Par connexit\'e, il
existe un point de l'un qui est dans l'adh\'erence de l'autre.
On en d\'eduit que pour tout couple $V_N,V_S$ de voisinages de $N$ et $S$
respectivement, il existe un it\'er\'e positif de $V_S$ qui rencontre
$V_N$ (en utilisant la continuit\'e des it\'er\'es de $h$).
Le Lemme~\ref{lemm.CRAS} affirme alors que l'ensemble $\wsn(F)$ n'est pas
vide, ce que l'on voulait d\'emontrer.

\emph{Montrons maintenant que l'un au moins des deux ensembles $\wsn(F)$
et
$\wns(F)$ est vide.}\footnote{Cette preuve simple, \emph{via}
l'argument de $G_\delta$, m'a \'et\'e sugg\'er\'ee par Sylvain Crovisier.}
Supposons que l'ensemble $\wsn(F)$ ne soit pas vide. Alors il est dense
dans $F$ (Corollaire~\ref{coro.NS}).
D'autre part, puisque $F=\wfn(F) \cup \wfs(F)$, on peut \'ecrire
l'ensemble  $\wfn(F)$
 de la mani\`ere suivante:
$$
\wfn(F)=\{x \in F : \forall V_N \mbox{ voisinage de } N,\ \exists n>0
\mbox{ tel que } h^n(x) \in V_N \}
$$
On en d\'eduit que l'ensemble $\wsn(F)$ est un $G_\delta$ de $F$
(Rappelons qu'un ensemble est un $G_\delta$ s'il peut s'\'ecrire comme
une intersection d\'enombrable d'ouverts).
De mani\`ere sym\'etrique, l'ensemble $\wns(F)$ est aussi un
$G_\delta$ de $F$
qui est soit vide, soit dense dans $\widehat F$.
Comme $\widehat F$ est un espace m\'etrique compact, il a  la
propri\'et\'e de
Baire: il ne peut pas contenir deux ensembles $G_\delta$ denses
disjoints. On en d\'eduit que l'un des deux ensembles  $\wsn(F)$ et
$\wns(F)$ est vide.

{On peut  maintenant montrer que l'ensemble $\wnn(F)$ est
ferm\'e}: dans le cas contraire, il existerait un point $x$ dans
l'adh\'erence de $\wnn(F)$ qui viendrait du Sud ou qui irait vers le
Sud; dans les deux cas, on montre que pour tous voisinages $V_N$ et
$V_S$ de $N$ et $S$ respectivement, il existerait un it\'er\'e positif
de $V_N$ rencontrant $V_S$ et un it\'er\'e positif de $V_S$
rencontrant $V_N$.
 On peut alors
appliquer deux fois le Lemme~\ref{lemm.CRAS}, et voir qu'aucun
des deux ensembles  $\wsn(F)$ et
$\wns(F)$ n'est vide. Ceci contredirait ce qui pr\'ec\`ede.
De m\^eme, {l'ensemble $\wss(F)$ est
ferm\'e}, ce qui prouve le troisi\`eme point de la proposition.

Puisque les quatre ensembles $\wnn(F)$, $\wss(F)$, $\wsn(F)$ et
$\wns(F)$ forment une partition de $F$, on en d\'eduit que celui des
ensembles $\wsn(F)$ et $\wns(F)$ qui n'est pas vide est
ouvert. Le Corollaire~\ref{coro.NS} montre qu'il est \'egalement
dense et connexe.
Ceci termine la preuve de la proposition.
\end{proof}

\appendix
\section*{Appendices}
\addcontentsline{toc}{section}{Appendices}
\small
\section{Topologie}
\label{appe.topo}
\subsection{Topologie g\'en\'erale}
\begin{prop}\label{prop.bumpwire}
Soit $X$ un espace topologique connexe, m\'etrisable,
qui est localement connexe ou compact.
Soit $F$ une partie ferm\'ee non vide de $X$, et $A$ une composante
connexe de $X
\setminus F$.
Alors $\partial A$ est non vide et inclus dans $\partial F$. En
particulier, si $F$ est connexe, $A \cup F$ est connexe.
\end{prop}
Le cas localement connexe est tr\`es facile, la preuve est laiss\'ee au
 lecteur (par ailleurs, nous ne
l'utilisons que dans le cas o\`u $X$ est le plan $\bbR^2$). Le cas
 compact est plus difficile; la preuve n\'ecessite la notion
 d'\emph{$\epsilon$--cha\^\i ne}. En fait, nous ne l'utilisons que dans
 la preuve du corollaire suivant, sous une version faible: il nous
 suffit de savoir que $\partial A$
 \emph{rencontre} $\partial F$. Sous cette forme faible, l'\'enonc\'e
 est classique
 (voir \cite{kura1}, chap. V, §42, III, th\'eor\`eme 2).

\begin{coro}\label{coro.ccc}
On se place dans un espace m\'etrique $X$ connexe et localement compact,
mais non
compact. Notons  $\hat X=X \cup \{\infty\}$  le compactifi\'e
d'Alexandroff de $X$, et supposons que $\hat X$ est m\'etrisable. Alors
pour toute composante compactement connexe $E$ de $X$, le point
$\infty$ est dans l'adh\'erence dans $\hat X$  de $E$.
\end{coro}
\begin{proof}[\proofname\ du Corollaire~\ref{coro.ccc}]
Comme $X$ est localement compact, l'espace $\hat X$ est compact. Soit
$x$ un point de $X$, et
$V$ un voisinage ferm\'e du point  $\infty$ ne contenant pas
$x$. Soit $C_x$ la composante connexe
de $x$ dans $\hat X \setminus V$. D'apr\`es la
Proposition~\ref{prop.bumpwire}, l'adh\'erence de $C_x$  rencontre
$V$. Or l'adh\'erence de $C_x$ est compacte et ne contient pas le
point $\infty$, elle est donc incluse dans la
composante compactement connexe de $x$ dans $X$. Celle-ci
rencontre donc tout voisinage du point $\infty$, donc son adh\'erence
contient le point $\infty$.
\end{proof}

\subsection{Topologie plane}
Voici une version du lemme d'Alexander (que l'on peut prouver avec un peu
d'hom\-o\-lo\-gie):
\begin{theo2}[Corollaire 2 du th\'eor\`eme V.9.2 de \cite{newm51}]
\label{theo.lemme-alex}
Soient $U_1$ et $U_2$ deux ouverts connexes du plan. Si $U_1 \cup U_2=
\bbR^2$, alors $U_1 \cap U_2$ est connexe.
\end{theo2}
Ensuite, un \'enonc\'e dans la m\^eme veine, plus fin:
\begin{theo2}[Th\'eor\`eme V.14.2 de \cite{newm51}]
\label{theo.connexite}
Si $D$ est un ouvert connexe du plan, alors les composantes connexes
de l'ouvert $\bbR^2 \setminus \adhe(D)$ ont toutes des fronti\`eres
connexes.
\end{theo2}
On en d\'eduit imm\'ediatement:
\begin{theo2}[Th\'eor\`eme V.14.3 de \cite{newm51}]
\label{theo.separation-plane}
Si $x$ et $y$ sont deux points du plan s\'e\-par\-\'es par une partie ferm\'ee
$F$, alors il existe une composante connexe de $F$ qui les s\'epare.
\end{theo2}

\section{Exemples}
\label{sec.exemples}

Dans cette section, nous d\'ecrivons les composantes de Reeb sur
quelques exemples (tr\`es classiques, voir en particulier
l'appendice de~\cite{homm1}). Rappelons que les exemples les plus
simples sont
obtenus en int\'egrant les champs de vecteurs (\cf\ les
Figures~\fref{fig.Reeb} et \fref{fig.feuilletage} de l'introduction).
Pour ceux-ci, toute composante de Reeb $(F,G)$ associ\'ee \`a un
couple de points $(x,y)$ est minimale (sauf \'eventuellement si elle est
d\'eg\'en\'er\'ee),
 et les bords $F$ et $G$ sont des
droites topologiques invariantes par la dynamique. De plus, les
fl\`eches verticales  associ\'ees \`a une suite g\'eod\'esique
infinie sont toujours alternativement $\uparrow$ et $\downarrow$
(Section~\ref{sec.fleches-verticales}).

Les autres exemples sont obtenus en modifiant des flots. Sur
la \figref{fig.ex-compo}, \`a gauche, le couple $(x,y)$
poss\`ede une composante de Reeb minimale $(F,G)$ dont le bord
n\'egatif $F$ est libre (et non pas invariant, \cf
Lemme~\ref{lemm.invariance-crg}). Sur l'exemple de
droite, le couple $(x,y)$ poss\`ede une infinit\'e de composantes de
Reeb (dont le bord n\'egatif contient
 la  droite verticale $F$ et les it\'er\'es d'un
segment horizontal de longueur variable).
 Par contre, il y a bien une unique composante minimale, dont
les bords sont encore des droites topologiques. Notons qu'on peut
bien s\^ur
compliquer les bords  des composantes de Reeb minimales pour obtenir
des bords  qui ne sont pas localement connexes. On peut trouver des
constructions formelles de ces exemples dans \cite{lero1,lero8}.
\begin{figure}[htp]\anchor{fig.ex-compo}
\centerline{\hbox{\input{fig-ex-compo.pstex_t}}}
\caption{\label{fig.ex-compo}Exemples de composantes de Reeb (I)}
\end{figure}
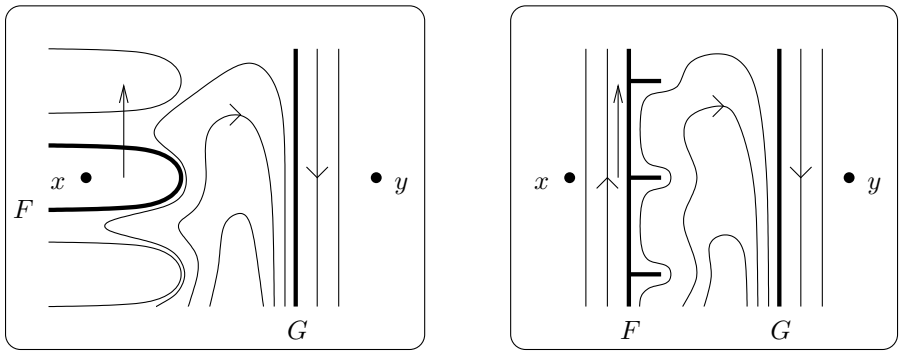
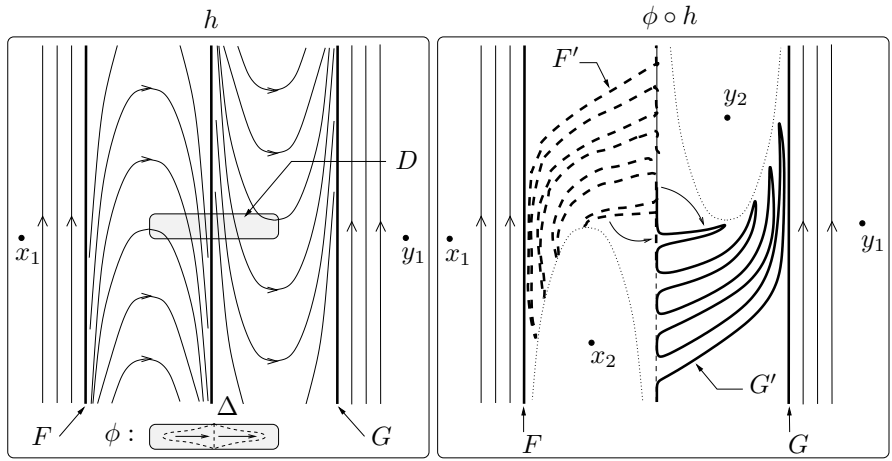
\begin{figure}[htp]\anchor{fig.flot-perturbe}
\centerline{\hbox{\input{fig-flot-perturbe2.pstex_t}}}
\caption{\label{fig.flot-perturbe}Exemples de composantes de Reeb (II)}
\end{figure}
Le dernier exemple est un peu plus compliqu\'e. Nous partons d'un
hom\'eo\-mor\-phisme $h$ temps $1$ d'un flot
(\figref{fig.flot-perturbe}, gauche). Les points $x_1$ et $y_1$
sont \`a $h$--distance $3$, deux
composantes de Reeb $(F, \Delta)$ et $(\Delta,G)$ leur sont
associ\'ees.
 On perturbe $h$ en le composant par un
hom\'eomorphisme $\phi$, \`a support dans un disque topologique
libre $D$ (comme indiqué sur la figure).
 Consid\'erons l'hom\'eomorphisme perturb\'e  $h'=\phi \circ h$.
La droite $\Delta$ est maintenant ``d\'edoubl\'ee'' en deux courbes
$F'$ et $G'$, d'adh\'erences respectives $F' \cup F$ et $G' \cup G$. La
distance de
translation entre les
deux points $x_1$ et $y_1$ a chut\'e \`a $2$; l'unique
composante  de Reeb minimale du couple $(x_1,y_1)$ est $(F,G)$.
Les points $x_2$ et $y_1$ sont aussi \`a $h'$--distance $2$; par
contre, la composante minimale de $(x_2,y_1)$ est $(F \cup
F',G)$. Ainsi, $(F \cup F',G)$ est une composante de Reeb pour chacun
des couple $(x_1,y_1)$ et $(x_2,y_1)$, cependant elle est minimale pour le
second mais pas pour le premier. De m\^eme, $(F, G' \cup G)$ est la
composante minimale associ\'ee au couple $(x_1,y_2)$. Les points
$x_2$ et $y_2$, eux, sont \`a $h'$--distance $1$.

Sur cet exemple, on peut  voir la n\'ecessit\'e de prendre l'ad\-h\'er\-ence
dans la formule d\'efinissant les composantes minimales
(d\'etaill\'ee \`a la Section~\ref{ssec.construction}; rappelons que
$F_1(x,y)$
est le bord de la composante connexe de $y$ dans l'ensemble
compl\'ementaire de
l'\emph{ad\-h\'er\-ence} de la boule $B_1(x)$, \textit{etc.}). En
effet, consid\'erons les boules pour la distance de translation
associ\'ee \`a	$h'$. La boule $B_1(x_1)$ est l'ouvert \`a gauche de
la courbe $G'$; sym\'etriquement, $B_1(y_1)$ est l'ouvert \`a
droite de la courbe $F'$. Ainsi, le compl\'ementaire de $B_1(x_1)$
est connexe, sa fronti\`ere est $G' \cup G$; de m\^eme, la
fronti\`ere du compl\'ementaire de $B_1(y_1)$ est $F' \cup F$; mais
le couple $(F' \cup F, G' \cup G)$ ne peut pas \^etre une composante
de Reeb
puisque $F'$ et $G'$ ne sont m\^eme pas disjoints.
Par contre le compl\'ementaire de l'\emph{adh\'erence} de $B_1(x_1)$
a deux composantes connexes (l'une contenant $y_1$, l'autre contenant
$y_2$), dont les fronti\`eres sont respectivement $G$ et $G' \cup G$.

On peut facilement \'etendre cet exemple pour qu'il existe une courbe
g\'eod\'esique infinie dont $x_1$ et $y_1$ sont deux sommets
(Section~\ref{sec.geodesiques-bout}). Dans
ce cas, le couple $(F,G)$ est une composante de Reeb associ\'ee \`a
cette courbe g\'eod\'esique, et les deux fl\`eches verticales indiquant
le sens de la dynamique dans $F$ et $G$ sont
\'egales \`a $\uparrow$ (rappelons que ceci n'arrive jamais dans le cas
du temps
$1$ d'un flot). On peut voir ici l'importance de la minimalit\'e des
composantes de Reeb consid\'er\'ees dans la d\'efinition des fl\`eches
verticales (Section~\ref{sec.fleches-verticales}): en effet, le
couple $(F,G' \cup G)$ est une composante de
Reeb (non minimale), et son bord $G' \cup G$ contient \`a la fois des
points allant de $N$ vers $S$ (dans $G'$) et d'autres allant de $S$
vers $N$ (dans $G$); ainsi, il n'est pas possible de d\'efinir le type
dynamique de ce bord (la Proposition~\ref{prop.continus-minimaux}
devient fausse pour les composantes non minimales).

De m\^eme, on peut \'etendre l'exemple pour qu'il
existe une courbe g\'eod\'esique infinie dont $x_1$ et $y_2$ sont des
sommets. Cette fois-ci, c'est $(F,G' \cup G)$ qui est une composante
de Reeb associ\'ee; les fl\`eches correspondantes sont $(\uparrow,
\downarrow)$. Confrontons cet exemple \`a l'\'etude de la dynamique
dans les continus minimaux (Proposition~\ref{prop.continus-minimaux}
et Section~\ref{sec.continus-minimaux}):
l'ensemble $G' \cup G \cup \{N,S\}$ est muni de la topologie de la
compactification Nord-Sud (pour laquelle l'adh\'erence de $G$ est $G
\cup \{S\}$, contrairement \`a ce que sugg\`ere la figure, \cf
Section~\ref{ssec.NS}). Cette
topologie en fait un espace compact, et les composantes compactement
connexes de $G \cup G'$ sont pr\'ecis\'ement $G$ et $G'$. Les points
de $G'$ vont tous du Nord au Sud, les points de $G$ vont tous du Sud
au Sud. On obtient ainsi un exemple avec deux types de comportements
dynamiques. On pourrait \'egalement obtenir un bord avec trois types
de comportements (ce qui est le maximum autoris\'e par la
Proposition~\ref{prop.continus-minimaux}; ceci est une remarque
commune avec C. Bonatti et
S. Crovisier). Terminons par deux questions.
\begin{ques*}
On consid\`ere la compactification Nord-Sud induite par une g\'eod\'esique
infinie, et $(F,G)$ une composante de Reeb associ\'ee \`a cette
g\'eod\'esique. L'ensemble $\wns(G)$ est-il toujours compactement connexe?
\end{ques*}
Dans l'exemple pr\'ec\'edent, l'ensemble $G' \cup G$ n'est pas
compactement connexe, mais $\wns(G' \cup G)=G'$ l'est.
\begin{ques*}
Soit $(F,G)$ une composante de Reeb minimale (associ\'ee \`a un couple
$(x,y)$ quelconque). Peut-il exister deux points $z$, $z'$ dans le
bord $G$ tels que $d_h(z,z')=3$?
\end{ques*}
Autrement dit, le ``$h$--diam\`etre'' d'un bord de  composante de Reeb
minimale peut-il \^etre \'egal \`a $3$? On voit facilement qu'il est
inf\'erieur ou \'egal \`a $3$, et dans l'exemple pr\'ec\'edent le $h$-diam\`etre
de $G'\cup G$ est $2$.

\end{document}

%% file: gtoutput.tex
\def\ifplaintex{\expandafter\ifx\csname documentclass\endcsname\relax}

\ifplaintex 
\else
\fi

\expandafter\ifx\csname beginpicture\endcsname\relax
\expandafter\ifx\csname documentclass\endcsname\relax
\input pictex \else
\input prepictex \input pictex \input postpictex \fi\fi

\def\gt{{\mathsurround=0pt\it $\cal G\mskip-2mu$eometry \&\ 
$\cal T\!\!$opology}}        

\def\gtp{{\mathsurround=0pt\it $\cal G\mskip-2mu$eometry \&\ 
$\cal T\!\!$opology $\cal P\!$ublications}}  

\def\lognumber#1{\def\thelognumber{#1}}
\def\volumenumber#1{\def\thevolumenumber{#1}}
\def\papernumber#1{\def\thepapernumber{#1}}
\def\volumeyear#1{\def\thevolumeyear{#1}}

\def\pagenumbers#1#2{\def\startpage{#1}\def\finishpage{#2}}
\def\published#1{\def\publishdate{#1}}
\def\proposed#1{\def\theproposer{#1}}
\def\seconded#1{\def\theseconders{#1}}
\def\received#1{\def\receiveddate{#1}}
\def\revised#1{\def\reviseddate{#1}}
\def\accepted#1{\def\accepteddate{#1}}
\def\asciititle#1{\def\theasciititle{#1}}
\def\covertitle#1{\def\thecovertitle{#1}}
\def\coverauthors#1{\def\thecoverauthors{#1}}
\def\asciiauthors#1{\def\theasciiauthors{#1}}
\def\asciiaddress#1{\def\theasciiaddress{#1}}

\long\def\asciiabstract#1{\long\def\theasciiabstract{#1}}

\let\\\par\let\thelognumber\relax
\let\thevolumenumber\relax\let\thepapernumber\relax
\let\thevolumeyear\relax\let\thesamplenumber\relax\let\startpage\relax
\let\finishpage\relax\let\publishdate\relax\let\receiveddate\relax
\let\reviseddate\relax\let\accepteddate\relax\let\theasciititle\relax
\let\thecovertitle\relax\let\theasciiauthors\relax\let\theasciiaddress\relax
\let\theasciiabstract\relax
\let\theasciiemail\relax\let\theshortauthors\relax\let\theshorttitle\relax
\let\thecoverauthors\relax

\long\def\maketitlep{   

\count0=\startpage

\gt\hfill      
\beginpicture
\setcoordinatesystem units <0.33truein, 0.33truein> point at 2.2 0.9
\setplotsymbol ({$\cal G$})
\plotsymbolspacing=9truept
\circulararc 315 degrees from 0 1 center at 0 0
\setplotsymbol ({$\cal T$})
\circulararc 315 degrees from 1 -1 center at 1 0
\endpicture
\break
{\small\ifx\thesamplenumber\relax 
Volume \else Sample
\fi\thevolumenumber\ (\thevolumeyear)
\startpage--\finishpage\nl
Published: \publishdate}
\vglue 0.5truein plus 0.4fil minus 0.1truein

{\parskip=0pt\leftskip 0pt plus 1fil\def\\{\par\smallskip}{\ifplaintex\large
\else\Large\fi\bf\thetitle}\par\medskip}   

\vglue 0pt plus 0.1fil 

{\parskip=0pt\leftskip 0pt plus 1fil\def\\{\par}{\sc\theauthors}
\par\medskip}

\vglue 0pt plus 0.1fil 

{\small\parskip=0pt\let\newline\\
{\leftskip 0pt plus 1fil\def\\{\par}{\sl\theaddress}\par}
\expandafter\ifx\theemail\relax    
\relax\else\vglue 5pt plus 0.02fil minus 2pt\def\\{\stdspace{\rm 
and}\stdspace} 
\cl{Email:\stdspace\tt\theemail}\fi
\ifx\theurl\relax                  
\relax\else\vglue 5pt plus 0.02fil minus 2pt\def\\{\stdspace{\rm 
and}\stdspace}
\cl{URL:\stdspace\tt\theurl}\fi\par}

\vglue 7pt plus 0.3fil minus 3pt

{\bf Abstract}
\vglue 5pt plus 0.1fil minus 2pt

\theabstract

\vglue 7pt plus 0.3fil minus 3pt

{\bf AMS Classification numbers}\quad Primary:\quad \theprimaryclass

Secondary:\quad \thesecondaryclass

\vglue 5pt plus 0.3fil minus 2pt

{\bf Keywords:}\quad \thekeywords

\vglue 10pt plus 0.5fil minus 5pt

{\small  Proposed: \theproposer\hfill Received: \receiveddate\nl
Seconded: \theseconders\hfill 
\ifx\reviseddate\relax                         
Accepted: \accepteddate                        
\else
Revised: \reviseddate                          
\fi}
\eject
}       

\let\maketitlepage\maketitlep
\let\maketitle\maketitlepage

\font\phead=cmsl9 scaled 950
\font=cmsl9 scaled 1050
\font\pnum=cmbx10 scaled 913
\font\lnum=cmbx10 
\font\pfoot=cmsl9 scaled 950
\font=cmsl9 scaled 1050
\ifplaintex
\headline{\vbox to 0pt{\vskip -4.5mm\line{\small\phead\ifnum
\count0=\startpage ISSN 1364-0380 (on line)
1465-3060 (printed) \hfill {\pnum\folio}\else\ifodd\count0\def\\{ }%
\ifx\theshorttitle\relax\thetitle\else\theshorttitle\fi\hfill{\pnum\folio}
\else\def\\{ and }{\pnum\folio}\hfill\ifx\theshortauthors\relax\theauthors
\else\theshortauthors\fi\fi\fi}\vss}}
\footline{\vbox to 0pt{\vglue 0mm\line{\small\pfoot\ifnum\count0=\startpage
\copyright\ \gtp\hfill\else
\gt, Volume \thevolumenumber\ (\thevolumeyear)\hfill\fi}\vss
}}
\else
\makeatletter
\def\@oddhead{{\small\lhead\ifnum\count0=\startpage ISSN 1364-0380 (on line)
1465-3060 (printed) \hfill {\lnum\number\count0}\else\ifodd\count0
\def\\{ }\ifx\theshorttitle\relax \thetitle \else\theshorttitle\fi\hfill
{\lnum\number\count0}\else\def\\{ and }{\lnum\number\count0}
\hfill\ifx\theshortauthors\relax 
\theauthors\else\theshortauthors\fi\fi\fi}}\def\@evenhead{\@oddhead}
\def\@oddfoot{\small\lfoot\ifnum\count0=\startpage\copyright\ \gtp\hfill\else
\gt, Volume \thevolumenumber\ (\thevolumeyear)\hfill\fi}
\def\@evenfoot{\@oddfoot}
\makeatother
\fi

\newwrite\gtoutfile
\long\gdef\makeheadfile{  
{\def\\{, }\def\s{ }
\immediate\openout\gtoutfile head.xxx
\immediate\write\gtoutfile{Proxy-for: \ifx\theasciiauthors\relax
\theauthors\else\theasciiauthors\fi\s<\ifx\theasciiemail\relax\theemail\else\theasciiemail\fi>}
\immediate\write\gtoutfile{\noexpand\\}
\immediate\write\gtoutfile{Authors: \ifx\theasciiauthors\relax
\theauthors\else\theasciiauthors\fi}
{\def\\{ }\immediate\write\gtoutfile{Title: \ifx\theasciititle\relax
\thetitle\else\theasciititle\fi}}
\immediate\write\gtoutfile{Subj-class: GT or SG or MG etc}
\immediate\write\gtoutfile{MSC-class: \theprimaryclass\ifx\thesecondaryclass\relax\else, \thesecondaryclass\fi}
\immediate\write\gtoutfile{Journal-ref: Geom. Topol. \thevolumenumber
(\thevolumeyear) \startpage-\finishpage}
\immediate\write\gtoutfile{Comments: Published by Geometry and Topology at}
\immediate\write\gtoutfile{\s\s http://www.maths.warwick.ac.uk/gt/GTVol\thevolumenumber/paper\thepapernumber.abs.html}
\immediate\write\gtoutfile{\noexpand\\}
\immediate\write\gtoutfile{}
\ifx\theasciiabstract\relax
\immediate\write\gtoutfile{\theabstract}\else
\immediate\write\gtoutfile{\theasciiabstract}\fi
\immediate\write\gtoutfile{}
\immediate\write\gtoutfile{\noexpand\\}
\immediate\write\gtoutfile{}
\immediate\closeout\gtoutfile}}  

\def\maketitlepage{\maketitlep\makeheadfile}
\let\maketitle\maketitlepage

%% file: fig-Reeb2.pstex_t
\begin{picture}(0,0)%
\includegraphics{fig-Reeb2.pstex}%
\end{picture}%
\setlength{\unitlength}{1579sp}%
\begingroup\makeatletter\ifx\SetFigFont\undefined%
\gdef\SetFigFont#1#2#3#4#5{%
  \reset@font\fontsize{#1}{#2pt}%
  \fontfamily{#3}\fontseries{#4}\fontshape{#5}%
  \selectfont}%
\fi\endgroup%
\begin{picture}(11982,5037)(3706,-6286)
\put(15673,-4016){\makebox(0,0)[lb]{\smash{{\SetFigFont{10}{12.0}{\rmdefault}{\mddefault}{\updefault}{\color[rgb]{0,0,0}$B$}%
}}}}
\end{picture}%

%% file: fig-feuilletage4.pstex_t
\begin{picture}(0,0)%
\includegraphics{fig-feuilletage4.pstex}%
\end{picture}%
\setlength{\unitlength}{1579sp}%
\begingroup\makeatletter\ifx\SetFigFont\undefined%
\gdef\SetFigFont#1#2#3#4#5{%
  \reset@font\fontsize{#1}{#2pt}%
  \fontfamily{#3}\fontseries{#4}\fontshape{#5}%
  \selectfont}%
\fi\endgroup%
\begin{picture}(10383,7317)(1165,-7635)
\put(6377,-1091){\makebox(0,0)[lb]{\smash{{\SetFigFont{14}{16.8}{\familydefault}{\mddefault}{\updefault}{\color[rgb]{0,0,0}\small$y$}%
}}}}
\put(2805,-7324){\makebox(0,0)[lb]{\smash{{\SetFigFont{14}{16.8}{\familydefault}{\mddefault}{\updefault}{\color[rgb]{0,0,0}\small$x$}%
}}}}
\end{picture}%

%% file: fig-feuilletage-geo2.pstex_t
\begin{picture}(0,0)%
\includegraphics{fig-feuilletage-geo2.pstex}%
\end{picture}%
\setlength{\unitlength}{1579sp}%
\begingroup\makeatletter\ifx\SetFigFont\undefined%
\gdef\SetFigFont#1#2#3#4#5{%
  \reset@font\fontsize{#1}{#2pt}%
  \fontfamily{#3}\fontseries{#4}\fontshape{#5}%
  \selectfont}%
\fi\endgroup%
\begin{picture}(10383,7317)(1165,-7635)
\put(2797,-5469){\makebox(0,0)[rb]{\smash{{\SetFigFont{10}{12.0}{\rmdefault}{\mddefault}{\updefault}{\color[rgb]{0,0,0}$h(\gamma)$}%
}}}}
\put(5002,-4119){\makebox(0,0)[rb]{\smash{{\SetFigFont{10}{12.0}{\rmdefault}{\mddefault}{\updefault}{\color[rgb]{0,0,0}$\gamma_2$}%
}}}}
\put(2805,-7324){\makebox(0,0)[lb]{\smash{{\SetFigFont{10}{12.0}{\rmdefault}{\mddefault}{\updefault}{\color[rgb]{0,0,0}$x$}%
}}}}
\put(3847,-7164){\makebox(0,0)[lb]{\smash{{\SetFigFont{10}{12.0}{\rmdefault}{\mddefault}{\updefault}{\color[rgb]{0,0,0}$\gamma_1$}%
}}}}
\put(5091,-5961){\makebox(0,0)[lb]{\smash{{\SetFigFont{10}{12.0}{\rmdefault}{\mddefault}{\updefault}{\color[rgb]{0,0,0}$x_1$}%
}}}}
\put(5811,-2661){\makebox(0,0)[lb]{\smash{{\SetFigFont{10}{12.0}{\rmdefault}{\mddefault}{\updefault}{\color[rgb]{0,0,0}$x_2$}%
}}}}
\put(6157,-2049){\makebox(0,0)[lb]{\smash{{\SetFigFont{10}{12.0}{\rmdefault}{\mddefault}{\updefault}{\color[rgb]{0,0,0}$\gamma_3$}%
}}}}
\put(6377,-1091){\makebox(0,0)[lb]{\smash{{\SetFigFont{10}{12.0}{\rmdefault}{\mddefault}{\updefault}{\color[rgb]{0,0,0}$y$}%
}}}}
\end{picture}%

%% file: fig-releve6.pstex_t
\begin{picture}(0,0)%
\includegraphics{fig-releve6.pstex}%
\end{picture}%
\setlength{\unitlength}{2842sp}%
\begingroup\makeatletter\ifx\SetFigFont\undefined%
\gdef\SetFigFont#1#2#3#4#5{%
  \reset@font\fontsize{#1}{#2pt}%
  \fontfamily{#3}\fontseries{#4}\fontshape{#5}%
  \selectfont}%
\fi\endgroup%
\begin{picture}(8642,10326)(804,-10133)
\put(2701,-2536){\makebox(0,0)[b]{\smash{{\SetFigFont{9}{10.8}{\rmdefault}{\mddefault}{\updefault}{\color[rgb]{0,0,0}\small $\la$}%
}}}}
\put(3601,-2536){\makebox(0,0)[b]{\smash{{\SetFigFont{9}{10.8}{\rmdefault}{\mddefault}{\updefault}{\color[rgb]{0,0,0}\small $\la$}%
}}}}
\put(5401,-2536){\makebox(0,0)[b]{\smash{{\SetFigFont{9}{10.8}{\rmdefault}{\mddefault}{\updefault}{\color[rgb]{0,0,0}\small $\la$}%
}}}}
\put(6301,-2536){\makebox(0,0)[b]{\smash{{\SetFigFont{9}{10.8}{\rmdefault}{\mddefault}{\updefault}{\color[rgb]{0,0,0}\small $\la$}%
}}}}
\put(7201,-2536){\makebox(0,0)[b]{\smash{{\SetFigFont{9}{10.8}{\rmdefault}{\mddefault}{\updefault}{\color[rgb]{0,0,0}\small $\la$}%
}}}}
\put(4501,-2536){\makebox(0,0)[b]{\smash{{\SetFigFont{9}{10.8}{\rmdefault}{\mddefault}{\updefault}{\color[rgb]{0,0,0}\small $\ra$}%
}}}}
\put(8101,-2536){\makebox(0,0)[b]{\smash{{\SetFigFont{9}{10.8}{\rmdefault}{\mddefault}{\updefault}{\color[rgb]{0,0,0}\small $\ra$}%
}}}}
\put(3976,-2536){\makebox(0,0)[b]{\smash{{\SetFigFont{9}{10.8}{\rmdefault}{\mddefault}{\updefault}{\color[rgb]{0,0,0}\small $\ua$}%
}}}}
\put(5776,-2536){\makebox(0,0)[b]{\smash{{\SetFigFont{9}{10.8}{\rmdefault}{\mddefault}{\updefault}{\color[rgb]{0,0,0}\small $\ua$}%
}}}}
\put(7576,-2536){\makebox(0,0)[b]{\smash{{\SetFigFont{9}{10.8}{\rmdefault}{\mddefault}{\updefault}{\color[rgb]{0,0,0}\small $\ua$}%
}}}}
\put(3076,-2536){\makebox(0,0)[b]{\smash{{\SetFigFont{9}{10.8}{\rmdefault}{\mddefault}{\updefault}{\color[rgb]{0,0,0}\small $\da$}%
}}}}
\put(4876,-2536){\makebox(0,0)[b]{\smash{{\SetFigFont{9}{10.8}{\rmdefault}{\mddefault}{\updefault}{\color[rgb]{0,0,0}\small $\da$}%
}}}}
\put(6676,-2536){\makebox(0,0)[b]{\smash{{\SetFigFont{9}{10.8}{\rmdefault}{\mddefault}{\updefault}{\color[rgb]{0,0,0}\small $\da$}%
}}}}
\put(1801,-2536){\makebox(0,0)[b]{\smash{{\SetFigFont{9}{10.8}{\rmdefault}{\mddefault}{\updefault}{\color[rgb]{0,0,0}\small $\la$}%
}}}}
\put(2176,-2536){\makebox(0,0)[b]{\smash{{\SetFigFont{9}{10.8}{\rmdefault}{\mddefault}{\updefault}{\color[rgb]{0,0,0}\small $\ua$}%
}}}}
\put(1801,-2836){\makebox(0,0)[b]{\smash{{\SetFigFont{9}{10.8}{\rmdefault}{\mddefault}{\updefault}{\color[rgb]{0,0,0}\small $\t x_1$}%
}}}}
\put(2701,-2836){\makebox(0,0)[b]{\smash{{\SetFigFont{9}{10.8}{\rmdefault}{\mddefault}{\updefault}{\color[rgb]{0,0,0}\small $\t x_2$}%
}}}}
\put(8101,-2836){\makebox(0,0)[b]{\smash{{\SetFigFont{9}{10.8}{\rmdefault}{\mddefault}{\updefault}{\color[rgb]{0,0,0}\small $\t x_8$}%
}}}}
\put(7201,-2836){\makebox(0,0)[b]{\smash{{\SetFigFont{9}{10.8}{\rmdefault}{\mddefault}{\updefault}{\color[rgb]{0,0,0}\small $\t x_7$}%
}}}}
\put(6301,-2836){\makebox(0,0)[b]{\smash{{\SetFigFont{9}{10.8}{\rmdefault}{\mddefault}{\updefault}{\color[rgb]{0,0,0}\small $\t x_6$}%
}}}}
\put(5401,-2836){\makebox(0,0)[b]{\smash{{\SetFigFont{9}{10.8}{\rmdefault}{\mddefault}{\updefault}{\color[rgb]{0,0,0}\small $\t x_5$}%
}}}}
\put(4501,-2836){\makebox(0,0)[b]{\smash{{\SetFigFont{9}{10.8}{\rmdefault}{\mddefault}{\updefault}{\color[rgb]{0,0,0}\small $\t x_4$}%
}}}}
\put(3601,-2836){\makebox(0,0)[b]{\smash{{\SetFigFont{9}{10.8}{\rmdefault}{\mddefault}{\updefault}{\color[rgb]{0,0,0}\small $\t x_3$}%
}}}}
\put(1501,-3661){\makebox(0,0)[b]{\smash{{\SetFigFont{9}{10.8}{\rmdefault}{\mddefault}{\updefault}{\color[rgb]{0,0,0}\small $F_1$}%
}}}}
\put(2101,-3661){\makebox(0,0)[b]{\smash{{\SetFigFont{9}{10.8}{\rmdefault}{\mddefault}{\updefault}{\color[rgb]{0,0,0}\small $G_1$}%
}}}}
\put(4201,-3661){\makebox(0,0)[b]{\smash{{\SetFigFont{9}{10.8}{\rmdefault}{\mddefault}{\updefault}{\color[rgb]{0,0,0}\small $F_4$}%
}}}}
\put(4801,-3661){\makebox(0,0)[b]{\smash{{\SetFigFont{9}{10.8}{\rmdefault}{\mddefault}{\updefault}{\color[rgb]{0,0,0}\small $G_4$}%
}}}}
\put(6001,-3661){\makebox(0,0)[b]{\smash{{\SetFigFont{9}{10.8}{\rmdefault}{\mddefault}{\updefault}{\color[rgb]{0,0,0}\small $F_6$}%
}}}}
\put(6601,-3661){\makebox(0,0)[b]{\smash{{\SetFigFont{9}{10.8}{\rmdefault}{\mddefault}{\updefault}{\color[rgb]{0,0,0}\small $G_6$}%
}}}}
\put(6901,-3661){\makebox(0,0)[b]{\smash{{\SetFigFont{9}{10.8}{\rmdefault}{\mddefault}{\updefault}{\color[rgb]{0,0,0}\small $F_7$}%
}}}}
\put(7501,-3661){\makebox(0,0)[b]{\smash{{\SetFigFont{9}{10.8}{\rmdefault}{\mddefault}{\updefault}{\color[rgb]{0,0,0}\small $G_7$}%
}}}}
\put(2401,-3661){\makebox(0,0)[b]{\smash{{\SetFigFont{9}{10.8}{\rmdefault}{\mddefault}{\updefault}{\color[rgb]{0,0,0}\small $F_2$}%
}}}}
\put(3001,-3661){\makebox(0,0)[b]{\smash{{\SetFigFont{9}{10.8}{\rmdefault}{\mddefault}{\updefault}{\color[rgb]{0,0,0}\small $G_2$}%
}}}}
\put(3301,-3661){\makebox(0,0)[b]{\smash{{\SetFigFont{9}{10.8}{\rmdefault}{\mddefault}{\updefault}{\color[rgb]{0,0,0}\small $F_3$}%
}}}}
\put(3901,-3661){\makebox(0,0)[b]{\smash{{\SetFigFont{9}{10.8}{\rmdefault}{\mddefault}{\updefault}{\color[rgb]{0,0,0}\small $G_3$}%
}}}}
\put(5101,-3661){\makebox(0,0)[b]{\smash{{\SetFigFont{9}{10.8}{\rmdefault}{\mddefault}{\updefault}{\color[rgb]{0,0,0}\small $F_5$}%
}}}}
\put(5701,-3661){\makebox(0,0)[b]{\smash{{\SetFigFont{9}{10.8}{\rmdefault}{\mddefault}{\updefault}{\color[rgb]{0,0,0}\small $G_5$}%
}}}}
\put(7801,-3661){\makebox(0,0)[b]{\smash{{\SetFigFont{9}{10.8}{\rmdefault}{\mddefault}{\updefault}{\color[rgb]{0,0,0}\small $F_8$}%
}}}}
\put(8401,-3661){\makebox(0,0)[b]{\smash{{\SetFigFont{9}{10.8}{\rmdefault}{\mddefault}{\updefault}{\color[rgb]{0,0,0}\small $G_8$}%
}}}}
\put(6451,-5086){\makebox(0,0)[b]{\smash{{\SetFigFont{9}{10.8}{\rmdefault}{\mddefault}{\updefault}{\color[rgb]{0,0,0}\small automorphisme de rev\^ etement}%
}}}}
\put(6076,-9511){\makebox(0,0)[b]{\smash{{\SetFigFont{9}{10.8}{\rmdefault}{\mddefault}{\updefault}{\color[rgb]{0,0,0}\small $x_2$}%
}}}}
\put(4201,-9511){\makebox(0,0)[b]{\smash{{\SetFigFont{9}{10.8}{\rmdefault}{\mddefault}{\updefault}{\color[rgb]{0,0,0}\small $x_1$}%
}}}}
\put(6301,-7336){\makebox(0,0)[b]{\smash{{\SetFigFont{9}{10.8}{\rmdefault}{\mddefault}{\updefault}{\color[rgb]{0,0,0}\small $x_3$}%
}}}}
\put(4126,-7336){\makebox(0,0)[b]{\smash{{\SetFigFont{9}{10.8}{\rmdefault}{\mddefault}{\updefault}{\color[rgb]{0,0,0}\small $x_4$}%
}}}}
\put(9001,-2386){\makebox(0,0)[b]{\smash{{\SetFigFont{14}{16.8}{\rmdefault}{\mddefault}{\updefault}{\color[rgb]{0,0,0}\normalsize $\Gamma$}%
}}}}
\put(1426,-2536){\makebox(0,0)[b]{\smash{{\SetFigFont{9}{10.8}{\rmdefault}{\mddefault}{\updefault}{\color[rgb]{0,0,0}\small $\da$}%
}}}}
\put(4951,-4561){\makebox(0,0)[b]{\smash{{\SetFigFont{9}{10.8}{\rmdefault}{\mddefault}{\updefault}{\color[rgb]{0,0,0}\normalsize $S$}%
}}}}
\put(4951,-61){\makebox(0,0)[b]{\smash{{\SetFigFont{9}{10.8}{\rmdefault}{\mddefault}{\updefault}{\color[rgb]{0,0,0}\normalsize $N$}%
}}}}
\put(3001,-8836){\makebox(0,0)[rb]{\smash{{\SetFigFont{9}{10.8}{\rmdefault}{\mddefault}{\updefault}{\color[rgb]{0,0,0}\small $\gamma$}%
}}}}
\put(3074,-6019){\makebox(0,0)[lb]{\smash{{\SetFigFont{14}{16.8}{\rmdefault}{\mddefault}{\updefault}{\color[rgb]{0,0,0}\normalsize rev\^etement}%
}}}}
\put(2776,-9886){\makebox(0,0)[lb]{\smash{{\SetFigFont{14}{16.8}{\rmdefault}{\mddefault}{\updefault}{\color[rgb]{0,0,0}\normalsize $h$}%
}}}}
\put(976,-5161){\makebox(0,0)[lb]{\smash{{\SetFigFont{14}{16.8}{\rmdefault}{\mddefault}{\updefault}{\color[rgb]{0,0,0}\normalsize $\t h$}%
}}}}
\end{picture}%

%% file: fig-decomposition.pstex_t
\begin{picture}(0,0)%
\includegraphics{fig-decomposition.pstex}%
\end{picture}%
\setlength{\unitlength}{3158sp}%
\begingroup\makeatletter\ifx\SetFigFont\undefined%
\gdef\SetFigFont#1#2#3#4#5{%
  \reset@font\fontsize{#1}{#2pt}%
  \fontfamily{#3}\fontseries{#4}\fontshape{#5}%
  \selectfont}%
\fi\endgroup%
\begin{picture}(4377,889)(1336,-2018)
\put(1351,-1936){\makebox(0,0)[lb]{\smash{{\SetFigFont{12}{14.4}{\rmdefault}{\mddefault}{\updefault}{\color[rgb]{0,0,0}\small$x_0$}%
}}}}
\put(2251,-1936){\makebox(0,0)[lb]{\smash{{\SetFigFont{12}{14.4}{\rmdefault}{\mddefault}{\updefault}{\color[rgb]{0,0,0}\small$x_1$}%
}}}}
\put(3151,-1936){\makebox(0,0)[lb]{\smash{{\SetFigFont{12}{14.4}{\rmdefault}{\mddefault}{\updefault}{\color[rgb]{0,0,0}\small$x_2$}%
}}}}
\put(4651,-1936){\makebox(0,0)[lb]{\smash{{\SetFigFont{12}{14.4}{\rmdefault}{\mddefault}{\updefault}{\color[rgb]{0,0,0}\small$x_{k-1}$}%
}}}}
\put(5551,-1936){\makebox(0,0)[lb]{\smash{{\SetFigFont{12}{14.4}{\rmdefault}{\mddefault}{\updefault}{\color[rgb]{0,0,0}\small$x_k$}%
}}}}
\put(3826,-1936){\makebox(0,0)[lb]{\smash{{\SetFigFont{12}{14.4}{\rmdefault}{\mddefault}{\updefault}{\color[rgb]{0,0,0}\small$\cdots$}%
}}}}
\put(2701,-1336){\makebox(0,0)[lb]{\smash{{\SetFigFont{12}{14.4}{\rmdefault}{\mddefault}{\updefault}{\color[rgb]{0,0,0}\small$\gamma_2$}%
}}}}
\put(3826,-1336){\makebox(0,0)[lb]{\smash{{\SetFigFont{12}{14.4}{\rmdefault}{\mddefault}{\updefault}{\color[rgb]{0,0,0}\small$\cdots$}%
}}}}
\put(5101,-1336){\makebox(0,0)[lb]{\smash{{\SetFigFont{12}{14.4}{\rmdefault}{\mddefault}{\updefault}{\color[rgb]{0,0,0}\small$\gamma_k$}%
}}}}
\put(1801,-1336){\makebox(0,0)[lb]{\smash{{\SetFigFont{12}{14.4}{\rmdefault}{\mddefault}{\updefault}{\color[rgb]{0,0,0}\small$\gamma_1$}%
}}}}
\end{picture}%

%% file: fig-simplification.pstex_t
\begin{picture}(0,0)%
\includegraphics{fig-simplification.pstex}%
\end{picture}%
\setlength{\unitlength}{3552sp}%
\begingroup\makeatletter\ifx\SetFigFont\undefined%
\gdef\SetFigFont#1#2#3#4#5{%
  \reset@font\fontsize{#1}{#2pt}%
  \fontfamily{#3}\fontseries{#4}\fontshape{#5}%
  \selectfont}%
\fi\endgroup%
\begin{picture}(5517,1207)(1411,-2159)
\put(3226,-2086){\makebox(0,0)[lb]{\smash{{\SetFigFont{11}{13.2}{\rmdefault}{\mddefault}{\updefault}{\color[rgb]{0,0,0}\small$x_1$}%
}}}}
\put(5026,-2086){\makebox(0,0)[lb]{\smash{{\SetFigFont{11}{13.2}{\rmdefault}{\mddefault}{\updefault}{\color[rgb]{0,0,0}\small$x_2$}%
}}}}
\put(6826,-2086){\makebox(0,0)[lb]{\smash{{\SetFigFont{11}{13.2}{\rmdefault}{\mddefault}{\updefault}{\color[rgb]{0,0,0}\small$x_3=y$}%
}}}}
\put(1426,-2086){\makebox(0,0)[lb]{\smash{{\SetFigFont{11}{13.2}{\rmdefault}{\mddefault}{\updefault}{\color[rgb]{0,0,0}\small$x=x_0$}%
}}}}
\put(1426,-1111){\makebox(0,0)[lb]{\smash{{\SetFigFont{11}{13.2}{\rmdefault}{\mddefault}{\updefault}{\color[rgb]{0,0,0}\small$x'_0=x_0$}%
}}}}
\put(3001,-1111){\makebox(0,0)[lb]{\smash{{\SetFigFont{11}{13.2}{\rmdefault}{\mddefault}{\updefault}{\color[rgb]{0,0,0}\small$x'_1$}%
}}}}
\put(4501,-1111){\makebox(0,0)[lb]{\smash{{\SetFigFont{11}{13.2}{\rmdefault}{\mddefault}{\updefault}{\color[rgb]{0,0,0}\small$x'_2$}%
}}}}
\put(6826,-1111){\makebox(0,0)[lb]{\smash{{\SetFigFont{11}{13.2}{\rmdefault}{\mddefault}{\updefault}{\color[rgb]{0,0,0}\small$x'_3=x_3$}%
}}}}
\end{picture}%

%% file: fig-22choseslune.pstex_t
\begin{picture}(0,0)%
\includegraphics{fig-22choseslune.pstex}%
\end{picture}%
\setlength{\unitlength}{3552sp}%
\begingroup\makeatletter\ifx\SetFigFont\undefined%
\gdef\SetFigFont#1#2#3#4#5{%
  \reset@font\fontsize{#1}{#2pt}%
  \fontfamily{#3}\fontseries{#4}\fontshape{#5}%
  \selectfont}%
\fi\endgroup%
\begin{picture}(5655,981)(1261,-2230)
\put(4126,-2161){\makebox(0,0)[lb]{\smash{{\SetFigFont{11}{13.2}{\rmdefault}{\mddefault}{\updefault}{\color[rgb]{0,0,0}\small$x_i$}%
}}}}
\put(3226,-1786){\makebox(0,0)[lb]{\smash{{\SetFigFont{11}{13.2}{\rmdefault}{\mddefault}{\updefault}{\color[rgb]{0,0,0}\small$\gamma_i$}%
}}}}
\put(5026,-1786){\makebox(0,0)[lb]{\smash{{\SetFigFont{11}{13.2}{\rmdefault}{\mddefault}{\updefault}{\color[rgb]{0,0,0}\small$\gamma_{i+1}$}%
}}}}
\put(6901,-1636){\makebox(0,0)[lb]{\smash{{\SetFigFont{11}{13.2}{\rmdefault}{\mddefault}{\updefault}{\color[rgb]{0,0,0}\small$D_2$}%
}}}}
\put(1276,-1636){\makebox(0,0)[lb]{\smash{{\SetFigFont{11}{13.2}{\rmdefault}{\mddefault}{\updefault}{\color[rgb]{0,0,0}\small$D_1$}%
}}}}
\end{picture}%

%% file: fig-fleches-hn.pstex_t
\begin{picture}(0,0)%
\includegraphics{fig-fleches-hn.pstex}%
\end{picture}%
\setlength{\unitlength}{3552sp}%
\begingroup\makeatletter\ifx\SetFigFont\undefined%
\gdef\SetFigFont#1#2#3#4#5{%
  \reset@font\fontsize{#1}{#2pt}%
  \fontfamily{#3}\fontseries{#4}\fontshape{#5}%
  \selectfont}%
\fi\endgroup%
\begin{picture}(5730,1430)(1111,-2009)
\put(3226,-1936){\makebox(0,0)[lb]{\smash{{\SetFigFont{11}{13.2}{\rmdefault}{\mddefault}{\updefault}{\color[rgb]{0,0,0}\small$\gamma_i$}%
}}}}
\put(5026,-1936){\makebox(0,0)[lb]{\smash{{\SetFigFont{11}{13.2}{\rmdefault}{\mddefault}{\updefault}{\color[rgb]{0,0,0}\small$\gamma_{i+1}$}%
}}}}
\put(6826,-1636){\makebox(0,0)[lb]{\smash{{\SetFigFont{11}{13.2}{\rmdefault}{\mddefault}{\updefault}{\color[rgb]{0,0,0}\small$D_{i+1}$}%
}}}}
\put(1351,-1636){\makebox(0,0)[lb]{\smash{{\SetFigFont{11}{13.2}{\rmdefault}{\mddefault}{\updefault}{\color[rgb]{0,0,0}\small$D_i$}%
}}}}
\put(1126,-886){\makebox(0,0)[lb]{\smash{{\SetFigFont{11}{13.2}{\rmdefault}{\mddefault}{\updefault}{\color[rgb]{0,0,0}\small$h(D_i)$}%
}}}}
\end{picture}%

%% file: fig-compo-Reeb.pstex_t
\begin{picture}(0,0)%
\includegraphics{fig-compo-Reeb.pstex}%
\end{picture}%
\setlength{\unitlength}{3552sp}%
\begingroup\makeatletter\ifx\SetFigFont\undefined%
\gdef\SetFigFont#1#2#3#4#5{%
  \reset@font\fontsize{#1}{#2pt}%
  \fontfamily{#3}\fontseries{#4}\fontshape{#5}%
  \selectfont}%
\fi\endgroup%
\begin{picture}(3848,2489)(1261,-2980)
\put(5026,-1786){\makebox(0,0)[lb]{\smash{{\SetFigFont{11}{13.2}{\rmdefault}{\mddefault}{\updefault}{\color[rgb]{0,0,0}\small$y$}%
}}}}
\put(3001,-2911){\makebox(0,0)[lb]{\smash{{\SetFigFont{11}{13.2}{\rmdefault}{\mddefault}{\updefault}{\color[rgb]{0,0,0}\small$F$}%
}}}}
\put(3601,-2911){\makebox(0,0)[lb]{\smash{{\SetFigFont{11}{13.2}{\rmdefault}{\mddefault}{\updefault}{\color[rgb]{0,0,0}\small$G$}%
}}}}
\put(1726,-1786){\makebox(0,0)[lb]{\smash{{\SetFigFont{11}{13.2}{\rmdefault}{\mddefault}{\updefault}{\color[rgb]{0,0,0}\small$x$}%
}}}}
\put(4501,-2461){\makebox(0,0)[lb]{\smash{{\SetFigFont{11}{13.2}{\rmdefault}{\mddefault}{\updefault}{\color[rgb]{0,0,0}\small bord n\'egatif}%
}}}}
\put(1276,-2461){\makebox(0,0)[lb]{\smash{{\SetFigFont{11}{13.2}{\rmdefault}{\mddefault}{\updefault}{\color[rgb]{0,0,0}\small bord positif}%
}}}}
\end{picture}%

%% file: fig-theoreme-crg.pstex_t
\begin{picture}(0,0)%
\includegraphics{fig-theoreme-crg.pstex}%
\end{picture}%
\setlength{\unitlength}{3552sp}%
\begingroup\makeatletter\ifx\SetFigFont\undefined%
\gdef\SetFigFont#1#2#3#4#5{%
  \reset@font\fontsize{#1}{#2pt}%
  \fontfamily{#3}\fontseries{#4}\fontshape{#5}%
  \selectfont}%
\fi\endgroup%
\begin{picture}(3398,2686)(1711,-2900)
\put(5026,-1786){\makebox(0,0)[lb]{\smash{{\SetFigFont{11}{13.2}{\rmdefault}{\mddefault}{\updefault}{\color[rgb]{0,0,0}$y$}%
}}}}
\put(1726,-1786){\makebox(0,0)[lb]{\smash{{\SetFigFont{11}{13.2}{\rmdefault}{\mddefault}{\updefault}{\color[rgb]{0,0,0}$x$}%
}}}}
\put(2401,-361){\makebox(0,0)[lb]{\smash{{\SetFigFont{8}{9.6}{\rmdefault}{\mddefault}{\updefault}{\color[rgb]{0,0,0}$F_1(x,y)$}%
}}}}
\put(3751,-2836){\makebox(0,0)[lb]{\smash{{\SetFigFont{8}{9.6}{\rmdefault}{\mddefault}{\updefault}{\color[rgb]{0,0,0}$F_1(y,x)$}%
}}}}
\put(2851,-2836){\makebox(0,0)[lb]{\smash{{\SetFigFont{8}{9.6}{\rmdefault}{\mddefault}{\updefault}{\color[rgb]{0,0,0}$F_2(y,x)$}%
}}}}
\put(1951,-2836){\makebox(0,0)[lb]{\smash{{\SetFigFont{8}{9.6}{\rmdefault}{\mddefault}{\updefault}{\color[rgb]{0,0,0}$F_3(y,x)$}%
}}}}
\put(3376,-361){\makebox(0,0)[lb]{\smash{{\SetFigFont{8}{9.6}{\rmdefault}{\mddefault}{\updefault}{\color[rgb]{0,0,0}$F_2(x,y)$}%
}}}}
\put(4276,-361){\makebox(0,0)[lb]{\smash{{\SetFigFont{8}{9.6}{\rmdefault}{\mddefault}{\updefault}{\color[rgb]{0,0,0}$F_3(x,y)$}%
}}}}
\end{picture}%

%% file: fig-rencontre-geodesique-crg.pstex_t
\begin{picture}(0,0)%
\includegraphics{fig-rencontre-geodesique-crg.pstex}%
\end{picture}%
\setlength{\unitlength}{3552sp}%
\begingroup\makeatletter\ifx\SetFigFont\undefined%
\gdef\SetFigFont#1#2#3#4#5{%
  \reset@font\fontsize{#1}{#2pt}%
  \fontfamily{#3}\fontseries{#4}\fontshape{#5}%
  \selectfont}%
\fi\endgroup%
\begin{picture}(4527,2274)(1186,-2698)
\put(1201,-1786){\makebox(0,0)[lb]{\smash{{\SetFigFont{11}{13.2}{\rmdefault}{\mddefault}{\updefault}{\color[rgb]{0,0,0}\small $x$}%
}}}}
\put(5626,-1786){\makebox(0,0)[lb]{\smash{{\SetFigFont{11}{13.2}{\rmdefault}{\mddefault}{\updefault}{\color[rgb]{0,0,0}\small $y$}%
}}}}
\put(2101,-2611){\makebox(0,0)[lb]{\smash{{\SetFigFont{11}{13.2}{\rmdefault}{\mddefault}{\updefault}{\color[rgb]{0,0,0}\small $F_{d-i}(y,x)$}%
}}}}
\put(3826,-2611){\makebox(0,0)[lb]{\smash{{\SetFigFont{11}{13.2}{\rmdefault}{\mddefault}{\updefault}{\color[rgb]{0,0,0}\small $F_i(x,y)$}%
}}}}
\put(2401,-1411){\makebox(0,0)[lb]{\smash{{\SetFigFont{11}{13.2}{\rmdefault}{\mddefault}{\updefault}{\color[rgb]{0,0,0}\small $\gamma_i$}%
}}}}
\put(3301,-1786){\makebox(0,0)[lb]{\smash{{\SetFigFont{11}{13.2}{\rmdefault}{\mddefault}{\updefault}{\color[rgb]{0,0,0}\small $x_i$}%
}}}}
\put(4201,-1411){\makebox(0,0)[lb]{\smash{{\SetFigFont{11}{13.2}{\rmdefault}{\mddefault}{\updefault}{\color[rgb]{0,0,0}\small $\gamma_{i+1}$}%
}}}}
\put(1208,-983){\makebox(0,0)[lb]{\smash{{\SetFigFont{11}{13.2}{\rmdefault}{\mddefault}{\updefault}{\color[rgb]{0,0,0}\small $\gamma$}%
}}}}
\end{picture}%

%% file: fig-iteres-disque.pstex_t
\begin{picture}(0,0)%
\includegraphics{fig-iteres-disque.pstex}%
\end{picture}%
\setlength{\unitlength}{3552sp}%
\begingroup\makeatletter\ifx\SetFigFont\undefined%
\gdef\SetFigFont#1#2#3#4#5{%
  \reset@font\fontsize{#1}{#2pt}%
  \fontfamily{#3}\fontseries{#4}\fontshape{#5}%
  \selectfont}%
\fi\endgroup%
\begin{picture}(6016,3002)(1261,-3062)
\put(3826,-2011){\makebox(0,0)[lb]{\smash{{\SetFigFont{11}{13.2}{\rmdefault}{\mddefault}{\updefault}{\color[rgb]{0,0,0}\small $D$}%
}}}}
\put(2701,-2386){\makebox(0,0)[lb]{\smash{{\SetFigFont{11}{13.2}{\rmdefault}{\mddefault}{\updefault}{\color[rgb]{0,0,0}\small $\lim^-D$}%
}}}}
\put(3751,-2611){\makebox(0,0)[lb]{\smash{{\SetFigFont{11}{13.2}{\rmdefault}{\mddefault}{\updefault}{\color[rgb]{0,0,0}\small $U(D)$}%
}}}}
\put(6376,-2611){\makebox(0,0)[lb]{\smash{{\SetFigFont{11}{13.2}{\rmdefault}{\mddefault}{\updefault}{\color[rgb]{0,0,0}\small $U^+(D)$}%
}}}}
\put(1276,-2611){\makebox(0,0)[lb]{\smash{{\SetFigFont{11}{13.2}{\rmdefault}{\mddefault}{\updefault}{\color[rgb]{0,0,0}\small $U^-(D)$}%
}}}}
\put(4801,-2386){\makebox(0,0)[lb]{\smash{{\SetFigFont{11}{13.2}{\rmdefault}{\mddefault}{\updefault}{\color[rgb]{0,0,0}\small $\lim^+D$}%
}}}}
\end{picture}%

%% file: fig-type1.pstex_t
\begin{picture}(0,0)%
\includegraphics{fig-type1.pstex}%
\end{picture}%
\setlength{\unitlength}{2763sp}%
\begingroup\makeatletter\ifx\SetFigFont\undefined%
\gdef\SetFigFont#1#2#3#4#5{%
  \reset@font\fontsize{#1}{#2pt}%
  \fontfamily{#3}\fontseries{#4}\fontshape{#5}%
  \selectfont}%
\fi\endgroup%
\begin{picture}(6774,2583)(2689,-3145)
\put(3861,-1926){\makebox(0,0)[lb]{\smash{{\SetFigFont{8}{9.6}{\rmdefault}{\mddefault}{\updefault}{\color[rgb]{0,0,0}\small $x$}%
}}}}
\put(5796,-733){\makebox(0,0)[lb]{\smash{{\SetFigFont{8}{9.6}{\rmdefault}{\mddefault}{\updefault}{\color[rgb]{0,0,0}\hspace{-2pt}\small $\partial U(D)$}%
}}}}
\put(7611,-1926){\makebox(0,0)[lb]{\smash{{\SetFigFont{8}{9.6}{\rmdefault}{\mddefault}{\updefault}{\color[rgb]{0,0,0}\small $x$}%
}}}}
\put(5776,-1186){\makebox(0,0)[lb]{\smash{{\SetFigFont{8}{9.6}{\rmdefault}{\mddefault}{\updefault}{\color[rgb]{0,0,0}\hspace{-4pt}\small $h^n(D)$}%
}}}}
\put(5589,-1711){\makebox(0,0)[lb]{\smash{{\SetFigFont{8}{9.6}{\rmdefault}{\mddefault}{\updefault}{\color[rgb]{0,0,0}\small $\alpha$}%
}}}}
\put(3268,-3058){\makebox(0,0)[lb]{\smash{{\SetFigFont{10}{12.0}{\rmdefault}{\mddefault}{\updefault}{\color[rgb]{0,0,0}\small type $\partial U(D)$}%
}}}}
\put(7176,-3058){\makebox(0,0)[lb]{\smash{{\SetFigFont{10}{12.0}{\rmdefault}{\mddefault}{\updefault}{\color[rgb]{0,0,0}\small type $D$}%
}}}}
\end{picture}%

%% file: fig-type2.pstex_t
\begin{picture}(0,0)%
\includegraphics{fig-type2.pstex}%
\end{picture}%
\setlength{\unitlength}{2763sp}%
\begingroup\makeatletter\ifx\SetFigFont\undefined%
\gdef\SetFigFont#1#2#3#4#5{%
  \reset@font\fontsize{#1}{#2pt}%
  \fontfamily{#3}\fontseries{#4}\fontshape{#5}%
  \selectfont}%
\fi\endgroup%
\begin{picture}(3895,2049)(2426,-2611)
\put(3861,-1926){\makebox(0,0)[lb]{\smash{{\SetFigFont{8}{9.6}{\rmdefault}{\mddefault}{\updefault}{\color[rgb]{0,0,0}\small $x$}%
}}}}
\put(5091,-1861){\makebox(0,0)[lb]{\smash{{\SetFigFont{8}{9.6}{\rmdefault}{\mddefault}{\updefault}{\color[rgb]{0,0,0}\small $\alpha$}%
}}}}
\put(2441,-1856){\makebox(0,0)[lb]{\smash{{\SetFigFont{8}{9.6}{\rmdefault}{\mddefault}{\updefault}{\color[rgb]{0,0,0}\small $\alpha'$}%
}}}}
\put(6306,-1258){\makebox(0,0)[lb]{\smash{{\SetFigFont{8}{9.6}{\rmdefault}{\mddefault}{\updefault}{\color[rgb]{0,0,0}\small $h^n(D)$}%
}}}}
\put(5796,-733){\makebox(0,0)[lb]{\smash{{\SetFigFont{8}{9.6}{\rmdefault}{\mddefault}{\updefault}{\color[rgb]{0,0,0}\small $\partial U(D)$}%
}}}}
\end{picture}%

%% file: fig-idee-crg.pstex_t
\begin{picture}(0,0)%
\includegraphics{fig-idee-crg.pstex}%
\end{picture}%
\setlength{\unitlength}{3552sp}%
\begingroup\makeatletter\ifx\SetFigFont\undefined%
\gdef\SetFigFont#1#2#3#4#5{%
  \reset@font\fontsize{#1}{#2pt}%
  \fontfamily{#3}\fontseries{#4}\fontshape{#5}%
  \selectfont}%
\fi\endgroup%
\begin{picture}(5702,3002)(1800,-3062)
\put(2776,-1036){\makebox(0,0)[lb]{\smash{{\SetFigFont{11}{13.2}{\rmdefault}{\mddefault}{\updefault}{\color[rgb]{0,0,0}\small $\lim^-z$}%
}}}}
\put(3976,-1036){\makebox(0,0)[lb]{\smash{{\SetFigFont{11}{13.2}{\rmdefault}{\mddefault}{\updefault}{\color[rgb]{0,0,0}\small $\lim^+z$}%
}}}}
\put(6856,-2108){\makebox(0,0)[lb]{\smash{{\SetFigFont{11}{13.2}{\rmdefault}{\mddefault}{\updefault}{\color[rgb]{0,0,0}\small $x$}%
}}}}
\put(4184,-1461){\makebox(0,0)[lb]{\smash{{\SetFigFont{11}{13.2}{\rmdefault}{\mddefault}{\updefault}{\color[rgb]{0,0,0}\small $z$}%
}}}}
\put(4108,-1918){\makebox(0,0)[lb]{\smash{{\SetFigFont{11}{13.2}{\rmdefault}{\mddefault}{\updefault}{\color[rgb]{0,0,0}\small $z'$}%
}}}}
\put(3076,-2386){\makebox(0,0)[lb]{\smash{{\SetFigFont{11}{13.2}{\rmdefault}{\mddefault}{\updefault}{\color[rgb]{0,0,0}\small $B_i(x)$}%
}}}}
\put(2251,-2686){\makebox(0,0)[lb]{\smash{{\SetFigFont{11}{13.2}{\rmdefault}{\mddefault}{\updefault}{\color[rgb]{0,0,0}\small $O_x(z)=O_x(z')$}%
}}}}
\end{picture}%

%% file: fig-calculs-distance2.pstex_t
\begin{picture}(0,0)%
\includegraphics{fig-calculs-distance2.pstex}%
\end{picture}%
\setlength{\unitlength}{3552sp}%
\begingroup\makeatletter\ifx\SetFigFont\undefined%
\gdef\SetFigFont#1#2#3#4#5{%
  \reset@font\fontsize{#1}{#2pt}%
  \fontfamily{#3}\fontseries{#4}\fontshape{#5}%
  \selectfont}%
\fi\endgroup%
\begin{picture}(3425,3002)(1241,-3062)
\put(4126,-1336){\makebox(0,0)[lb]{\smash{{\SetFigFont{11}{13.2}{\rmdefault}{\mddefault}{\updefault}{\color[rgb]{0,0,0}\small $D$}%
}}}}
\put(2551,-586){\makebox(0,0)[lb]{\smash{{\SetFigFont{11}{13.2}{\rmdefault}{\mddefault}{\updefault}{\color[rgb]{0,0,0}\small $h^{n_0}(D)$}%
}}}}
\put(4206,-1836){\makebox(0,0)[lb]{\smash{{\SetFigFont{11}{13.2}{\rmdefault}{\mddefault}{\updefault}{\color[rgb]{0,0,0}\small $y$}%
}}}}
\put(4651,-2536){\makebox(0,0)[lb]{\smash{{\SetFigFont{11}{13.2}{\rmdefault}{\mddefault}{\updefault}{\color[rgb]{0,0,0}\small $x$}%
}}}}
\put(2436,-2621){\makebox(0,0)[lb]{\smash{{\SetFigFont{11}{13.2}{\rmdefault}{\mddefault}{\updefault}{\color[rgb]{0,0,0}\small $p$}%
}}}}
\put(1386,-1581){\makebox(0,0)[lb]{\smash{{\SetFigFont{11}{13.2}{\rmdefault}{\mddefault}{\updefault}{\color[rgb]{0,0,0}\small $\gamma$}%
}}}}
\put(1256,-2361){\makebox(0,0)[lb]{\smash{{\SetFigFont{11}{13.2}{\rmdefault}{\mddefault}{\updefault}{\color[rgb]{0,0,0}\small $\alpha$}%
}}}}
\put(2861,-2641){\makebox(0,0)[lb]{\smash{{\SetFigFont{11}{13.2}{\rmdefault}{\mddefault}{\updefault}{\color[rgb]{0,0,0}\small $s$}%
}}}}
\end{picture}%

%% file: fig-distances.pstex_t
\begin{picture}(0,0)%
\includegraphics{fig-distances.pstex}%
\end{picture}%
\setlength{\unitlength}{3552sp}%
\begingroup\makeatletter\ifx\SetFigFont\undefined%
\gdef\SetFigFont#1#2#3#4#5{%
  \reset@font\fontsize{#1}{#2pt}%
  \fontfamily{#3}\fontseries{#4}\fontshape{#5}%
  \selectfont}%
\fi\endgroup%
\begin{picture}(3555,2556)(1711,-2980)
\put(5026,-1786){\makebox(0,0)[lb]{\smash{{\SetFigFont{11}{13.2}{\rmdefault}{\mddefault}{\updefault}{\color[rgb]{0,0,0}\small $y$}%
}}}}
\put(3001,-2911){\makebox(0,0)[lb]{\smash{{\SetFigFont{11}{13.2}{\rmdefault}{\mddefault}{\updefault}{\color[rgb]{0,0,0}\small $F$}%
}}}}
\put(3601,-2911){\makebox(0,0)[lb]{\smash{{\SetFigFont{11}{13.2}{\rmdefault}{\mddefault}{\updefault}{\color[rgb]{0,0,0}\small $G$}%
}}}}
\put(1726,-1786){\makebox(0,0)[lb]{\smash{{\SetFigFont{11}{13.2}{\rmdefault}{\mddefault}{\updefault}{\color[rgb]{0,0,0}\small $x$}%
}}}}
\put(3076,-1988){\makebox(0,0)[lb]{\smash{{\SetFigFont{11}{13.2}{\rmdefault}{\mddefault}{\updefault}{\color[rgb]{0,0,0}\small $z'$}%
}}}}
\put(3878,-1913){\makebox(0,0)[lb]{\smash{{\SetFigFont{11}{13.2}{\rmdefault}{\mddefault}{\updefault}{\color[rgb]{0,0,0}\small $z$}%
}}}}
\put(5251,-953){\makebox(0,0)[lb]{\smash{{\SetFigFont{11}{13.2}{\rmdefault}{\mddefault}{\updefault}{\color[rgb]{0,0,0}\small $\gamma$}%
}}}}
\end{picture}%

%% file: fig-geodesique-infinie2.pstex_t
\begin{picture}(0,0)%
\includegraphics{fig-geodesique-infinie2.pstex}%
\end{picture}%
\setlength{\unitlength}{3552sp}%
\begingroup\makeatletter\ifx\SetFigFont\undefined%
\gdef\SetFigFont#1#2#3#4#5{%
  \reset@font\fontsize{#1}{#2pt}%
  \fontfamily{#3}\fontseries{#4}\fontshape{#5}%
  \selectfont}%
\fi\endgroup%
\begin{picture}(4599,3769)(1039,-3509)
\put(3076,-2911){\makebox(0,0)[lb]{\smash{{\SetFigFont{11}{13.2}{\rmdefault}{\mddefault}{\updefault}{\color[rgb]{0,0,0}\small $F_0$}%
}}}}
\put(2101,-2911){\makebox(0,0)[lb]{\smash{{\SetFigFont{11}{13.2}{\rmdefault}{\mddefault}{\updefault}{\color[rgb]{0,0,0}\small $F_{-1}$}%
}}}}
\put(3901,-2911){\makebox(0,0)[lb]{\smash{{\SetFigFont{11}{13.2}{\rmdefault}{\mddefault}{\updefault}{\color[rgb]{0,0,0}\small $F_1$}%
}}}}
\put(3601,-361){\makebox(0,0)[lb]{\smash{{\SetFigFont{11}{13.2}{\rmdefault}{\mddefault}{\updefault}{\color[rgb]{0,0,0}\small $G_0$}%
}}}}
\put(4426,-361){\makebox(0,0)[lb]{\smash{{\SetFigFont{11}{13.2}{\rmdefault}{\mddefault}{\updefault}{\color[rgb]{0,0,0}\small $G_1$}%
}}}}
\put(3301,-1861){\makebox(0,0)[lb]{\smash{{\SetFigFont{11}{13.2}{\rmdefault}{\mddefault}{\updefault}{\color[rgb]{0,0,0}\small $x_0$}%
}}}}
\put(4126,-1861){\makebox(0,0)[lb]{\smash{{\SetFigFont{11}{13.2}{\rmdefault}{\mddefault}{\updefault}{\color[rgb]{0,0,0}\small $x_1$}%
}}}}
\put(4876,-1861){\makebox(0,0)[lb]{\smash{{\SetFigFont{11}{13.2}{\rmdefault}{\mddefault}{\updefault}{\color[rgb]{0,0,0}\small $x_2$}%
}}}}
\put(2326,-1861){\makebox(0,0)[lb]{\smash{{\SetFigFont{11}{13.2}{\rmdefault}{\mddefault}{\updefault}{\color[rgb]{0,0,0}\small $x_{-1}$}%
}}}}
\put(1651,-1861){\makebox(0,0)[lb]{\smash{{\SetFigFont{11}{13.2}{\rmdefault}{\mddefault}{\updefault}{\color[rgb]{0,0,0}\small $x_{-2}$}%
}}}}
\put(3976,-3436){\makebox(0,0)[lb]{\smash{{\SetFigFont{11}{13.2}{\rmdefault}{\mddefault}{\updefault}{\color[rgb]{0,0,0}\small $O_+(F_0)$}%
}}}}
\put(2701,89){\makebox(0,0)[lb]{\smash{{\SetFigFont{11}{13.2}{\rmdefault}{\mddefault}{\updefault}{\color[rgb]{0,0,0}\small $O_-(G_0)$}%
}}}}
\put(4351,89){\makebox(0,0)[lb]{\smash{{\SetFigFont{11}{13.2}{\rmdefault}{\mddefault}{\updefault}{\color[rgb]{0,0,0}\small $O_+(G_0)$}%
}}}}
\put(2551,-361){\makebox(0,0)[lb]{\smash{{\SetFigFont{11}{13.2}{\rmdefault}{\mddefault}{\updefault}{\color[rgb]{0,0,0}\small $G_{-1}$}%
}}}}
\put(2101,-3436){\makebox(0,0)[lb]{\smash{{\SetFigFont{11}{13.2}{\rmdefault}{\mddefault}{\updefault}{\color[rgb]{0,0,0}\small $O_-(F_0)$}%
}}}}
\end{picture}%

%% file: fig-compo-bouts2.pstex_t
\begin{picture}(0,0)%
\includegraphics{fig-compo-bouts2.pstex}%
\end{picture}%
\setlength{\unitlength}{3158sp}%
\begingroup\makeatletter\ifx\SetFigFont\undefined%
\gdef\SetFigFont#1#2#3#4#5{%
  \reset@font\fontsize{#1}{#2pt}%
  \fontfamily{#3}\fontseries{#4}\fontshape{#5}%
  \selectfont}%
\fi\endgroup%
\begin{picture}(5055,2982)(293,-3509)
\put(901,-1756){\makebox(0,0)[lb]{\smash{{\SetFigFont{10}{12.0}{\rmdefault}{\mddefault}{\updefault}{\color[rgb]{0,0,0}\small $x_{-k}$}%
}}}}
\put(1396,-2453){\makebox(0,0)[lb]{\smash{{\SetFigFont{10}{12.0}{\rmdefault}{\mddefault}{\updefault}{\color[rgb]{0,0,0}\small ${x_{-k}}'$}%
}}}}
\put(4876,-1696){\makebox(0,0)[lb]{\smash{{\SetFigFont{10}{12.0}{\rmdefault}{\mddefault}{\updefault}{\color[rgb]{0,0,0}\small ${x_k}'$}%
}}}}
\put(2701,-2236){\makebox(0,0)[lb]{\smash{{\SetFigFont{10}{12.0}{\rmdefault}{\mddefault}{\updefault}{\color[rgb]{0,0,0}\small $F_0$}%
}}}}
\put(3826,-2236){\makebox(0,0)[lb]{\smash{{\SetFigFont{10}{12.0}{\rmdefault}{\mddefault}{\updefault}{\color[rgb]{0,0,0}\small $G_0$}%
}}}}
\put(1501,-886){\makebox(0,0)[lb]{\smash{{\SetFigFont{10}{12.0}{\rmdefault}{\mddefault}{\updefault}{\color[rgb]{0,0,0}\small $B_p(x_{-k})$}%
}}}}
\put(2026,-3436){\makebox(0,0)[lb]{\smash{{\SetFigFont{10}{12.0}{\rmdefault}{\mddefault}{\updefault}{\color[rgb]{0,0,0}\small $O_-(F_0)$}%
}}}}
\put(4313,-698){\makebox(0,0)[lb]{\smash{{\SetFigFont{10}{12.0}{\rmdefault}{\mddefault}{\updefault}{\color[rgb]{0,0,0}\small $O_+(G_0)$}%
}}}}
\end{picture}%

%% file: fig-compactification.pstex_t
\begin{picture}(0,0)%
\includegraphics{fig-compactification.pstex}%
\end{picture}%
\setlength{\unitlength}{3158sp}%
\begingroup\makeatletter\ifx\SetFigFont\undefined%
\gdef\SetFigFont#1#2#3#4#5{%
  \reset@font\fontsize{#1}{#2pt}%
  \fontfamily{#3}\fontseries{#4}\fontshape{#5}%
  \selectfont}%
\fi\endgroup%
\begin{picture}(3702,2703)(1789,-2830)
\put(5476,-1636){\makebox(0,0)[lb]{\smash{{\SetFigFont{10}{12.0}{\rmdefault}{\mddefault}{\updefault}{\color[rgb]{0,0,0}\small $\Gamma$}%
}}}}
\put(4801,-1261){\makebox(0,0)[lb]{\smash{{\SetFigFont{10}{12.0}{\rmdefault}{\mddefault}{\updefault}{\color[rgb]{0,0,0}\small $O_N(\Gamma)$}%
}}}}
\put(3676,-2761){\makebox(0,0)[lb]{\smash{{\SetFigFont{10}{12.0}{\rmdefault}{\mddefault}{\updefault}{\color[rgb]{0,0,0}\small $S$}%
}}}}
\put(3676,-286){\makebox(0,0)[lb]{\smash{{\SetFigFont{10}{12.0}{\rmdefault}{\mddefault}{\updefault}{\color[rgb]{0,0,0}\small $N$}%
}}}}
\put(4801,-2011){\makebox(0,0)[lb]{\smash{{\SetFigFont{10}{12.0}{\rmdefault}{\mddefault}{\updefault}{\color[rgb]{0,0,0}\small $O_S(\Gamma)$}%
}}}}
\put(2926,-2236){\makebox(0,0)[lb]{\smash{{\SetFigFont{10}{12.0}{\rmdefault}{\mddefault}{\updefault}{\color[rgb]{0,0,0}\small $F$}%
}}}}
\end{picture}%

%% file: fig-bandeB.pstex_t
\begin{picture}(0,0)%
\includegraphics{fig-bandeB.pstex}%
\end{picture}%
\setlength{\unitlength}{3158sp}%
\begingroup\makeatletter\ifx\SetFigFont\undefined%
\gdef\SetFigFont#1#2#3#4#5{%
  \reset@font\fontsize{#1}{#2pt}%
  \fontfamily{#3}\fontseries{#4}\fontshape{#5}%
  \selectfont}%
\fi\endgroup%
\begin{picture}(5574,2928)(1489,-3055)
\put(4126,-286){\makebox(0,0)[lb]{\smash{{\SetFigFont{10}{12.0}{\rmdefault}{\mddefault}{\updefault}{\color[rgb]{0,0,0}\small $N$}%
}}}}
\put(4126,-2986){\makebox(0,0)[lb]{\smash{{\SetFigFont{10}{12.0}{\rmdefault}{\mddefault}{\updefault}{\color[rgb]{0,0,0}\small $S$}%
}}}}
\put(3526,-511){\makebox(0,0)[lb]{\smash{{\SetFigFont{10}{12.0}{\rmdefault}{\mddefault}{\updefault}{\color[rgb]{0,0,0}\small $\Delta_-$}%
}}}}
\put(4726,-511){\makebox(0,0)[lb]{\smash{{\SetFigFont{10}{12.0}{\rmdefault}{\mddefault}{\updefault}{\color[rgb]{0,0,0}\small $\Delta_+$}%
}}}}
\put(2401,-2986){\makebox(0,0)[lb]{\smash{{\SetFigFont{10}{12.0}{\rmdefault}{\mddefault}{\updefault}{\color[rgb]{0,0,0}\small $O'_-$}%
}}}}
\put(6001,-2986){\makebox(0,0)[lb]{\smash{{\SetFigFont{10}{12.0}{\rmdefault}{\mddefault}{\updefault}{\color[rgb]{0,0,0}\small $O'_+$}%
}}}}
\put(6518,-1478){\makebox(0,0)[lb]{\smash{{\SetFigFont{10}{12.0}{\rmdefault}{\mddefault}{\updefault}{\color[rgb]{0,0,0}\small $\Gamma$}%
}}}}
\put(1576,-886){\makebox(0,0)[lb]{\smash{{\SetFigFont{10}{12.0}{\rmdefault}{\mddefault}{\updefault}{\color[rgb]{0,0,0}\small $B_{d+1}(y)$}%
}}}}
\put(1801,-1186){\makebox(0,0)[lb]{\smash{{\SetFigFont{10}{12.0}{\rmdefault}{\mddefault}{\updefault}{\color[rgb]{0,0,0}\small $B_d(y)$}%
}}}}
\end{picture}%

%% file: fig-theo-verticales.pstex_t
\begin{picture}(0,0)%
\includegraphics{fig-theo-verticales.pstex}%
\end{picture}%
\setlength{\unitlength}{3947sp}%
\begingroup\makeatletter\ifx\SetFigFont\undefined%
\gdef\SetFigFont#1#2#3#4#5{%
  \reset@font\fontsize{#1}{#2pt}%
  \fontfamily{#3}\fontseries{#4}\fontshape{#5}%
  \selectfont}%
\fi\endgroup%
\begin{picture}(3905,2508)(1636,-2698)
\put(1651,-361){\makebox(0,0)[lb]{\smash{{\SetFigFont{11}{13.2}{\rmdefault}{\mddefault}{\updefault}{\color[rgb]{0,0,0}\small $M((x_k))=$}%
}}}}
\put(2926,-361){\makebox(0,0)[lb]{\smash{{\SetFigFont{11}{13.2}{\rmdefault}{\mddefault}{\updefault}{\color[rgb]{0,0,0}\small $\downarrow$}%
}}}}
\put(3301,-361){\makebox(0,0)[lb]{\smash{{\SetFigFont{11}{13.2}{\rmdefault}{\mddefault}{\updefault}{\color[rgb]{0,0,0}\small $\rightarrow$}%
}}}}
\put(3751,-361){\makebox(0,0)[lb]{\smash{{\SetFigFont{11}{13.2}{\rmdefault}{\mddefault}{\updefault}{\color[rgb]{0,0,0}\small $\uparrow$}%
}}}}
\put(4126,-361){\makebox(0,0)[lb]{\smash{{\SetFigFont{11}{13.2}{\rmdefault}{\mddefault}{\updefault}{\color[rgb]{0,0,0}\small $\rightarrow$}%
}}}}
\put(4501,-361){\makebox(0,0)[lb]{\smash{{\SetFigFont{11}{13.2}{\rmdefault}{\mddefault}{\updefault}{\color[rgb]{0,0,0}\small $\uparrow$}%
}}}}
\put(4801,-361){\makebox(0,0)[lb]{\smash{{\SetFigFont{11}{13.2}{\rmdefault}{\mddefault}{\updefault}{\color[rgb]{0,0,0}\small $\cdots$}%
}}}}
\put(5526,-1651){\makebox(0,0)[lb]{\smash{{\SetFigFont{12}{14.4}{\rmdefault}{\mddefault}{\updefault}{\color[rgb]{0,0,0}\small $\Gamma$}%
}}}}
\put(3321,-1806){\makebox(0,0)[lb]{\smash{{\SetFigFont{11}{13.2}{\rmdefault}{\mddefault}{\updefault}{\color[rgb]{0,0,0}\small $x_k$}%
}}}}
\put(4211,-1806){\makebox(0,0)[lb]{\smash{{\SetFigFont{11}{13.2}{\rmdefault}{\mddefault}{\updefault}{\color[rgb]{0,0,0}\small $x_{k+1}$}%
}}}}
\put(5096,-956){\makebox(0,0)[lb]{\smash{{\SetFigFont{12}{14.4}{\rmdefault}{\mddefault}{\updefault}{\color[rgb]{0,0,0}\small $O_N(\Gamma)$}%
}}}}
\put(5101,-2311){\makebox(0,0)[lb]{\smash{{\SetFigFont{12}{14.4}{\rmdefault}{\mddefault}{\updefault}{\color[rgb]{0,0,0}\small $O_S(\Gamma)$}%
}}}}
\put(2551,-361){\makebox(0,0)[lb]{\smash{{\SetFigFont{11}{13.2}{\rmdefault}{\mddefault}{\updefault}{\color[rgb]{0,0,0}\small $\cdots$}%
}}}}
\end{picture}%

%% file: fig-image-gamma0-indif2.pstex_t
\begin{picture}(0,0)%
\includegraphics{fig-image-gamma0-indif2.pstex}%
\end{picture}%
\setlength{\unitlength}{3158sp}%
\begingroup\makeatletter\ifx\SetFigFont\undefined%
\gdef\SetFigFont#1#2#3#4#5{%
  \reset@font\fontsize{#1}{#2pt}%
  \fontfamily{#3}\fontseries{#4}\fontshape{#5}%
  \selectfont}%
\fi\endgroup%
\begin{picture}(3627,3025)(2089,-3439)
\put(3976,-1486){\makebox(0,0)[lb]{\smash{{\SetFigFont{10}{12.0}{\rmdefault}{\mddefault}{\updefault}{\color[rgb]{0,0,0}\small $d$}%
}}}}
\put(2916,-1816){\makebox(0,0)[lb]{\smash{{\SetFigFont{10}{12.0}{\rmdefault}{\mddefault}{\updefault}{\color[rgb]{0,0,0}\small $x_0$}%
}}}}
\put(4491,-1801){\makebox(0,0)[lb]{\smash{{\SetFigFont{10}{12.0}{\rmdefault}{\mddefault}{\updefault}{\color[rgb]{0,0,0}\small $x_1$}%
}}}}
\put(3396,-1466){\makebox(0,0)[lb]{\smash{{\SetFigFont{10}{12.0}{\rmdefault}{\mddefault}{\updefault}{\color[rgb]{0,0,0}\small $p$}%
}}}}
\put(2881,-2566){\makebox(0,0)[lb]{\smash{{\SetFigFont{10}{12.0}{\rmdefault}{\mddefault}{\updefault}{\color[rgb]{0,0,0}\small $h^n(x_0)$}%
}}}}
\put(4566,-1061){\makebox(0,0)[lb]{\smash{{\SetFigFont{10}{12.0}{\rmdefault}{\mddefault}{\updefault}{\color[rgb]{0,0,0}\small $h^n(x_1)$}%
}}}}
\put(3796,-606){\makebox(0,0)[lb]{\smash{{\SetFigFont{10}{12.0}{\rmdefault}{\mddefault}{\updefault}{\color[rgb]{0,0,0}\small $h^n(d)$}%
}}}}
\put(5486,-1651){\makebox(0,0)[lb]{\smash{{\SetFigFont{10}{12.0}{\rmdefault}{\mddefault}{\updefault}{\color[rgb]{0,0,0}\small $\Gamma$}%
}}}}
\put(5701,-961){\makebox(0,0)[lb]{\smash{{\SetFigFont{10}{12.0}{\rmdefault}{\mddefault}{\updefault}{\color[rgb]{0,0,0}\small $O_N(\Gamma)$}%
}}}}
\put(5701,-2311){\makebox(0,0)[lb]{\smash{{\SetFigFont{10}{12.0}{\rmdefault}{\mddefault}{\updefault}{\color[rgb]{0,0,0}\small $O_S(\Gamma)$}%
}}}}
\put(3636,-2796){\makebox(0,0)[lb]{\smash{{\SetFigFont{10}{12.0}{\rmdefault}{\mddefault}{\updefault}{\color[rgb]{0,0,0}\small $G$}%
}}}}
\put(4351,-3361){\makebox(0,0)[lb]{\smash{{\SetFigFont{10}{12.0}{\rmdefault}{\mddefault}{\updefault}{\color[rgb]{0,0,0}\small $O_+(G)$}%
}}}}
\put(3226,-2161){\makebox(0,0)[lb]{\smash{{\SetFigFont{10}{12.0}{\rmdefault}{\mddefault}{\updefault}{\color[rgb]{0,0,0}\small $h^n(p)$}%
}}}}
\put(2776,-3286){\makebox(0,0)[lb]{\smash{{\SetFigFont{10}{12.0}{\rmdefault}{\mddefault}{\updefault}{\color[rgb]{0,0,0}\small $O_-(G)$}%
}}}}
\end{picture}%

%% file: fig-gamma-prime.pstex_t
\begin{picture}(0,0)%
\includegraphics{fig-gamma-prime.pstex}%
\end{picture}%
\setlength{\unitlength}{3158sp}%
\begingroup\makeatletter\ifx\SetFigFont\undefined%
\gdef\SetFigFont#1#2#3#4#5{%
  \reset@font\fontsize{#1}{#2pt}%
  \fontfamily{#3}\fontseries{#4}\fontshape{#5}%
  \selectfont}%
\fi\endgroup%
\begin{picture}(5130,2275)(1036,-2689)
\put(3396,-1466){\makebox(0,0)[lb]{\smash{{\SetFigFont{10}{12.0}{\rmdefault}{\mddefault}{\updefault}{\color[rgb]{0,0,0}\small $p$}%
}}}}
\put(3796,-606){\makebox(0,0)[lb]{\smash{{\SetFigFont{10}{12.0}{\rmdefault}{\mddefault}{\updefault}{\color[rgb]{0,0,0}\small $h(d)$}%
}}}}
\put(5701,-2611){\makebox(0,0)[lb]{\smash{{\SetFigFont{10}{12.0}{\rmdefault}{\mddefault}{\updefault}{\color[rgb]{0,0,0}\small $O_S(\Gamma)$}%
}}}}
\put(5701,-736){\makebox(0,0)[lb]{\smash{{\SetFigFont{10}{12.0}{\rmdefault}{\mddefault}{\updefault}{\color[rgb]{0,0,0}\small $O_N(\Gamma)$}%
}}}}
\put(3976,-1486){\makebox(0,0)[lb]{\smash{{\SetFigFont{10}{12.0}{\rmdefault}{\mddefault}{\updefault}{\color[rgb]{0,0,0}\small $d=A$}%
}}}}
\put(6151,-1636){\makebox(0,0)[lb]{\smash{{\SetFigFont{10}{12.0}{\rmdefault}{\mddefault}{\updefault}{\color[rgb]{0,0,0}\small $\Gamma^+$}%
}}}}
\put(3796,-2521){\makebox(0,0)[lb]{\smash{{\SetFigFont{10}{12.0}{\rmdefault}{\mddefault}{\updefault}{\color[rgb]{0,0,0}\small $h(p)=B$}%
}}}}
\put(6151,-2011){\makebox(0,0)[lb]{\smash{{\SetFigFont{10}{12.0}{\rmdefault}{\mddefault}{\updefault}{\color[rgb]{0,0,0}\small ${\Gamma'}^+$}%
}}}}
\put(1051,-1111){\makebox(0,0)[lb]{\smash{{\SetFigFont{10}{12.0}{\rmdefault}{\mddefault}{\updefault}{\color[rgb]{0,0,0}\small ${\Gamma'}^-$}%
}}}}
\put(1051,-1636){\makebox(0,0)[lb]{\smash{{\SetFigFont{10}{12.0}{\rmdefault}{\mddefault}{\updefault}{\color[rgb]{0,0,0}\small $\Gamma^-$}%
}}}}
\end{picture}%

%% file: fig-B-a-droite.pstex_t
\begin{picture}(0,0)%
\includegraphics{fig-B-a-droite.pstex}%
\end{picture}%
\setlength{\unitlength}{3158sp}%
\begingroup\makeatletter\ifx\SetFigFont\undefined%
\gdef\SetFigFont#1#2#3#4#5{%
  \reset@font\fontsize{#1}{#2pt}%
  \fontfamily{#3}\fontseries{#4}\fontshape{#5}%
  \selectfont}%
\fi\endgroup%
\begin{picture}(4677,2270)(1489,-3130)
\put(2401,-1411){\makebox(0,0)[lb]{\smash{{\SetFigFont{10}{12.0}{\rmdefault}{\mddefault}{\updefault}{\color[rgb]{0,0,0}\small $C$}%
}}}}
\put(5176,-1411){\makebox(0,0)[lb]{\smash{{\SetFigFont{10}{12.0}{\rmdefault}{\mddefault}{\updefault}{\color[rgb]{0,0,0}\small $D$}%
}}}}
\put(6151,-1636){\makebox(0,0)[lb]{\smash{{\SetFigFont{10}{12.0}{\rmdefault}{\mddefault}{\updefault}{\color[rgb]{0,0,0}\small $\Gamma^+$}%
}}}}
\put(6151,-2011){\makebox(0,0)[lb]{\smash{{\SetFigFont{10}{12.0}{\rmdefault}{\mddefault}{\updefault}{\color[rgb]{0,0,0}\small ${\Gamma'}^+$}%
}}}}
\put(3826,-3061){\makebox(0,0)[lb]{\smash{{\SetFigFont{10}{12.0}{\rmdefault}{\mddefault}{\updefault}{\color[rgb]{0,0,0}\small $B$}%
}}}}
\put(3826,-1411){\makebox(0,0)[lb]{\smash{{\SetFigFont{10}{12.0}{\rmdefault}{\mddefault}{\updefault}{\color[rgb]{0,0,0}\small $A$}%
}}}}
\put(4051,-1861){\makebox(0,0)[lb]{\smash{{\SetFigFont{10}{12.0}{\rmdefault}{\mddefault}{\updefault}{\color[rgb]{0,0,0}\small $A'$}%
}}}}
\put(4051,-2611){\makebox(0,0)[lb]{\smash{{\SetFigFont{10}{12.0}{\rmdefault}{\mddefault}{\updefault}{\color[rgb]{0,0,0}\small $B'$}%
}}}}
\put(1651,-2161){\makebox(0,0)[lb]{\smash{{\SetFigFont{10}{12.0}{\rmdefault}{\mddefault}{\updefault}{\color[rgb]{0,0,0}\small $U$}%
}}}}
\end{picture}%

%% file: fig-gammaprime.pstex_t
\begin{picture}(0,0)%
\includegraphics{fig-gammaprime.pstex}%
\end{picture}%
\setlength{\unitlength}{3158sp}%
\begingroup\makeatletter\ifx\SetFigFont\undefined%
\gdef\SetFigFont#1#2#3#4#5{%
  \reset@font\fontsize{#1}{#2pt}%
  \fontfamily{#3}\fontseries{#4}\fontshape{#5}%
  \selectfont}%
\fi\endgroup%
\begin{picture}(4912,2451)(589,-2865)
\put(1801,-1786){\makebox(0,0)[lb]{\smash{{\SetFigFont{10}{12.0}{\rmdefault}{\mddefault}{\updefault}{\color[rgb]{0,0,0}\small $x_{-1}$}%
}}}}
\put(5486,-1651){\makebox(0,0)[lb]{\smash{{\SetFigFont{10}{12.0}{\rmdefault}{\mddefault}{\updefault}{\color[rgb]{0,0,0}\small $\Gamma$}%
}}}}
\put(3636,-2796){\makebox(0,0)[lb]{\smash{{\SetFigFont{10}{12.0}{\rmdefault}{\mddefault}{\updefault}{\color[rgb]{0,0,0}\small $G$}%
}}}}
\put(2851,-1861){\makebox(0,0)[lb]{\smash{{\SetFigFont{10}{12.0}{\rmdefault}{\mddefault}{\updefault}{\color[rgb]{0,0,0}\small $x_0$}%
}}}}
\put(2513,-2483){\makebox(0,0)[lb]{\smash{{\SetFigFont{10}{12.0}{\rmdefault}{\mddefault}{\updefault}{\color[rgb]{0,0,0}\small $\gamma'$}%
}}}}
\put(4801,-1111){\makebox(0,0)[lb]{\smash{{\SetFigFont{10}{12.0}{\rmdefault}{\mddefault}{\updefault}{\color[rgb]{0,0,0}\small $D$}%
}}}}
\put(4801,-2461){\makebox(0,0)[lb]{\smash{{\SetFigFont{10}{12.0}{\rmdefault}{\mddefault}{\updefault}{\color[rgb]{0,0,0}\small $\gamma'_0$}%
}}}}
\put(2851,-811){\makebox(0,0)[lb]{\smash{{\SetFigFont{10}{12.0}{\rmdefault}{\mddefault}{\updefault}{\color[rgb]{0,0,0}\small $V$}%
}}}}
\put(4576,-811){\makebox(0,0)[lb]{\smash{{\SetFigFont{10}{12.0}{\rmdefault}{\mddefault}{\updefault}{\color[rgb]{0,0,0}\small $d$}%
}}}}
\put(1051,-1786){\makebox(0,0)[lb]{\smash{{\SetFigFont{10}{12.0}{\rmdefault}{\mddefault}{\updefault}{\color[rgb]{0,0,0}\small $x'_{-1}$}%
}}}}
\put(4576,-1786){\makebox(0,0)[lb]{\smash{{\SetFigFont{10}{12.0}{\rmdefault}{\mddefault}{\updefault}{\color[rgb]{0,0,0}\small $x_1$}%
}}}}
\put(3376,-1936){\makebox(0,0)[lb]{\smash{{\SetFigFont{10}{12.0}{\rmdefault}{\mddefault}{\updefault}{\color[rgb]{0,0,0}\small $x'_0$}%
}}}}
\end{picture}%

%% file: fig-gammaprime2.pstex_t
\begin{picture}(0,0)%
\includegraphics{fig-gammaprime2.pstex}%
\end{picture}%
\setlength{\unitlength}{3158sp}%
\begingroup\makeatletter\ifx\SetFigFont\undefined%
\gdef\SetFigFont#1#2#3#4#5{%
  \reset@font\fontsize{#1}{#2pt}%
  \fontfamily{#3}\fontseries{#4}\fontshape{#5}%
  \selectfont}%
\fi\endgroup%
\begin{picture}(4912,1712)(589,-2865)
\put(3636,-2796){\makebox(0,0)[lb]{\smash{{\SetFigFont{10}{12.0}{\rmdefault}{\mddefault}{\updefault}{\color[rgb]{0,0,0}\small $G$}%
}}}}
\put(5486,-1651){\makebox(0,0)[lb]{\smash{{\SetFigFont{10}{12.0}{\rmdefault}{\mddefault}{\updefault}{\color[rgb]{0,0,0}\small $\Gamma'$}%
}}}}
\put(2513,-2483){\makebox(0,0)[lb]{\smash{{\SetFigFont{10}{12.0}{\rmdefault}{\mddefault}{\updefault}{\color[rgb]{0,0,0}\small $\gamma'_{-1}$}%
}}}}
\put(751,-1786){\makebox(0,0)[lb]{\smash{{\SetFigFont{10}{12.0}{\rmdefault}{\mddefault}{\updefault}{\color[rgb]{0,0,0}\small $\gamma'_{-2}$}%
}}}}
\put(1126,-1411){\makebox(0,0)[lb]{\smash{{\SetFigFont{10}{12.0}{\rmdefault}{\mddefault}{\updefault}{\color[rgb]{0,0,0}\small $x'_{-1}$}%
}}}}
\put(3376,-1786){\makebox(0,0)[lb]{\smash{{\SetFigFont{10}{12.0}{\rmdefault}{\mddefault}{\updefault}{\color[rgb]{0,0,0}\small $x'_0$}%
}}}}
\put(4126,-1936){\makebox(0,0)[lb]{\smash{{\SetFigFont{10}{12.0}{\rmdefault}{\mddefault}{\updefault}{\color[rgb]{0,0,0}\small $\gamma'_0$}%
}}}}
\put(4426,-1336){\makebox(0,0)[lb]{\smash{{\SetFigFont{10}{12.0}{\rmdefault}{\mddefault}{\updefault}{\color[rgb]{0,0,0}\small $x'_1=x_1$}%
}}}}
\end{picture}%

%% file: fig-quadrillage.pstex_t
\begin{picture}(0,0)%
\includegraphics{fig-quadrillage.pstex}%
\end{picture}%
\setlength{\unitlength}{3158sp}%
\begingroup\makeatletter\ifx\SetFigFont\undefined%
\gdef\SetFigFont#1#2#3#4#5{%
  \reset@font\fontsize{#1}{#2pt}%
  \fontfamily{#3}\fontseries{#4}\fontshape{#5}%
  \selectfont}%
\fi\endgroup%
\begin{picture}(3627,2424)(1189,-2773)
\put(4801,-1336){\makebox(0,0)[lb]{\smash{{\SetFigFont{10}{12.0}{\rmdefault}{\mddefault}{\updefault}{\color[rgb]{0,0,0}\small $D(x)$}%
}}}}
\put(4801,-2311){\makebox(0,0)[lb]{\smash{{\SetFigFont{10}{12.0}{\rmdefault}{\mddefault}{\updefault}{\color[rgb]{0,0,0}\small $C$}%
}}}}
\put(4801,-1861){\makebox(0,0)[lb]{\smash{{\SetFigFont{10}{12.0}{\rmdefault}{\mddefault}{\updefault}{\color[rgb]{0,0,0}\small $x$}%
}}}}
\end{picture}%

%% file: fig-chaine2.pstex_t
\begin{picture}(0,0)%
\includegraphics{fig-chaine2.pstex}%
\end{picture}%
\setlength{\unitlength}{3158sp}%
\begingroup\makeatletter\ifx\SetFigFont\undefined%
\gdef\SetFigFont#1#2#3#4#5{%
  \reset@font\fontsize{#1}{#2pt}%
  \fontfamily{#3}\fontseries{#4}\fontshape{#5}%
  \selectfont}%
\fi\endgroup%
\begin{picture}(5523,2420)(-14,-2704)
\put(1191,-2451){\makebox(0,0)[lb]{\smash{{\SetFigFont{10}{12.0}{\rmdefault}{\mddefault}{\updefault}{\color[rgb]{0,0,0}\small $h(C)$}%
}}}}
\put(4768,-651){\makebox(0,0)[lb]{\smash{{\SetFigFont{10}{12.0}{\rmdefault}{\mddefault}{\updefault}{\color[rgb]{0,0,0}\small $h^n(C)$}%
}}}}
\put(2993,-443){\makebox(0,0)[lb]{\smash{{\SetFigFont{10}{12.0}{\rmdefault}{\mddefault}{\updefault}{\color[rgb]{0,0,0}\small $E$}%
}}}}
\put(4590,-1839){\makebox(0,0)[lb]{\smash{{\SetFigFont{10}{12.0}{\rmdefault}{\mddefault}{\updefault}{\color[rgb]{0,0,0}\small $h^n(E)$}%
}}}}
\put(3286,-2626){\makebox(0,0)[lb]{\smash{{\SetFigFont{10}{12.0}{\rmdefault}{\mddefault}{\updefault}{\color[rgb]{0,0,0}\small $h(E)$}%
}}}}
\put(3091,-1515){\makebox(0,0)[lb]{\smash{{\SetFigFont{10}{12.0}{\rmdefault}{\mddefault}{\updefault}{\color[rgb]{0,0,0}\small $D_1$}%
}}}}
\put(1276,-661){\makebox(0,0)[lb]{\smash{{\SetFigFont{10}{12.0}{\rmdefault}{\mddefault}{\updefault}{\color[rgb]{0,0,0}\small $C$}%
}}}}
\put(  1,-1411){\makebox(0,0)[lb]{\smash{{\SetFigFont{10}{12.0}{\rmdefault}{\mddefault}{\updefault}{\color[rgb]{0,0,0}\small $h^{-n}(D_1)$}%
}}}}
\end{picture}%

%% file: fig-ex-compo.pstex_t
\begin{picture}(0,0)%
\includegraphics{fig-ex-compo.pstex}%
\end{picture}%
\setlength{\unitlength}{3552sp}%
\begingroup\makeatletter\ifx\SetFigFont\undefined%
\gdef\SetFigFont#1#2#3#4#5{%
  \reset@font\fontsize{#1}{#2pt}%
  \fontfamily{#3}\fontseries{#4}\fontshape{#5}%
  \selectfont}%
\fi\endgroup%
\begin{picture}(6249,2424)(964,-2773)
\put(5251,-2686){\makebox(0,0)[lb]{\smash{{\SetFigFont{14}{16.8}{\rmdefault}{\mddefault}{\updefault}{\color[rgb]{0,0,0}\small $F$}%
}}}}
\put(6301,-2686){\makebox(0,0)[lb]{\smash{{\SetFigFont{14}{16.8}{\rmdefault}{\mddefault}{\updefault}{\color[rgb]{0,0,0}\small $G$}%
}}}}
\put(4651,-1636){\makebox(0,0)[lb]{\smash{{\SetFigFont{14}{16.8}{\rmdefault}{\mddefault}{\updefault}{\color[rgb]{0,0,0}\small $x$}%
}}}}
\put(6976,-1636){\makebox(0,0)[lb]{\smash{{\SetFigFont{14}{16.8}{\rmdefault}{\mddefault}{\updefault}{\color[rgb]{0,0,0}\small $y$}%
}}}}
\put(1013,-1846){\makebox(0,0)[lb]{\smash{{\SetFigFont{12}{14.4}{\rmdefault}{\mddefault}{\updefault}{\color[rgb]{0,0,0}\small $F$}%
}}}}
\put(1276,-1636){\makebox(0,0)[lb]{\smash{{\SetFigFont{12}{14.4}{\rmdefault}{\mddefault}{\updefault}{\color[rgb]{0,0,0}\small $x$}%
}}}}
\put(3676,-1636){\makebox(0,0)[lb]{\smash{{\SetFigFont{12}{14.4}{\rmdefault}{\mddefault}{\updefault}{\color[rgb]{0,0,0}\small $y$}%
}}}}
\put(2926,-2686){\makebox(0,0)[lb]{\smash{{\SetFigFont{12}{14.4}{\rmdefault}{\mddefault}{\updefault}{\color[rgb]{0,0,0}\small $G$}%
}}}}
\end{picture}%

%% file: fig-flot-perturbe2.pstex_t
\begin{picture}(0,0)%
\includegraphics[scale=0.9]{fig-flot-perturbe2.pstex}%
\end{picture}%
\setlength{\unitlength}{1421sp}%
\begingroup\makeatletter\ifx\SetFigFont\undefined%
\gdef\SetFigFont#1#2#3#4#5{%
  \reset@font\fontsize{#1}{#2pt}%
  \fontfamily{#3}\fontseries{#4}\fontshape{#5}%
  \selectfont}%
\fi\endgroup%
\begin{picture}(15419,7729)(1029,-7208)
\put(9976,-7111){\makebox(0,0)[lb]{\smash{{\SetFigFont{10}{24.0}{\rmdefault}{\mddefault}{\updefault}{\color[rgb]{0,0,0}\small $F$}%
}}}}
\put(2701,-6886){\makebox(0,0)[lb]{\smash{{\SetFigFont{10}{24.0}{\rmdefault}{\mddefault}{\updefault}{\color[rgb]{0,0,0}\small $\phi:$}%
}}}}
\put(13951,-6061){\makebox(0,0)[lb]{\smash{{\SetFigFont{10}{24.0}{\rmdefault}{\mddefault}{\updefault}{\color[rgb]{0,0,0}\small $G'$}%
}}}}
\put(11176,-5536){\makebox(0,0)[lb]{\smash{{\SetFigFont{10}{24.0}{\rmdefault}{\mddefault}{\updefault}{\color[rgb]{0,0,0}\small $x_2$}%
}}}}
\put(8641,-3751){\makebox(0,0)[lb]{\smash{{\SetFigFont{10}{24.0}{\rmdefault}{\mddefault}{\updefault}{\color[rgb]{0,0,0}\small $x_1$}%
}}}}
\put(10501,-361){\makebox(0,0)[lb]{\smash{{\SetFigFont{10}{24.0}{\rmdefault}{\mddefault}{\updefault}{\color[rgb]{0,0,0}\small $F'$}%
}}}}
\put(15901,-3511){\makebox(0,0)[lb]{\smash{{\SetFigFont{10}{24.0}{\rmdefault}{\mddefault}{\updefault}{\color[rgb]{0,0,0}\small $y_1$}%
}}}}
\put(12076,389){\makebox(0,0)[lb]{\smash{{\SetFigFont{10}{24.0}{\rmdefault}{\mddefault}{\updefault}{\color[rgb]{0,0,0}\small $\phi \circ h$}%
}}}}
\put(13501,-961){\makebox(0,0)[lb]{\smash{{\SetFigFont{10}{24.0}{\rmdefault}{\mddefault}{\updefault}{\color[rgb]{0,0,0}\small $y_2$}%
}}}}
\put(7876,-3721){\makebox(0,0)[lb]{\smash{{\SetFigFont{10}{24.0}{\rmdefault}{\mddefault}{\updefault}{\color[rgb]{0,0,0}\small $y_1$}%
}}}}
\put(4426,314){\makebox(0,0)[lb]{\smash{{\SetFigFont{10}{24.0}{\rmdefault}{\mddefault}{\updefault}{\color[rgb]{0,0,0}\small $h$}%
}}}}
\put(7771,-2146){\makebox(0,0)[lb]{\smash{{\SetFigFont{10}{24.0}{\rmdefault}{\mddefault}{\updefault}{\color[rgb]{0,0,0}\small $D$}%
}}}}
\put(4651,-6436){\makebox(0,0)[lb]{\smash{{\SetFigFont{10}{24.0}{\rmdefault}{\mddefault}{\updefault}{\color[rgb]{0,0,0}\small $\Delta$}%
}}}}
\put(1426,-6961){\makebox(0,0)[lb]{\smash{{\SetFigFont{10}{24.0}{\rmdefault}{\mddefault}{\updefault}{\color[rgb]{0,0,0}\small $F$}%
}}}}
\put(7351,-6961){\makebox(0,0)[lb]{\smash{{\SetFigFont{10}{24.0}{\rmdefault}{\mddefault}{\updefault}{\color[rgb]{0,0,0}\small $G$}%
}}}}
\put(1141,-3736){\makebox(0,0)[lb]{\smash{{\SetFigFont{10}{24.0}{\rmdefault}{\mddefault}{\updefault}{\color[rgb]{0,0,0}\small $x_1$}%
}}}}
\put(14626,-7111){\makebox(0,0)[lb]{\smash{{\SetFigFont{10}{24.0}{\rmdefault}{\mddefault}{\updefault}{\color[rgb]{0,0,0}\small $G$}%
}}}}
\end{picture}%